\documentclass[11pt]{article}

\usepackage[T1]{fontenc}        
\usepackage[utf8]{inputenc}     
\usepackage[hidelinks, colorlinks=true, linkcolor=refcolor, citecolor=bibcolor]{hyperref}
\usepackage{lmodern}
\usepackage{amsmath, amssymb, amsthm, amsfonts, bbm}
\usepackage{amsbsy,stmaryrd}
\usepackage{mathtools}
\usepackage{color}
\usepackage{float}
\usepackage{tikz}
\usetikzlibrary{patterns}
\usepackage{thmtools}
\usepackage{mathtools}
\usepackage{mathrsfs}              
\usepackage{bbold}              
\usepackage{biblatex}
\usepackage{xcolor}
\usepackage{enumitem}
\usepackage{pgfplots}
\usepackage{pgfplotstable}
\usepackage{subcaption}
\usepackage{tikz-3dplot}        
\usetikzlibrary{hobby, patterns, fadings}
\usepgfplotslibrary{groupplots, units}

\usetikzlibrary{calc}
\usetikzlibrary{math}
\usetikzlibrary{fit, positioning}
\usetikzlibrary{decorations.pathreplacing,calligraphy}
\usetikzlibrary{arrows.meta}

\declaretheoremstyle[
    spaceabove=6pt, spacebelow=6pt,
    headfont=\normalfont\bfseries, %
    notefont=\normalfont\bfseries,
    notebraces={}{},
    bodyfont=\normalfont,
    postheadspace=1em, qed=$\blacksquare$
]{demstyle}

\declaretheorem[name=Theorem,numberwithin=section]{thm}
\declaretheorem[name=Lemma,sibling=thm]{lem}
\declaretheorem[name=Proposition,sibling=thm]{prop}
\declaretheorem[name=Corollary,sibling=thm]{cor}

\declaretheorem[name=Definition, style=definition, sibling=thm]{defi}

\declaretheorem[name=Remark,sibling=thm]{rmk}

\declaretheorem[name=Proof,numbered=no,style=demstyle]{dem}

\numberwithin{equation}{section}
\pgfplotsset{width=0.98\textwidth,
        compat=1.9 }
\definecolor{bleuDoux}{HTML}{3333ff}
\definecolor{orangeTerre}{HTML}{e76d55}
\definecolor{rougeTerre}{RGB}{235, 86, 75}
\definecolor{surligneur}{RGB}{0, 0, 0}

\pgfplotsset{
 colormap={parula}{
  rgb255=(53,42,135)
  rgb255=(15,92,221)
  rgb255=(18,125,216)
  rgb255=(7,156,207)
	rgb255=(21,177,180)
	rgb255=(89,189,140)
	rgb255=(165,190,107)
	rgb255=(225,185,82)
	rgb255=(252,206,46)
	rgb255=(249,251,14)
 },
}
\definecolor{refcolor}{RGB}{255, 0, 143}
\definecolor{bibcolor}{RGB}{0, 0, 150}
\pgfdeclarefading{degrade bas vers haut}
 {\tikz
    \fill[bottom color=transparent!0, top color=transparent!100] (0,0) rectangle (50bp,50bp);}

\newcommand{\nc}{\newcommand}

\nc{\ret}{\\[12pt]}
\nc{\rett}{\\[18pt]}

\newcommand{\eps}{\varepsilon}

\newcommand{\Real}{\mathop{\mathfrak{Re}}\nolimits}
\newcommand{\Imag}{\mathop{\mathfrak{Im}}\nolimits}

\DeclareMathOperator{\Ker}{Ker}
\DeclareMathOperator{\vect}{Span}

\newcommand*{\supp}{\operatorname{supp}}

\newcommand*{\euler}{\mathrm{e}}

\newcommand{\spforall}{\forall\;}
\newcommand{\aeforall}{\operatorname{a.e.}\;}
\newcommand{\spexists}{\exists\;}

\newcommand\restr[2]{{
  \kern-\nulldelimiterspace 
  #1 
  |_{#2} 
  }}

\nc{\R}{\mathbb{R}}
\nc{\N}{\mathbb{N}}
\nc{\C}{\mathbb{C}}
\nc{\Z}{\mathbb{Z}}
\nc{\Q}{\mathbb{Q}}

\newcommand{\esplignestableau}{\\[3pt]}

\newcommand{\itbf}[1]{\textit{\textbf{#1}}}
\newcommand{\xv}{x}       
\newcommand{\yv}{{\textit{\textbf{y}}}}       
\newcommand{\yi}{y}       
\newcommand{\zv}{{\textit{\textbf{z}}}}       
\newcommand{\zi}{z}       
\newcommand{\sv}{{\textit{\textbf{s}}}}       
\newcommand{\si}{s}       
\newcommand{\kv}{{\textit{\textbf{k}}}}       
\newcommand{\ki}{k}       
\newcommand{\ev}{{\textit{\textbf{e}}}}

\newcommand{\hv}{{\underline{\textit{\textbf{h}}}}}
\newcommand{\ord}{n}      
\newcommand{\cut}{{\boldsymbol{\theta}}} 
\newcommand{\cutslope}{{\boldsymbol{\delta}}}
\newcommand{\cuti}{\theta} 
\newcommand*{\Dt}[1]{D^{#1\vphantom{+}}_\cut\,} 
\newcommand*{\diese}{{\smash{\scriptstyle \sharp}}}
\newcommand{\icplx}{\ensuremath{\mathrm{i}}}

\graphicspath{{figures/}}
\pgfplotsset{table/search path={figures},}
\input{./figures/figures.tex}

\newboolean{afficherGraphes}
\setboolean{afficherGraphes}{true}
\title{Wave propagation in one-dimensional quasiperiodic media}
\author{Pierre Amenoagbadji, Sonia Fliss, Patrick Joly}
\date{}

\begin{document}
%
%
%
\maketitle

\begin{abstract}
	This work is devoted to the resolution of the Helmholtz equation $-(\mu\, u')' - \rho\, \omega^2 u = f$ in a one-dimensional unbounded medium. We assume the coefficients of this equation to be local perturbations of \textit{quasiperiodic} functions, namely the traces along a particular line of higher-dimensional periodic functions. Using the definition of quasiperiodicity, the problem is lifted onto a higher-dimensional problem with periodic coefficients. The periodicity of the augmented problem allows us to extend the ideas of the DtN-based method developed in \cite{flissthese, jolyLiFliss} for the elliptic case. However, the associated mathematical and numerical analysis of the method are more delicate because the augmented PDE is degenerate, in the sense that the principal part of its operator is no longer elliptic. We also study the numerical resolution of this PDE, which relies on the resolution of Dirichlet cell problems as well as a constrained Riccati equation.
\end{abstract}

\section{Introduction and motivation}\label{sec:introduction_motivation}
We consider the Helmholtz equation
\begin{equation}
\label{eq:whole_line_problem}
\displaystyle
- \frac{d}{d \xv} \Big( \mu \; \frac{d u}{d \xv} \Big) - \rho \; \omega^2 \; u = f \quad \textnormal{in} \quad \R,
\end{equation}
where the coefficients $\mu$ and $\rho$ have positive upper and lower bounds:
	\begin{equation}\label{eq:coef_ellipt}
		\displaystyle
		\spexists \mu_\pm, \rho_\pm, \quad \spforall x \in \R, \qquad 0 < \mu_- \leq \mu(\xv) \leq \mu_+ \quad \textnormal{and} \quad 0 < \rho_- \leq \rho(\xv) \leq \rho_+.
	\end{equation}
 The source term $f$ belongs to $L^2(\R)$ and is assumed to have a compact support:
\begin{equation}
	\label{eq:source_terme}
	\spexists a>0,\quad \supp f \subset (-a, a).
\end{equation}

	\noindent

	\noindent
	Equation \eqref{eq:whole_line_problem} is encountered when one is looking for time-harmonic solutions $u(x)\, e^{\icplx \omega t}$ of the linear wave equation in heterogeneous media. For real frequencies $\omega$, the well-posedness of this problem is unclear. In fact, on one hand, one expects that the physical solution $u$, if it exists, may not belong to $H^1(\R)$ due to possible wave propagation phenomena and a lack of decay at infinity. On the other hand, uniqueness of a solution in $H^1_{\textit{loc}}(\R)$ does not hold in general. In this case, one needs a so-called a \emph{radiation condition} that imposes the behaviour of the solution at infinity. Such a condition can be obtained in practice using the \emph{limiting absorption principle}, which consists in (\emph{i}) adding some absorption -- that is some imaginary part to $\omega$: $\Imag \omega$, and (\emph{ii}) studying the limit of the solution $u \equiv u(\omega)$ as the absorption tends to $0$. The limiting absorption principle is a classical approach to study time-harmonic wave propagation problems in unbounded domains; see for instance \cite{agmon1975spectral, eidus1986limiting, wilcox1966wave}. More recently, it has been successfully applied for locally perturbed periodic media \cite{flissthese, hoang2011limiting, kirsch2018radiation, radosz2015new}.

	\vspace{1\baselineskip} \noindent
	In this paper, we will only address the case with absorption, that is
 \begin{equation}\label{eq:dissipation}
 	\text{the frequency }\omega \text{ satisfies }\Imag \omega > 0.
 \end{equation}
 Under these assumptions,  \eqref{eq:whole_line_problem} admits a unique solution in $H^1(\R)$ by Lax-Milgram's theorem. Moreover, it can be shown (using for instance an argument similar to the one in \cite{combes1973asymptotic}) that this solution satisfies a sharp exponential decay property
 \begin{equation}\label{eq:exp_decay}
 	\spexists c,\,\alpha > 0, \quad \spforall \xv \in \R, \quad |u(\xv)|\leq c \,\euler^{-\alpha \Imag \omega|\xv|}.
 \end{equation}
Exploiting \eqref{eq:exp_decay}, a naive numerical method for treating the unboundedness would consist in truncating the computational domain (with homogeneous Dirichlet boundary conditions for instance) at a certain distance related to $\Imag \omega$. However the cost and the accuracy of the method would deteriorate when $\Imag \omega$ tends to $0$.
Our objective in this paper is to develop a numerical method which is robust when $\Imag \omega$ tends to $0$, in the particular case of locally perturbed quasiperiodic media. More precisely, we solve the problem in the bounded domain $(-a,a)$ (which is independent of $\Imag \omega$) by constructing transparent boundary conditions of Dirichlet-to-Neumann type:
	\begin{equation}\label{eq:DtN_coef}
		\displaystyle
		\pm\; \mu\; \frac{d u}{d \xv} + \lambda^\pm \; u = 0 \quad \textnormal{on} \quad \xv = \pm a,
	\end{equation}
	where $\lambda^\pm$ are called \emph{Dirichlet-to-Neumann} (DtN) coefficients. These coefficients are defined by
	\begin{equation}
		\displaystyle
		\lambda^\pm = \mp\; \Big[\mu\; \frac{d u^\pm}{d \xv}\Big](\pm a), \label{eq:DtN_coefficients}
	\end{equation}
where  $u^\pm$ is the unique solution in $H^1(\pm a, \pm \infty)$ of
	\begin{equation}
	\left|
	\begin{array}{r@{\ }c@{\ }l@{\quad}l}
	\displaystyle - \frac{d}{d \xv} \Big( \mu \; \frac{d u^\pm}{d \xv} \Big) - \rho \; \omega^2 \; u^\pm &=& 0, \quad \textnormal{for}&   \pm \xv > a,
	\\[8pt]
	\displaystyle u^\pm (\pm a) &=& 1.
	\end{array}
	\right.
	\label{eq:half_line_problems_0}
	\end{equation}
Knowing $\lambda^\pm$, one is then reduced to compute $u|_{(-a,a)}$ by solving the problem
	\begin{equation}
	\left|
	\begin{array}{r@{\ }c@{\ }l@{\quad}l}
	\displaystyle - \frac{d}{d \xv} \Big( \mu \; \frac{d u^i}{d \xv} \Big) - \rho \; \omega^2 \; u^i &=& f, \quad \textnormal{for} & \xv \in (-a, a),
	\\[10pt]
	\displaystyle \Big[\pm \mu\; \frac{d u^i}{d \xv} + \lambda^\pm \; u^i\Big](\pm a) &=& 0.
	\end{array}
	\right.
	\label{eq:interior_problem}
	\end{equation}
	The well-posedness of this problem is a direct consequence of the sign property
	\[
		\Imag \lambda^\pm < 0,
	\]
	which, through a Green's formula, results itself from the presence of dissipation \eqref{eq:dissipation} in \eqref{eq:half_line_problems_0}. Then the solution $u$ of \eqref{eq:whole_line_problem} is given by
		\begin{equation}\label{eq:solution_of_whole_line_problem}
		\spforall \xv \in \R, \quad u(\xv) = \left\{
		\begin{array}{c@{\quad}l}
			\displaystyle u^i(-a)\; u^-(\xv), & \xv < -a,\\[6pt]
			\displaystyle u^i(\xv), & \xv \in (-a, a),\\[6pt]
			\displaystyle u^i(a)\; u^+(\xv), & \xv > a.
		\end{array}
		\right.
		\end{equation}
		In general, the problem is that computing $\lambda^\pm$, that is to say solving \eqref{eq:half_line_problems_0}, is as difficult as the original problem. However, this is no longer true when the exterior medium (\emph{i.e.} outside $(-a,a)$) has a certain structure:
		\begin{itemize}
			\item if the exterior medium is homogeneous ($\rho$ and $\mu$ are constant), these coefficients can be computed explicitly;
			\item if the exterior medium is periodic ($\rho$ and $\mu$ are periodic), several methods for the computation of these DtN coefficients are developed in \cite{flissthese, jolyLiFliss, kirsch2018radiation};
			\item if the exterior medium is a weakly random perturbation of a periodic medium, the coefficients can be approximated via an asymptotic analysis; see \cite{fliss_giovangigli}.
		\end{itemize}
Our main objective in this paper is to compute the DtN coefficients for a quasiperiodic exterior medium, 
in order to develop a numerical method according to \eqref{eq:half_line_problems_0}, \eqref{eq:interior_problem}, \eqref{eq:solution_of_whole_line_problem}.

\vspace{1\baselineskip}\noindent
The outline of the rest of the paper is as follows. In Section \ref{sec:quasiperiodicity}, we introduce the fundamental notion of quasiperiodic functions (in $1$D) and define what is a locally perturbed quasiperiodic medium in the context of the problem \eqref{eq:whole_line_problem}. Sections \ref{sec:the_half_line_quasiperiodic_problem} and \ref{sec:resolution_half_guide_problem} are the most important sections of the paper. In Section \ref{sec:the_half_line_quasiperiodic_problem}, we link
	the solution of the 1D half-line problem with quasiperiodic coefficients to the solution of a degenerate directional Helmholtz equation posed in dimension $\ord$, with $\ord >1$ defined as in Section \ref{sec:quasiperiodicity}. This is the so-called lifting approach whose principle is presented in Section \ref{sec:lifting_in_a_higher_dimensional_periodic_problem}. More precisely, in Section \ref{sec:link_with_a_periodic_half_guide_problem}, we characterize the solution of the 1D quasiperiodic problem as the trace along a (broken) line of a $\ord$D problem posed in a domain  with the geometry of a half-waveguide: $(0,1)^{\ord-1} \times \R_+$. In between, we need to dedicate the (rather long) Section \ref{sec:preliminary_material} to fix the notations used in the rest of the paper and present some useful preliminary material about an adapted functional framework for the rigorous setting of our method. This concerns anisotropic Sobolev spaces with an emphasis on trace theorems and related Green's formula.
	In Section \ref{sec:resolution_half_guide_problem}, we provide a method for solving the half-waveguide problem of Section \ref{sec:link_with_a_periodic_half_guide_problem}. In Section \ref{sec:structure_of_the_solution}, we describe the structure of the solution with the help of a propagation operator ${\cal P}$ and local cell problems. In Section \ref{sec:Riccati}, we show that the operator ${\cal P}$ is characterized as a particular solution of a Riccati equation. In Section \ref{sec:the_DtN_operator_and_the_DtN_coefficient}, we first build a directional DtN operator $\Lambda$ for the half-waveguide problem, from which we  deduce the DtN coefficients $\lambda^\pm$ we are looking for (\emph{cf.} \eqref{eq:DtN_coefficients}). Finally, in  Section \ref{sec:about_Riccati_equation}, we analyze the Riccati equation from a spectral point of view and in Section \ref{sec:propagation_operator} we describe the spectrum of ${\cal P}$. In Section \ref{sec:resolution_algorithm} devoted to numerical results,  we restrict ourselves to $\ord=2$ for the sake of simplicity. The first two subsections are devoted to the discretization of the cell problems evoked above. We have considered two approaches: one, natural but naive, consists \textcolor{surligneur}{in using 2D} Lagrange finite elements (Section \ref{sec:methode_2D}) while the other, called the quasi-1D method, is better fitted to the anisotropy of the problem (Section \ref{sec:methode_quasi1D}). In Section \ref{sec:discrete_Riccati_equation}, we explain how we can construct a discrete propagation operator from a discrete Riccati equation that we choose to solve via a spectral approach, while Section \ref{sec:discrete_DtN_coefficient} simply mimics Section \ref{sec:the_DtN_operator_and_the_DtN_coefficient} at the discrete level. Section \ref{sec:numerical_results} is devoted to numerical results. In the first three subsections, we provide various validations of our method for the half-line problem (Sections \ref{sec:results:half_line_guide} and \ref{sec:results:absorption}) and the whole line problem (Section \ref{sec:results:whole_line_problem}). At last, in Section \ref{sec:results:spectral_approximation_P}, we address the question of the approximation of the spectrum of the propagation operator ${\cal P}$ by the one of its discrete approximation.

	\paragraph{Particular notation used throughout the paper.} In what follows, 
	\begin{enumerate}
		\item the equality modulo $1$ is denoted by
		\[
		\spforall y\in\R, \quad \  z = y\,[1] \quad \Longleftrightarrow \quad z\in[0,1) \ \  \text{and} \ \  y-z \in \Z.
		\]
		and for all $p, q \in \N,\  p < q$, we set $\llbracket p, q\rrbracket:=\{j\in\N,\ p\leq j\leq q\}$.
		\item We introduce $\mathscr{C}_{\textit{per}}(\R^\ord)$ as the space of continuous functions $F : \R^\ord \to \R$ that are $1$--periodic with respect to each variable, and $\mathscr{C}^\infty_0(\mathcal{O})$ as the space of smooth functions that are compactly supported in $\mathcal{O} \subset \R^\ord$.
		\item For $i \in\llbracket 1,\ord\rrbracket$, we denote by $\vec{\ev}_i$ the $i$-th unit vector from the canonical basis of $\R^\ord$. For any element $\yv = (\yi_1, \dots, \yi_\ord)$ in $\R^\ord$, we define $\hat{\yv}$ as the vector $(\yi_1, \dots, \yi_{\ord-1}) \in \R^{\ord-1}$, so that $\yv = (\hat{\yv}, \yi_\ord)$. For $\yv = (\yi_1, \dots, \yi_\ord)$ and $\itbf{z} = (z_1, \dots, z_\ord)$, the Euclidean inner product of $\yv$ and $\itbf{z}$ is denoted $\yv \cdot \itbf{z} := \yi_1\,z_1 + \cdots \yi_\ord\,z_\ord$, and the associated norm is $|\yv| := \sqrt{\yv \cdot \yv}$.
	\end{enumerate}

\section{Quasiperiodicity}
\label{sec:quasiperiodicity}
\subsection{Quasiperiodic functions of one real variable} 
\begin{figure}[H]
	\centering
	\begin{tikzpicture}
		\def\cote{4.5cm}
		\begin{scope}[xshift=0\textwidth]
			\begin{axis}[
				enlargelimits=false,
				axis on top,
				width=\cote, height=\cote,
				xtick = {0, 0.4, 0.8}, ytick = {0, 0.4, 0.8},
				disabledatascaling,
				axis equal,
				colorbar,
				point meta min=-2, point meta max=2,
				colorbar style={
					width=0.05*\pgfkeysvalueof{/pgfplots/parent axis width},
					yticklabel style={
						/pgf/number format/precision=2,
						/pgf/number format/fixed,
					},
					ytick={0},
					scaled y ticks=false,
					major tick length=0.05*\pgfkeysvalueof{/pgfplots/parent axis width},
				},
				unbounded coords=jump,  
			]
				\ifthenelse{\boolean{afficherGraphes}}{
					\addplot graphics [xmin=0, xmax=1, ymin=0, ymax=1] {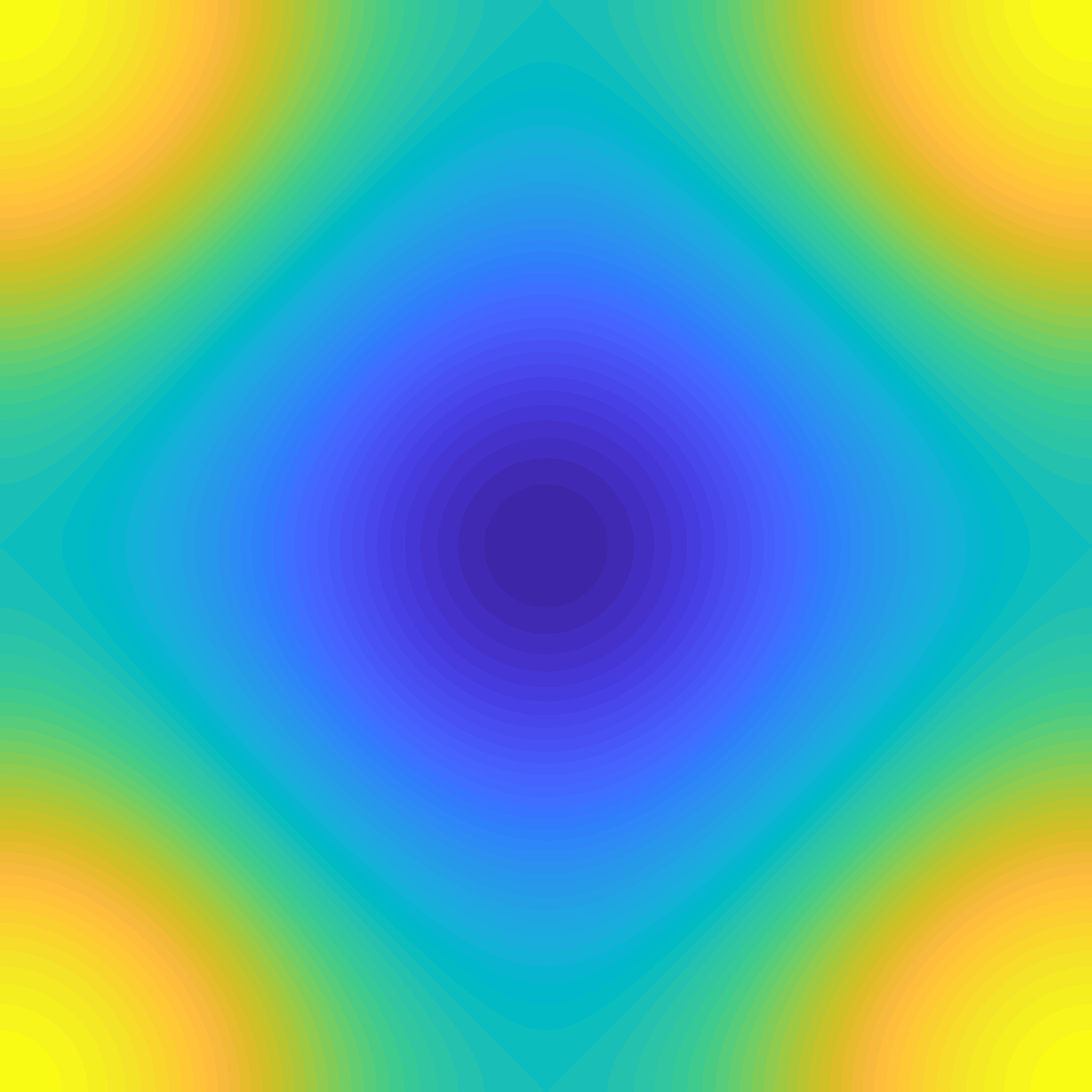};
				}{}
				\draw [dashed, thick] (axis cs:0, 0) -- (axis cs:0.7071, 1);
				\draw[->, thick] (axis cs:0, 0) -- (axis cs:0.1768, 0.25);
				\draw (axis cs:0.1, 0.095) node[right] {$\cut$};
			\end{axis}
		\end{scope}
		\begin{scope}[xshift=\cote]
			\begin{axis}[
				width=\textwidth-\cote,
				height=\cote,
				xmin=-4, xmax=4,
				ymin=-2.1, ymax=2.1,
				xtick = {-4, -2, 0, 2, 4},
				ytick = {-2, 0, 2},
				enlargelimits=false,
				axis on top,
				]
				\ifthenelse{\boolean{afficherGraphes}}{
					\addplot+[domain=-5:5, smooth, no marks, samples=400, color=bleuDoux!75, thick] {cos(deg(2*pi*x)) + cos(deg(2*pi*sqrt(2)*x))};
				}{}
				\draw[thick, color=black] (axis cs: 0, -1.8) -- (axis cs:0, -2) -- (axis cs:0.7071, -2) -- (axis cs:0.7071, -1.8);
				\draw (axis cs:0.7071, -1.9) node[right] {\textit{\footnotesize Size of periodicity cell}};
			\end{axis}
		\end{scope}
	\end{tikzpicture}
	\caption{Function $F: (\yi_1, \yi_2) \mapsto  \cos 2\pi \yi_1 + \cos 2\pi \yi_2$ in its periodicity cell (left), and whose trace along $\cut = (1, \sqrt{2})$ leads to a quasiperiodic function (right). \label{fig:example_quasiperiodic_function}}
\end{figure}
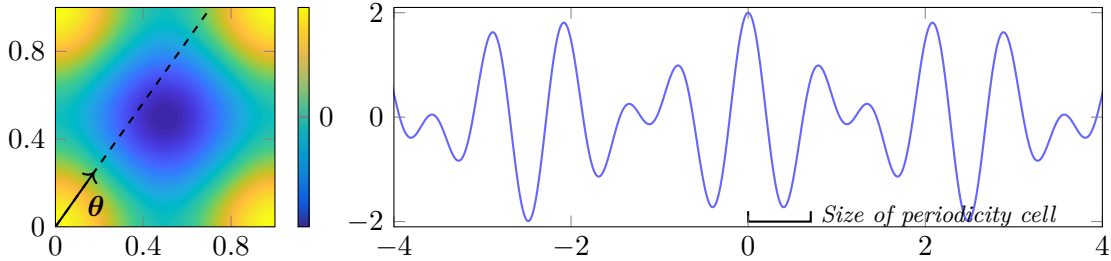

In this section, we present a brief overview of the main properties of quasiperiodic functions. We refer to \cite{besicovitch, bohr, levitan} for more complete presentations. Quasiperiodicity is defined as follows.
\begin{defi}
	\label{def:quasiperiodic_function}
	A continuous function $f : \R \to \R$ is said to be \emph{quasiperiodic of order $n > 1$} if there exist a constant real vector $\cut = (\cuti_1, \dots, \cuti_\ord)$, with $\cuti_i > 0$ for all $i \in\llbracket 1,\ord\rrbracket$, and a continuous function $F : \R^\ord \to \R$, $1$--periodic with respect to each variable, such that
	\begin{equation}
		\label{eq:def_quasiperiodic_function}
		\displaystyle
		\spforall \xv \in \R, \quad f(\xv) = F(\xv\,\cut).
	\end{equation}
	The vector $\cut$ is called a \emph{cut direction}, and $F$ is a periodic \emph{extension} of $f$.

	\vspace{1\baselineskip} \noindent
	A geometrical interpretation of this definition is to see the one-dimensional function $f$ as the trace of a $\ord$-dimensional function $F$ along the line passing through $(0, 0)$ and parallel to the vector $\cut$. This is illustrated in Figure \ref{fig:example_quasiperiodic_function} for $\ord = 2$ and $\cut = (1, \sqrt{2})$.
\end{defi}

%

\vspace{0\baselineskip} \noindent
	Periodic functions are obviously quasiperiodic. Other examples of quasiperiodic functions are finite sums or products of periodic functions: if $f_1$ and $f_2$ are periodic, then $f_1 + f_2$ and $f_1 f_2$ can be expressed under the form \eqref{eq:def_quasiperiodic_function}. Note that $f_1 + f_2$ and $f_1 f_2$ are \textit{not} periodic if $f_1$ and $f_2$ are \textit{continuous} functions with non-commensurable least periods. For instance, with $f_1(x) = \cos 2\pi x$ and $f_2(x) = \cos 2\pi\sqrt{2} x$, one easily checks that the sum $f_1 + f_2$, represented in Figure \ref{fig:example_quasiperiodic_function}, is not periodic since it equals $2$ only when $x = 0$.
%
%

\vspace{1\baselineskip}\noindent
In Definition \ref{def:quasiperiodic_function}, it is easy to see that neither the periodic extension nor the cut direction are uniquely defined. Given  $(F, \cut)$, it is always possible to lower the value of  $n$, and change the function $F$ accordingly, so that the coefficients $\cuti_1,\dots,\cuti_\ord$ are \emph{linearly independent} \emph{over the integers} (see \cite[Chapter 2]{levitan}), that is
\begin{equation}
	\label{eq:linear_independance}
	\displaystyle
	\spforall \textit{\textbf{k}} \in \Z^\ord, \qquad \textit{\textbf{k}}\cdot\cut = 0 \quad \Longleftrightarrow \quad \textit{\textbf{k}}= 0.
\end{equation}
For $\ord = 2$ and $\cut = (\cuti_1, \cuti_2)$, the above condition amounts to saying that the ratio $\cuti_1 / \cuti_2$ is irrational. Due to this observation, vectors that satisfy \eqref{eq:linear_independance} will be abusively referred to as \emph{irrational vectors}. 
 A consequence of \eqref{eq:linear_independance} is given by Kronecker's approximation theorem.
\begin{thm}[\textnormal{\cite[Theorem 444]{hardy}}]
	\label{thm:kronecker}
	If $\cut$ is an irrational vector, then the set $\cut\, \R + \N^\ord$ is dense in $\R^\ord$. 
\end{thm}
\noindent
If $\cut$ is an irrational vector, and if $F \in \mathscr{C}_{\textit{per}}(\R^\ord)$ satisfies $F(\cut\, \R) = 0$, then Theorem \ref{thm:kronecker} ensures that $F = 0$. In other words, under the linear independence assumption, $F$ is uniquely determined by its restriction on the line $\cut\, \R$.

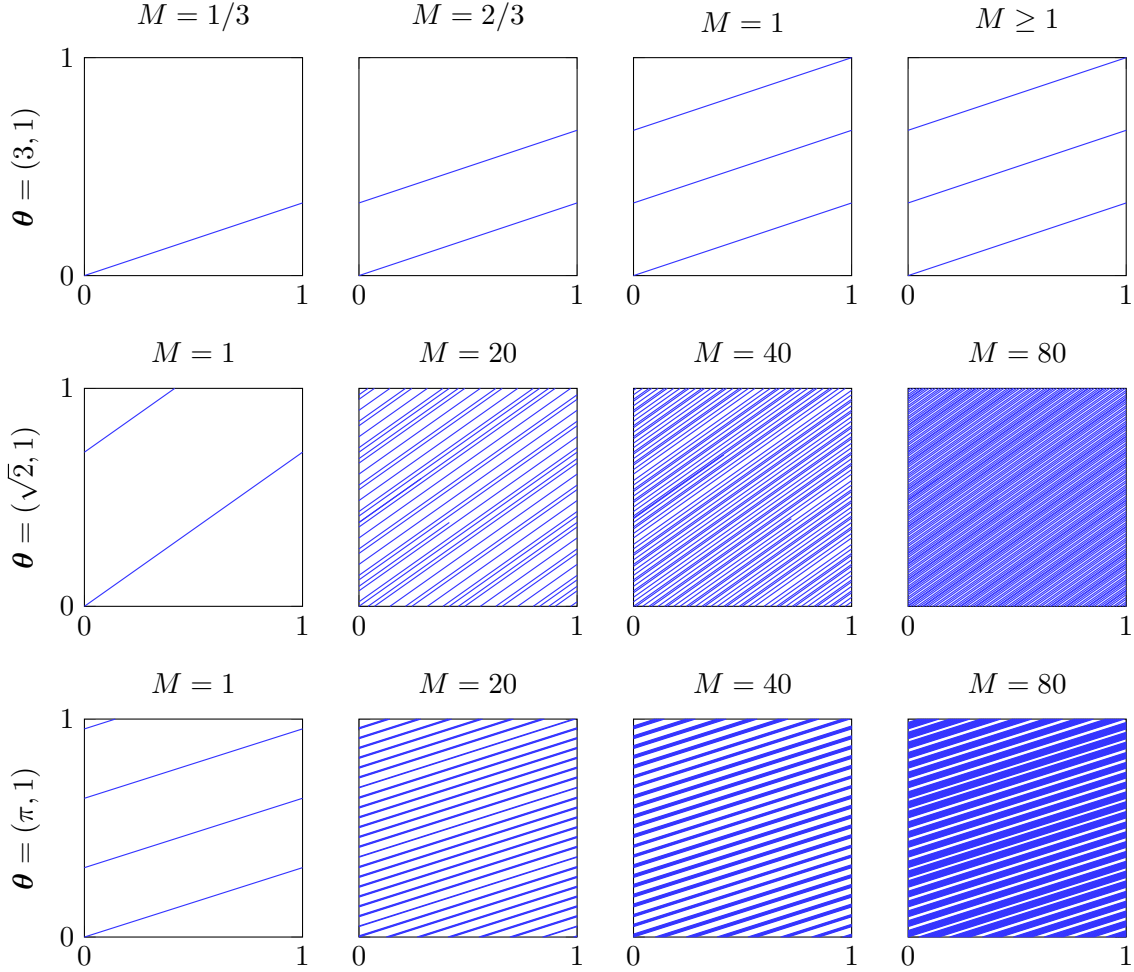
\begin{figure}[ht!]
	\begin{tikzpicture}
		\pgfmathsetmacro\coteplot{2.9}
		\begin{groupplot}[
			group style={
				group name=broken_line_plots,
		    group size=4 by 3,
		    horizontal sep=0.75cm,
				vertical sep=1.5cm,
				ylabels at=edge left,
		    yticklabels at=edge left,
			},
      width =\coteplot cm,
			height=\coteplot cm,
      xmin=0, xmax=1, ymin=0, ymax=1,
      xtick={0, 1}, ytick={0, 1},
      enlargelimits=false,
      axis equal,
      unbounded coords=jump,
			scale only axis,
			]
			\nextgroupplot[ylabel={$\cut = (3, 1)$}, title={$M = 1/3$}]
			\nextgroupplot[title={$M = 2/3$}]
			\nextgroupplot[title={$M = 1$}]
			\nextgroupplot[title={$M \geq 1$}]
			\nextgroupplot[ylabel={$\cut = (\sqrt{2}, 1)$}, title={$M = 1$}]
			\nextgroupplot[title={$M = 20$}]
			\nextgroupplot[title={$M = 40$}]
			\nextgroupplot[title={$M = 80$}]
			\nextgroupplot[ylabel={$\cut = (\pi, 1)$}, title={$M = 1$}]
			\nextgroupplot[title={$M = 20$}]
			\nextgroupplot[title={$M = 40$}]
			\nextgroupplot[title={$M = 80$}]
		\end{groupplot}

		\ifthenelse{\boolean{afficherGraphes}}{
			\node [] at (broken_line_plots c1r1.center) {%
				\begin{tikzpicture}%
					\draw[opacity=0] (0, 0) rectangle ++(\coteplot, \coteplot);%
					\genererFibrage{\coteplot}{0}{0}{3}{1}{1}{1}{line width=0.15mm, bleuDoux}
				\end{tikzpicture}%
			};
			\node [] at (broken_line_plots c2r1.center) {%
				\begin{tikzpicture}%
					\draw[opacity=0] (0, 0) rectangle ++(\coteplot, \coteplot);%
					\genererFibrage{\coteplot}{0}{0}{3}{1}{2}{1}{line width=0.15mm, bleuDoux}
				\end{tikzpicture}%
			};
			\node [] at (broken_line_plots c3r1.center) {%
				\begin{tikzpicture}%
					\draw[opacity=0] (0, 0) rectangle ++(\coteplot, \coteplot);%
					\genererFibrage{\coteplot}{0}{0}{3}{1}{3}{1}{line width=0.15mm, bleuDoux}
				\end{tikzpicture}%
			};
			\node [] at (broken_line_plots c4r1.center) {%
				\begin{tikzpicture}%
					\draw[opacity=0] (0, 0) rectangle ++(\coteplot, \coteplot);%
					\genererFibrage{\coteplot}{0}{0}{3}{1}{4}{1}{line width=0.15mm, bleuDoux}
				\end{tikzpicture}%
			};
			\node [] at (broken_line_plots c1r2.center) {%
				\begin{tikzpicture}%
					\draw[opacity=0] (0, 0) rectangle ++(\coteplot, \coteplot);%
					\genererFibrage{\coteplot}{0}{0}{sqrt(2)}{1}{3}{1}{line width=0.15mm, bleuDoux}
				\end{tikzpicture}%
			};
			\node [] at (broken_line_plots c2r2.center) {%
				\begin{tikzpicture}%
					\draw[opacity=0] (0, 0) rectangle ++(\coteplot, \coteplot);%
					\genererFibrage{\coteplot}{0}{0}{sqrt(2)}{1}{200}{20}{line width=0.15mm, bleuDoux}
				\end{tikzpicture}%
			};
			\node [] at (broken_line_plots c3r2.center) {%
				\begin{tikzpicture}%
					\draw[opacity=0] (0, 0) rectangle ++(\coteplot, \coteplot);%
					\genererFibrage{\coteplot}{0}{0}{sqrt(2)}{1}{200}{40}{line width=0.15mm, bleuDoux}
				\end{tikzpicture}%
			};
			\node [] at (broken_line_plots c4r2.center) {%
				\begin{tikzpicture}%
					\draw[opacity=0] (0, 0) rectangle ++(\coteplot, \coteplot);%
					\genererFibrage{\coteplot}{0}{0}{sqrt(2)}{1}{500}{80}{line width=0.15mm, bleuDoux}
				\end{tikzpicture}%
			};
			\node [] at (broken_line_plots c1r3.center) {%
				\begin{tikzpicture}%
					\draw[opacity=0] (0, 0) rectangle ++(\coteplot, \coteplot);%
					\genererFibrage{\coteplot}{0}{0}{pi}{1}{4}{1}{line width=0.15mm, bleuDoux}
				\end{tikzpicture}%
			};
			\node [] at (broken_line_plots c2r3.center) {%
				\begin{tikzpicture}%
					\draw[opacity=0] (0, 0) rectangle ++(\coteplot, \coteplot);%
					\genererFibrage{\coteplot}{0}{0}{pi}{1}{200}{20}{line width=0.15mm, bleuDoux}
				\end{tikzpicture}%
			};
			\node [] at (broken_line_plots c3r3.center) {%
				\begin{tikzpicture}%
					\draw[opacity=0] (0, 0) rectangle ++(\coteplot, \coteplot);%
					\genererFibrage{\coteplot}{0}{0}{pi}{1}{200}{40}{line width=0.15mm, bleuDoux}
				\end{tikzpicture}%
			};
			\node [] at (broken_line_plots c4r3.center) {%
			\begin{tikzpicture}%
				\draw[opacity=0] (0, 0) rectangle ++(\coteplot, \coteplot);%
				\genererFibrage{\coteplot}{0}{0}{pi}{1}{500}{80}{line width=0.15mm, bleuDoux}
			\end{tikzpicture}%
		};
		}{}
	\end{tikzpicture}
	\caption{Representation of the set $\big\{(x\,\cuti_1[1],x\,\cuti_2[1]),\;x \in (0, M)\big\}$ in $(0, 1)^2$ for different values of M, when $\cuti_1/\cuti_2 \in \Q$ (first row), and when $\cuti_1/\cuti_2 \in \R \setminus \Q$ (second row for $\cut = (\sqrt{2}, 1)$ and third row for $\cut = (\pi, 1)$).\label{fig:fibrage}}
\end{figure}

\vspace{1\baselineskip}\noindent
For $\ord = 2$, Theorem \ref{thm:kronecker} implies that the broken line $\big\{(x\,\cuti_1[1],x\,\cuti_2[1]),\;x\in\R\big\}$ is dense in the unit cell $(0,1)^2$. To illustrate this, Figure \ref{fig:fibrage} represents the set $\big\{(x\,\cuti_1[1],x\,\cuti_2[1]),\;x\in(0,M)\big\}$ in the unit cell for different values of $M$, when (\emph{1}) $\theta_1/\theta_2 \in \Q$ (see the first row), and when (\emph{2}) $\theta_1/\theta_2 \in \R \setminus \Q$ (see the second row for $\cut = (\sqrt{2}, 1)$ and the third one for $\cut = (\pi, 1)$). For $M$ large enough, in the first case, this set is reduced to a \textit{finite} union of segments, whereas in the second case, it seems to fill the unit cell without ever passing through the same positions. It is also interesting to see that for $\cut = (\sqrt{2}, 1)$, the unit cell is somehow filled uniformly, contrary to the case where $\cut = (\pi, 1)$.
%
%
%

%

\vspace{1\baselineskip} \noindent
Finally, it is worth mentioning that Definition \ref{def:quasiperiodic_function} extends to higher-dimensional continuous functions as well. Moreover, the notion of quasiperiodicty can be defined at a discrete level, to describe the properties of tilings that are cuts and projections of higher-dimensional periodic tilings. These quasiperiodic tilings have been extensively studied \cite{gardner_1977_aperiodic, meyer1995quasicrystals, penrose_pentaplexity_1979, senechal1996quasicrystals}, and are used for modelling quasicrystals \cite{shechtmanAl}.
%
%
%
%
%
%
%
%
%
%
%
%
%
%
\subsection{Locally perturbed quasiperiodic media}
A locally perturbed quasiperiodic medium is a medium corresponding to functions $\mu$ and $\rho$ that satisfy \eqref{eq:coef_ellipt} and that are quasiperiodic outside a bounded interval, which can be supposed to be $(-a,a)$ (see \eqref{eq:source_terme}) without any loss of generality. More precisely,
	\[
		\mu(x) = \left|
		\begin{array}{c l}
			\mu_i(x) 			  & x \in (-a,a)\esplignestableau
			\mu_p(\xv\,\cut ) & x \in \R \setminus (-a,a)
		\end{array}
		\right.
		\quad \textnormal{and} \quad
		\rho(x) = \left|
		\begin{array}{c l}
			\rho_i(x) 			 & x \in (-a,a)\esplignestableau
			\rho_p( \xv\,\cut) & x \in \R \setminus (-a,a),
		\end{array}
		\right.
	\]
	where the functions $\mu_p,\, \rho_p$ belong to $\mathscr{C}_{\textit{per}}(\R^\ord)$ with $\ord > 1$, and $\cut \in \R^\ord$ is an irrational vector (see Condition \eqref{eq:linear_independance}).

	\begin{rmk}\label{rem}
		(a). Since $\cut$ is an irrational vector, Kronecker's approximation theorem \ref{thm:kronecker} ensures that the functions $\mu_p$ and $\rho_p$ are entirely determined by their restrictions on the line $\R\, \cut$. Therefore, $\mu_p$ and $\rho_p$ satisfy \eqref{eq:coef_ellipt}
		with respectively the same bounds as $\mu$ and $\rho$.

		\vspace{1\baselineskip}
		(b). The present study \textcolor{surligneur}{can be extended} without difficulty to the case where $\mu$ (resp. $\rho$) coincides with two different quasiperiodic functions in $(-\infty, -a)$ and in $(a, +\infty)$:
		\[
			\text{for}\; \pm x>\pm a, \quad\  \mu(x) =
				\mu_p^\pm(\xv\,\cut^\pm\, )   \quad \text{and} \quad \rho(x) =
				\rho_p^\pm(\xv\,\cut^\pm\, ),
		\]
		where $\mu_p^\pm,\, \rho_p^\pm$ belong to $\mathscr{C}_{\textit{per}}(\R^{\ord^\pm})$ with $\ord^\pm > 1$, and where $\cut^\pm \in \R^{\ord^\pm}$ are irrational vectors.
		\end{rmk}
	\section{The half-line quasiperiodic problem} 
	\label{sec:the_half_line_quasiperiodic_problem}
	We now focus on the half-line quasiperiodic problems \eqref{eq:half_line_problems_0}. As these problems are very similar to each other, it is sufficient to study the half-line problem set on $(a,+\infty)$ and suppose without loss of generality that $a = 0$.  Let $\mu_\cut := \mu_p(\cut\,\cdot)$ and $\rho_\cut := \rho_p(\cut\,\cdot)$. Therefore, the problem we consider in this section is the following:
	\begin{equation}
	\left|
	\begin{array}{r@{\ }c@{\ }l@{\quad}l}
	\displaystyle - \frac{d}{d \xv} \Big( \mu_\cut \; \frac{d u^+_\cut}{d \xv} \Big) - \rho_\cut \; \omega^2 \; u^+_\cut &=& 0, \quad \textnormal{in} & \R_+, \\[8pt]
	\displaystyle u^+_\cut(0) &=& 1.
	\end{array}
	\right.
	\label{eq:half_line_problem}
	\end{equation}
	\begin{rmk}
		\color{surligneur}
		The function $u^+_\cut$ corresponds exactly to the solution $u^+$ of \eqref{eq:half_line_problems_0} that was introduced in Section \ref{sec:introduction_motivation} for very general media. The reason why this solution is relabeled $u^+_\cut$ is due to the fact that, because we consider here quasiperiodic media, the coefficients $\mu$ and $\rho$ that appear in \eqref{eq:half_line_problems_0} have been replaced by $\mu_\cut$ and $\rho_\cut$.
	\end{rmk}

	\subsection{Lifting in a higher-dimensional periodic problem}
	\label{sec:lifting_in_a_higher_dimensional_periodic_problem}
\begin{figure}
		\centering
		\setlength{\belowcaptionskip}{-10pt}
		\begin{tikzpicture}[thick, scale=0.8]
			\def\cotecellule{2}
			\def\nbcellules{3}
			%
			%
			\begin{scope}[xshift=0]
				\pgfmathsetmacro\lenX{4*\cotecellule}
				\pgfmathsetmacro\bandeY{\nbcellules*\cotecellule+0.5*\cotecellule}
				\pgfmathsetmacro\lenY{\nbcellules*\cotecellule+0.5*\cotecellule}
				\def\anglecoupe{60}
				\pgfmathsetmacro\chgtvarS{0*\cotecellule}
				\def\chgtvarX{3.0}
				\fill[black!30, opacity=0.5] (-4.875*\cotecellule , 0) -- (-4.875*\cotecellule , \bandeY) -- (3.875*\cotecellule , \bandeY) -- (3.875*\cotecellule , 0) -- cycle;
				\fill[black!50, opacity=0.5] (-1*\cotecellule , 0) -- (-1*\cotecellule , \bandeY) -- (0 , \bandeY) -- (0 , 0) -- cycle;
				%
				%
				%
				\foreach \idI in {-4,-3,-2,-1,1,2,3} {
				\draw[black!40] (\idI*\cotecellule, 0) -- +(0, \lenY);
				}
				\draw[->] (-\lenX-\cotecellule, 0) -- (\lenX, 0) node[right] {$\yi_1$};
				\draw[->] (-1*\cotecellule, -0.25*\cotecellule) -- (-1*\cotecellule, \lenY) node[below left]  {$\yi_2$};
				\draw (0, 0) -- +(0, \lenY);
				\draw[rougeTerre] (-1*\cotecellule, 0) node {$\bullet$} -- +({\bandeY*cos(\anglecoupe)/sin(\anglecoupe)}, \bandeY);
				\draw ({-0.25-1*\cotecellule+\bandeY*cos(\anglecoupe)/sin(\anglecoupe)}, \lenY) node[rougeTerre, below left]{$\cut\, \R_+$};
				%
				%
				\draw[-latex] (0.875*\lenX, 0) -- +({1.0*cos(\anglecoupe)}, {1.0*sin(\anglecoupe)}) node [pos=1.25] {$\cut$};
				%
				%
				\draw (-\cotecellule, 0) node[below left] {$0$};
				\node[rougeTerre] (pt1Eq1D) at ({-1*\cotecellule + 0.6*\bandeY*cos(\anglecoupe)/sin(\anglecoupe)}, 0.6*\bandeY) {};
				\node[rougeTerre, scale=0.85] (pt2Eq1D) at (2.25*\cotecellule, 2.75*\cotecellule) {};
				\node[rougeTerre] (pt1BC1D) at (-1*\cotecellule, 0) {};
				\node[rougeTerre, scale=0.85] (pt2BC1D) at (-2*\cotecellule, \cotecellule) {}; 
				\node (pt1Eq2D) at (-0.5*\cotecellule, 2*\cotecellule) {};
				\node[scale=0.85] (pt2Eq2D) at (-3.25*\cotecellule, 2.75*\cotecellule) {}; 
				\node (pt1BC2D) at (1.5*\cotecellule, 0) {};
				\node[scale=0.85] (pt2BC2D) at (1.5*\cotecellule, \cotecellule) {}; 
				\draw[latex-, rougeTerre] (pt1Eq1D) to[bend right=50] (pt2Eq1D);
				\draw[latex-, rougeTerre] (pt1BC1D) to[bend left=30] (pt2BC1D);
				\draw[latex-, black!75] (pt1Eq2D) to[bend right=30] (pt2Eq2D);
				\draw[latex-, black!75] (1.5*\cotecellule, 0) to (pt2BC2D);
				\draw[fill=white, draw=rougeTerre] ({2.25*\cotecellule-1.4*\cotecellule}, {2.75*\cotecellule-0.3*\cotecellule}) rectangle +(2.8*\cotecellule, 0.6*\cotecellule) node[pos=0.5, scale=0.8, rougeTerre] {$\displaystyle - \frac{d}{d \xv} \Big( \mu_\cut \; \frac{d u^+_\cut}{d \xv} \Big) - \rho_\cut \; \omega^2 \; u^+_\cut = 0$};
				\draw[fill=white, draw=black!80] ({-3.25*\cotecellule-1.5*\cotecellule}, {2.75*\cotecellule-0.3*\cotecellule}) rectangle +(3*\cotecellule, 0.6*\cotecellule) node[pos=0.5, scale=0.8] {$- \Dt{} \big( \mu_p \; \Dt{} \widetilde{U}^+_\cut \big) - \rho_p \; \omega^2 \; \widetilde{U}^+_\cut = 0$};;
				\draw[fill=white, draw=rougeTerre] ({-2*\cotecellule-0.6*\cotecellule}, {\cotecellule-0.25*\cotecellule}) rectangle +(1.2*\cotecellule, 0.5*\cotecellule) node[pos=0.5, scale=0.8, rougeTerre] {$\displaystyle u^+_\cut(0) = 1$};
				\draw[fill=white, draw=black!80] ({1.5*\cotecellule-0.6*\cotecellule}, {\cotecellule-0.25*\cotecellule}) rectangle +(1.2*\cotecellule, 0.5*\cotecellule) node[pos=0.5, scale=0.8] {$\displaystyle \widetilde{U}^+_\cut = \widetilde{\varphi}$};
			\end{scope}
		\end{tikzpicture}
		\caption{Illustration of the lifting approach for $\ord = 2$\label{fig:illustration_lifting_approach}}
	\end{figure}
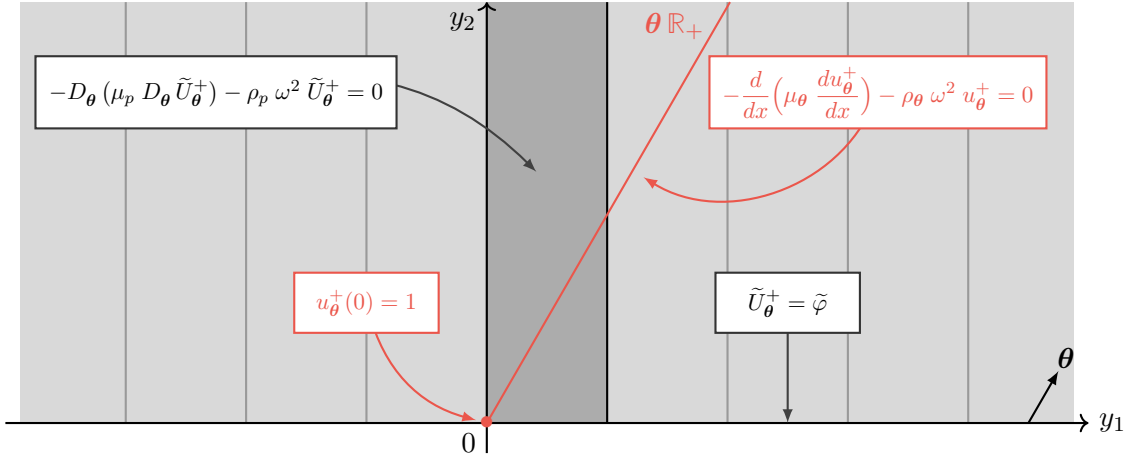

	We wish to exhibit some structure of the solution $u^+_\cut$. As the coefficients $\mu_\cut$ and $\rho_\cut$ in \eqref{eq:half_line_problem} are by definition traces of $\ord$--dimensional functions along the half-line $\cut\, \R_+$, it is natural to seek $u^+_\cut$ as the trace along the same line of a {\color{surligneur}function $\yv \in \R^\ord \mapsto \widetilde{U}^+_\cut(\yv)$, that is to say:
	\begin{equation}
		\aeforall \xv \in \R, \quad u^+_\cut(\xv) = \widetilde{U}^+_\cut (\xv\, \cut),
	\end{equation}
	where $\widetilde{U}^+_\cut$ shall be characterized as the solution of a $\ord$--dimensional PDE (in some sense, an “augmented” problem in which $\yv$ is the augmented space variable) with periodic coefficients, as illustrated in Figure \ref{fig:illustration_lifting_approach}}. This so-called \textit{lifting approach} has been used in the homogenization setting for the analysis of some correctors in presence of periodic halfspaces \cite{gerard2012, gerard2011} or periodic structures separated by an interface \cite{blancLeBrisLions}, as well as for the homogenization of quasicrystals and Penrose tilings \cite{bouchitte, wellanderAl}. However, to our knowledge, very little seems to have been done in other contexts (such as wave propagation), and in particular for numerical analysis and simulation purposes.

	\vspace{1\baselineskip} \noindent
	To build a higher-dimensional PDE, one has to exploit the correspondence between the derivative of $u^+_\cut$ and the partial derivatives of $\widetilde{U}^+_\cut$: according to the chain rule, for any smooth enough function $F: \R^\ord \to \C$, one has
	\begin{equation}
		\label{eq:chain_rule}
		\displaystyle
		\spforall \xv \in \R, \quad \frac{d}{d \xv}[ F(\cut\, \xv) ] = ( \Dt{} F ) (\cut\, \xv), \quad \text{with} \quad \Dt{} = \cut \cdot \nabla = \sum_{i = 1}^\ord \cuti_i\, \frac{\partial}{\partial \yi_i}.
	\end{equation}
	This leads us to introduce the $\ord$--dimensional PDE set on a half-space (see Remark \ref{rmk:lifting_approach})
	\begin{subequations}
		\label{eq:half-space_problem}
	  \begin{align}
	   \displaystyle - \Dt{} \big( \mu_p \; \Dt{} \widetilde{U}^+_\cut \big) - \rho_p \; \omega^2 \; \widetilde{U}^+_\cut &= 0, \quad \textnormal{for} \quad \yi_\ord > 0,
	 \intertext{where we recall that the coefficients $\mu_p,\, \rho_p : \R^\ord \to \R$ are continuous and $1$--periodic with respect to each variable. In addition, the boundary condition in \eqref{eq:half_line_problem} can be lifted onto the inhomogeneous Dirichlet boundary condition}
	   \displaystyle \widetilde{U}^+_\cut &= \widetilde{\varphi}, \quad \textnormal{on} \quad \yi_\ord = 0,
	  \end{align}
	\end{subequations}
	\noindent
	where the data $\widetilde{\varphi} : \R^{\ord-1} \to \C$ could be chosen continuous and must satisfy $\widetilde{\varphi}(0) = 1$, for the sake of consistency with the fact that $u^+_\cut(0) = 1$. Furthermore, to exploit the periodicity of the coefficients $\mu_p$ and $\rho_p$ with respect to the transverse variables $y_j,\,j<n$, we \textcolor{surligneur}{can} impose the following:
	\begin{equation}
		\label{eq:varphi_per}
		\text{$\widetilde{\varphi}$ is $1$--periodic,}
	\end{equation}
	so that it is natural to impose that
	\begin{equation}
		\label{eq:U_theta_periodic}
		\text{$\widetilde{U}^+_\cut(\varphi)$ is $1$--periodic with respect to the transverse variables $\yi_j,\ j<n$.}
	\end{equation}
	In Section \ref{sec:link_with_a_periodic_half_guide_problem}, we show how to reduce the above to a half-guide problem with periodic coefficients. In order to do so, we shall need some preliminary materials, which is the object of the next section.
	\begin{rmk}\label{rmk:lifting_approach}
		(a). One could have defined the augmented problem \eqref{eq:half-space_problem} on other half-spaces $\{\yv \in \R^\ord,\ \yi_i > 0\}$. The choice of the half-space is purely arbitrary.

		\vspace{1\baselineskip}
		(b). At first glance, one could imagine restricting the whole study to a constant boundary data $\widetilde{\varphi} = 1$. Though, in practice, this can be the case, the method used to solve the half-guide problem requires to investigate the structure of $\widetilde{U}^+_\cut(\widetilde{\varphi})$ for any $\widetilde{\varphi}$ in an appropriate function space (see Section \ref{sec:resolution_half_guide_problem} for more details).
	\end{rmk}

	\subsection{Preliminary material}
	\label{sec:preliminary_material}
	The main objective of this section is to establish rigorously some Green's formulas that are formally obvious, such as the one of Proposition \ref{prop:Green_formula_H1per}. This requires first to introduce the adapted functional framework and, since Green's formulas involve boundary integrals, to establish relevant trace theorems. Section \ref{sub:trace} is devoted to these trace theorems, while we present the corresponding Green's formulas in Section \ref{sub:Green}. Finally, Section \ref{sub:oblique_cov} highlights a simple but useful link between the derivative $\Dt{}$ and a single partial derivative with respect to one real variable, through a so-called oblique change of variables.

	\subsubsection{Anisotropic Sobolev spaces and trace theorems}\label{sub:trace}
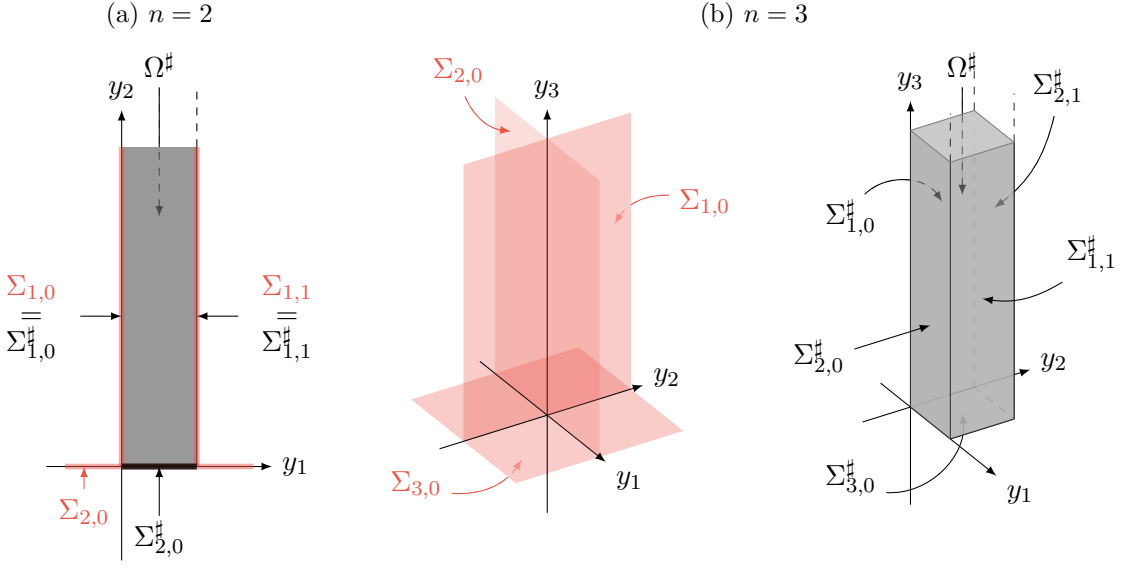
\begin{figure}
		\centering
		\begin{subfigure}[t]{0.34\textwidth}
			\centering
			\caption{$\ord = 2$}
			\def\lenX{1.0}
			\def\lenY{5.0}
			\def\hppXmax{1.75*\lenX}\def\hppXmin{-1.75*\lenX+1}
			\def\hppYmax{0.85*\lenY}
			\begin{tikzpicture}[scale=1.0, every node/.style={scale=1.0}]
				\begin{scope}
					\node (ptOmegaSharp1) at (0.5, 0.75*\hppYmax) {};
					\node (ptOmegaSharp2) at (0.5, 1.25*\hppYmax) {$\Omega^\diese$};
					\draw[latex-] (ptOmegaSharp1) to (ptOmegaSharp2);
					\fill[white] (0, 0) -- (1, 0) -- (1, \hppYmax) -- (0, \hppYmax) -- cycle;
					\fill[black!50, fill opacity=0.8] (0, 0) -- (1, 0) -- (1, \hppYmax) -- (0, \hppYmax) -- cycle;
					%
					\node (zaxis) at (0, \lenY) {$y_2$};
					\draw[-latex] (-\lenX, 0) -- (\lenX+1, 0) node[right] {$\yi_1$};
					\draw[-latex] (0, -0.25*\lenY) -- (zaxis);
					\draw (1, 0) -- (1, \hppYmax);
					\draw[dashed] (1, \hppYmax) -- (1, \lenY);
					%
					\draw[rougeTerre, line width=0.75mm, opacity=0.5] (0, 0) -- (0, \hppYmax);
					\draw[rougeTerre, line width=0.75mm, opacity=0.5] (1, 0) -- (1, \hppYmax);
					\draw[rougeTerre, line width=0.75mm, opacity=0.5] (\hppXmin, 0) -- (\hppXmax, 0);
					\draw[line width=0.8mm, opacity=0.8] (0, 0) -- (1, 0);
					%
					%
					\draw[latex-, opacity=0.5, dashed] (ptOmegaSharp1) to (ptOmegaSharp2);
					\node (ptSigma10) at ({-2.2*\lenX+1}, 0.4*\lenY) {$%
						\begin{array}{c}
							\textcolor{rougeTerre}{\Sigma_{1, 0}}\\[-3pt] = \\[-3pt]\Sigma^\diese_{1, 0}
						\end{array}$};%
					\node (ptSigma11) at (2.2*\lenX, 0.4*\lenY) {$%
						\begin{array}{c}
							\textcolor{rougeTerre}{\Sigma_{1, 1}}\\[-3pt] = \\[-3pt]\Sigma^\diese_{1, 1}
						\end{array}$};%
					\node[rougeTerre] (ptSigma20) at (-0.5*\lenX, -0.125*\lenY) {$\Sigma_{2, 0}$};
					\node (ptSigmaSharp20) at (0.5, -0.2*\lenY) {$\Sigma^\diese_{2, 0}$};
					\draw[latex-] (0, 0.4*\lenY) -- (ptSigma10);
					\draw[latex-] (1, 0.4*\lenY) -- (ptSigma11);
					\draw[latex-, rougeTerre] (-0.5*\lenX, 0) to (ptSigma20);
					\draw[latex-] (0.5, 0) to (ptSigmaSharp20);

				\end{scope}
			\end{tikzpicture}
		\end{subfigure}
		\hfill
		\begin{subfigure}[t]{0.64\textwidth}
			\centering
			\caption{$\ord = 3$}
			\def\lenX{0.75}
			\def\lenY{0.75}
			\def\lenZ{5.0}
			\def\hppXmax{1.75*\lenX}\def\hppXmin{-1.75*\lenX}\def\hppXlen{\hppXmax-\hppXmin}
			\def\hppYmax{1.75*\lenY}\def\hppYmin{-1.75*\lenY}\def\hppYlen{\hppYmax-\hppYmin}
			\def\hppZmax{0.85*\lenZ}\def\hppZmin{0}\def\hppZlen{\hppZmax-\hppZmin}
			\tdplotsetmaincoords{60}{58}
			\begin{tikzpicture}[tdplot_main_coords, scale=1.,]
				\begin{scope}
					\node (ptSigma301) at (0.5*\hppXmax, 0.5*\hppYmin, 0) {};
					\node[rougeTerre] (ptSigma302) at (0.5*\hppXmax, -0.5+1.5*\hppYmin, 0.5*\hppZmin) {$\Sigma_{3, 0}$};
					\node (ptSigma101) at (0, 0.25+0.5*\hppYmax, 0.5*\lenZ) {};
					\node[rougeTerre] (ptSigma102) at (0, 0.5+1.5*\hppYmax, 0.5*\lenZ) {$\Sigma_{1, 0}$};
					\node (ptSigma201) at (0.5*\hppXmin, 0, 0.25+0.75*\lenZ) {};
					\node[rougeTerre] (ptSigma202) at (0.5*\hppXmin, 0.75*\hppYmin, 0.25+1.0*\lenZ) {$\Sigma_{2, 0}$};
					\draw[latex-, rougeTerre] (ptSigma101) to[bend left=30] (ptSigma102);
					\draw[latex-, rougeTerre] (ptSigma201) to[bend left=30] (ptSigma202);
					\draw[latex-, rougeTerre] (ptSigma301) to[bend left=30] (ptSigma302);
					%
					%
					%
					\fill[fill=white] (0, \hppYmin, 0) -- (0, \hppYmin, \hppZmax) -- (0, \hppYmax, \hppZmax) -- (0, \hppYmax, 0) -- cycle;
					\fill[fill=white] (\hppXmin, \hppYmin, 0) -- (\hppXmin, \hppYmax, 0) -- (\hppXmax, \hppYmax, 0) -- (\hppXmax, \hppYmin, 0) -- cycle;
					\fill[fill=rougeTerre, fill opacity=0.3] (0, \hppYmin, \hppZmin) -- (0, \hppYmin, \hppZmax) -- (0, \hppYmax, \hppZmax) -- (0, \hppYmax, \hppZmin) -- cycle;
					\fill[fill=rougeTerre, fill opacity=0.3] (\hppXmin, \hppYmin, 0) -- (\hppXmin, \hppYmax, 0) -- (\hppXmax, \hppYmax, 0) -- (\hppXmax, \hppYmin, 0) -- cycle;
					%
					%
					\node (xaxis) at (2.75*\lenX, 0, 0) {$y_1$};
					\node (yaxis) at (0, 2.5*\lenY, 0) {$y_2$};
					\node (zaxis) at (0, 0, \lenZ) {$y_3$};

					\draw[-latex] (-2.25*\lenX, 0, 0) -- (xaxis);
					\draw[-latex] (0, -2.25*\lenY, 0) -- (yaxis);
					\draw[-latex] (0, 0, -0.3*\lenZ) -- (zaxis);
					\fill[fill=rougeTerre, fill opacity=0.2] (\hppXmin, 0, \hppZmin) -- (\hppXmin, 0, \hppZmax) -- (\hppXmax, 0, \hppZmax) -- (\hppXmax, 0, \hppZmin) -- cycle;
					%
					\draw[latex-, rougeTerre, opacity=0.5, dashed] (ptSigma301) to[bend left=30] (ptSigma302);
					\draw[latex-, rougeTerre, opacity=0.5, dashed] (ptSigma101) to[bend left=30] (ptSigma102);
					\draw[latex-, rougeTerre, opacity=0.5, dashed] (ptSigma201) to[bend left=30] (ptSigma202);
				\end{scope}
			\end{tikzpicture}
			\hfill
			\def\lenX{2.5}
			\def\lenY{1.0}
			\def\hppXmax{0.75*\lenX}\def\hppXmin{-1.25}\def\hppXlen{\hppXmax-\hppXmin}
			\def\hppYmax{1.50*\lenY}\def\hppYmin{-0.5}\def\hppYlen{\hppYmax-\hppYmin}
			\begin{tikzpicture}[tdplot_main_coords, scale=1.,]
				\begin{scope}
					\node (ptOmegaSharp1) at (0.5, 0.5, 0.75*\hppZmax) {};
					\node (ptOmegaSharp2) at (0.5, 0.5, 1.25*\hppZmax) {$\Omega^\diese$};
					\node (ptSigma20sharp) at (0.5, -1.75, 1.5) {$\Sigma^\diese_{2, 0}$};
					\node (ptSigma21sharp) at (0.5, 2, 0.9*\lenZ) {$\Sigma^\diese_{2, 1}$};
					\node (ptSigma10sharp) at (0.5, -1.25, 3.5) {$\Sigma^\diese_{1, 0}$};
					\node (ptSigma11sharp) at (1.5*\lenX, 0.5, 0.8*\lenZ) {$\Sigma^\diese_{1, 1}$};
					\node (ptSigma30sharp) at (0.5, -1.25, -0.5) {$\Sigma^\diese_{3, 0}$};
					\draw[latex-] (ptOmegaSharp1) to (ptOmegaSharp2);
					\draw[latex-] (0.5, 1, 0.6*\lenZ) to[bend right=20] (ptSigma21sharp);
					\draw[latex-] (0, 0.5, 0.6*\lenZ) to[bend right=50] (ptSigma10sharp);
					\draw[latex-] (0.5, 0.5, 0) to[bend left=50] (ptSigma30sharp);
					\node (xaxis) at (1.1*\lenX, 0, 0) {$y_1$};
					\node (yaxis) at (0, 2.25*\lenY, 0) {$y_2$};
					\node (zaxis) at (0, 0, \lenZ) {$y_3$};

					\draw[-latex] (-0.75, 0, 0) -- (xaxis);
					\draw[-latex] (0, -0.75, 0) -- (yaxis);
					\draw[-latex] (0, 0, -0.3*\lenZ) -- (zaxis);
					\draw[dashed, black!80] (0, 1, 0) -- (0, 1, \hppZmax); \draw (1, 1, 0) -- (1, 1, \hppZmax);
					\draw[dashed, black!80] (0, 0, 0) -- (0, 1, 0); \draw (1, 0, 0) -- (1, 1, 0);
					\draw[dashed, black!80] (0, 1, 0) -- (1, 1, 0);
					%
					\fill[black!30, fill opacity=0.8] (0, 0, 0) -- (0, 1, 0) -- (0, 1, \hppZmax) -- (0, 0, \hppZmax) -- cycle;
					\fill[black!30, fill opacity=0.8] (0, 1, 0) -- (1, 1, 0) -- (1, 1, \hppZmax) -- (0, 1, \hppZmax) -- cycle;
					\fill[black!30, fill opacity=0.8] (0, 0, 0) -- (1, 0, 0) -- (1, 0, \hppZmax) -- (0, 0, \hppZmax) -- cycle;
					\fill[black!30, fill opacity=0.8] (0, 0, 0) -- (1, 0, 0) -- (1, 1, 0) -- (0, 1, 0) -- cycle;
					\fill[black!20, fill opacity=0.8] (0, 0, \hppZmax) -- (1, 0, \hppZmax) -- (1, 1, \hppZmax) -- (0, 1, \hppZmax) -- cycle;
					\fill[black!30, fill opacity=0.5] (1, 0, 0) -- (1, 1, 0) -- (1, 1, \hppZmax) -- (1, 0, \hppZmax) -- cycle;
					\draw[dashed, black!80] (0, 1, \hppZmax) -- (0, 1, \lenZ); \draw[dashed, black!80] (1, 1, \hppZmax) -- (1, 1, \lenZ); \draw[dashed, black!80] (1, 0, \hppZmax) -- (1, 0, \lenZ);
					\draw[black!80] (0, 0, 0) -- (1, 0, 0);
					\draw[black!80] (1, 0, 0) -- (1, 1, 0);
					\draw[black!50] (0, 0, \hppZmax) -- (1, 0, \hppZmax); \draw[black!50] (0, 1, \hppZmax) -- (1, 1, \hppZmax);
					\draw[black!80] (0, 0, 0) -- (0, 0, \hppZmax);
					\draw[black!80] (1, 0, 0) -- (1, 0, \hppZmax);
					\draw[black!80] (1, 1, 0) -- (1, 1, \hppZmax);
					\draw[black!50] (0, 0, \hppZmax) -- (0, 1, \hppZmax); \draw[black!50] (1, 0, \hppZmax) -- (1, 1, \hppZmax);

					%
					%
					\draw[latex-, opacity=0.5, dashed] (ptOmegaSharp1) to (ptOmegaSharp2);
					\draw[latex-] (0.5, 0, 1.5) to (ptSigma20sharp);
					\draw[latex-, opacity=0.5, dashed] (0, 0.5, 0.6*\lenZ) to[bend right=50] (ptSigma10sharp);
					\draw[latex-] (1, 0.5, 0.4*\lenZ) to[bend right=20] (ptSigma11sharp);

					\draw[latex-, opacity=0.3, dashed] (0.5, 1, 0.6*\lenZ) to[bend right=20] (ptSigma21sharp);
				\end{scope}
			\end{tikzpicture}
		\end{subfigure}
	  \caption{Domains $\Omega^\diese$, $\Sigma_{i, a}$ and $\Sigma^\diese_{i, a}$ for $\ord = 2$ (a) and $\ord = 3$ (b).\label{fig:domains_3D}}
	\end{figure}

 For any open set $\mathcal{O} \subset \R^\ord$, let us first define the directional Sobolev space
	\begin{equation}
		 H^1_\cut(\mathcal{O}) := \big\{U \in L^2(\mathcal{O})\ /\ \Dt{} U \in L^2(\mathcal{O}) \big\},
	\end{equation}
	which is a Hilbert space, provided with the scalar product
	\[
	(U, V)_{H^1_\cut(\mathcal{O})} := \int_{\mathcal{O}} \Big(\Dt{} U\, \Dt{} \overline{V} + U\, \overline{V}\Big).
	\]
	Let us denote $\|\cdot\|_{H^1_\cut(\mathcal{O})}$ the induced norm. We begin with the following density property, whose proof can be found in \cite[Appendix 1]{temam1968stabilite}.
	\begin{lem}
	 \label{prop:density_Cinfy_H1cut}
	 The space $\mathscr{C}^\infty_0(\overline{\mathcal{O}})$ is dense in $H^1_\cut(\mathcal{O})$.
	\end{lem}
	\vspace{1\baselineskip}
	\noindent
	We denote the half-space $\R^\ord_+:=\{\yv\in\R^\ord,\;\yi_n>0 \}$  and the half-cylinder $\Omega^\diese:=(0,1)^{\ord-1}\times\R^+$ in the following. Let us introduce also the sets, for $a \in \{0, 1\}$ and for any integer $i\in\llbracket 1, \ord\rrbracket $,
	\[
		\Sigma_{i,a} = \{\yv \in \R^\ord_+,\ \yi_i = a\}\quad \text{and}\quad  \Sigma_{i,a}^\diese = \{\yv \in \Sigma_{i,a},\  \yi_j \in(0,1),\; j\in \llbracket 1, \ord-1\rrbracket,\; j\neq i\}.
	\]
	This definition is illustrated in Figure \ref{fig:domains_3D}. Note that $\Sigma_{\ord,a}^\diese$ is bounded whereas $\Sigma_{i,a}^\diese$ for $i\neq \ord$ is unbounded in the direction $\yi_\ord$. Moreover,
	\[
		\partial\Omega^\diese =  \Sigma_{n,0}^\diese \cup \bigg[ \bigcup_{i=1}^{n-1} \big(\overline{\Sigma}_{i,0}^\diese\cup\overline{\Sigma}_{i,1}^\diese\big)  \bigg].
	\]

	\noindent
	A trace operator can be defined from $H^1_{\cut}(\R^n_+)$ on $\Sigma_{i,a}$. The main idea for doing so consists in using a one-dimensional trace theorem on the $\cut$--oriented line that starts from a point $(\zi_1,\dots,\zi_{i-1}, a, \zi_{i+1},\dots,\zi_\ord) \in \Sigma_{i, a}$, to obtain an inequality which will be integrated with respect to $\zi_j$, $j \neq i$. The 1D trace theorem which will be used is the following.
	\begin{prop}
		\label{lem:1Dtrace_theorem}
		Let \textcolor{surligneur}{$L \in \R^*_+ \cup \{+\infty\}$}. Then the mapping $\gamma_L : u \mapsto u(0)$ is continuous from $H^1(0, L)$ to $\C$. Moreover, the operator norm of $\gamma_L$ is given by\color{surligneur}{%
		\begin{equation}
			\displaystyle
			\|\gamma_L\|^2 = \frac{\euler^{L} + \euler^{-L}}{\euler^{L} - \euler^{-L}} =: [\tanh L]^{-1}\ \ \textnormal{for}\ L > 0, \quad \textnormal{and} \quad \|\gamma_\infty\|^2 = 1.
		\end{equation}
		}
	\end{prop}
	\begin{dem}
		The continuity property is a classical result which can be proved by density.

		By definition, $\|\gamma_L\| := \sup\{|u(0)|,\ \|u\|_{H^1(0, L)} = 1\}$. This corresponds to a constrained optimization problem. Using the standard theory, this leads to introduce a Lagrange multiplier $\lambda$ and to find a pair $(\lambda, u_L) \in \C\setminus\{0\}\times H^1(0, L) $ such that $\|u_L\|_{H^1(0, L)} = 1$ and
		\begin{equation}
			\spforall v \in H^1(0, L) \quad \lambda\, u_L(0)\, \overline{v(0)} = \int_0^L \Big( \frac{d u_L}{d \xv}\,\frac{d \overline{v}}{d \xv} + u_L\, \overline{v} \Big) \; d\xv,
		\end{equation}
		in which case, we have $\|\gamma_L\|^2 = \lambda$. The explicit solution of this problem leads to the result.
	\end{dem}
	\noindent
	Note that, in particular, $\smash{\displaystyle \|\gamma_L\|^2 \underset{L \to 0}{\sim} L^{-1}}$.

	\vspace{1\baselineskip}\noindent
	We are now able to define traces on $\Sigma_{i, a}$ in the following sense.
	\begin{prop}
		\label{prop:trace_H1cut_demi_espace}
		Fix $a \in \{0, 1\}$ and $i\in \llbracket 1, \ord\rrbracket $. The mapping $\gamma_{i,a} : \mathscr{C}^\infty_0(\overline{\R^\ord_+}) \to \mathscr{C}^\infty_0(\Sigma_{i,a})$ defined by $\gamma_{i,a} U = \restr{U}{\Sigma_{i,a}}$ extends by continuity to a linear mapping still denoted $\gamma_{i,a}$, from $H^1_{\cut}(\R^{\ord}_+)$ to $L^2(\Sigma_{i,a})$, and which satisfies the estimate
		\begin{equation}
			\displaystyle
			\spforall U \in H^1_{\cut}(\R^{\ord}_+),\quad \|\gamma_{i,a} U\|^2_{L^2(\Sigma_{i,a})} \leq \frac{1}{\cuti_i}\; \|U\|^2_{H^1_{\cut}(\R^{\ord}_+)}.\label{eq:trace_H1cut_demi_espace}
		\end{equation}
	\end{prop}
	\begin{dem}
		 One can simply prove the continuity estimate \eqref{eq:trace_H1cut_demi_espace} for any function $U \in \mathscr{C}^\infty_0(\overline{\R^\ord_+})$ and conclude using the density result of Proposition \ref{prop:density_Cinfy_H1cut}.

		 \vspace{1\baselineskip} \noindent
		 %
		 ($i$) \smash{\underline{\textit{Case} $i \in \llbracket 1, \ord-1 \rrbracket$}}: Without loss of generality, we set $i = 1$. 
		 Define
		 \begin{equation}\displaystyle
		 		\Gamma_{1, a} := \{\zv = (\zi_2, \dots, \zi_\ord),\quad (a, \zv) \in \Sigma_{1, a}\} \equiv \R^{\ord-1}_+, \quad \textnormal{where} \quad (a, \zv) = (a, \zi_2, \dots, \zi_\ord).
				\label{eq:Gamma_1a}
		 \end{equation}
		 For $U \in \mathscr{C}^\infty_0(\overline{\R^\ord_+})$ and given any $\zv = (\zi_2, \dots, \zi_\ord) \in \Gamma_{1, a}$, consider the function
		 \begin{equation}
			 \displaystyle
			 \spforall \xv > 0, \quad u_{\zv, \cut}(\xv) = U(\xv\,\cut + (a, \itbf{z})).
		 \end{equation}
		 As $u_{\zv, \cut}$ belongs to $H^1(\R^*_+)$, Lemma \ref{lem:1Dtrace_theorem} for $L = +\infty$ combined with an integration with respect to $\zv \in \Gamma_{1, a}$ leads to 
		 \begin{equation}
			 \displaystyle
			 \int_{\Gamma_{1, a}} |u_{\zv, \cut}(0)|^2\; d \itbf{z} \leq \int_{\Gamma_{1, a}} \|u_{\zv, \cut}\|^2_{H^1(\R^*_+)} d \itbf{z}.
			 \label{eq:preuve_trace_1}
		 \end{equation}
		 On the other hand, let us introduce the transformation 
		 \begin{equation}
			 \label{eq:preuve_trace_3}
			 \displaystyle
			 \operatorname{T} : \yv \mapsto \big((\yi_1-a)/\cuti_1, \yi_2 - (\yi_1-a)\, \cuti_2/\cuti_1, \cdots, \yi_\ord - (\yi_1 - a)\, \cuti_\ord/\cuti_1\big),
		 \end{equation}
		 which defines a $\mathscr{C}^1$--diffeomorphism with a Jacobian determinant $\det \mathbf{J}_{\operatorname{T}} = 1/\theta_1 \neq 0$. Since the inverse image $\{\operatorname{T}^{-1}(\xv, \zv),\ \zv \in \Gamma_{1, a},\; \xv > 0\}$ is nothing but the polyhedron
		 \[
		 	\mathcal{Q}_{1, a} := \{\yv \in \R^\ord_+,\ \yi_1 > a,\ \yi_\ord > (\yi_1 - a)\, \cuti_\ord/\cuti_1\} \subset \R^\ord_+,
		 \]
		 it follows from the chain rule and from the change of variables $\yv \mapsto \operatorname{T}\yv$ that
		 \begin{equation}
			 \displaystyle
			 \frac{d u_{\itbf{z}, \cut}}{d\xv}(\xv) = \Dt{} U(x\,\cut + (a, \itbf{z})) \quad \textnormal{and} \quad \int_{\Gamma_{1, a}} \|u_{\zv, \cut}\|^2_{H^1(\R^*_+)}\; d \itbf{z} = \frac{1}{\cuti_1} \, \|U\|^2_{H^1_\cut(\mathcal{Q}_{1, a})}.
			 \label{eq:preuve_trace_2}
		 \end{equation}
		 Finally, since $u_{\zv, \cut}(0) = U(a, z_2, \cdots, z_\ord)$, Equations \eqref{eq:preuve_trace_1} and \eqref{eq:preuve_trace_2} imply
		 \begin{equation}\label{eq:preuve_trace_4}
			 \|U\|^2_{L^2(\Sigma_{1, a})} \leq \frac{1}{\cuti_1} \, \|U\|^2_{H^1_\cut(\mathcal{Q}_{1, a})} \leq \frac{1}{\cuti_1} \, \|U\|^2_{H^1_\cut(\R^\ord_+)},
		 \end{equation}
		 which is exactly the desired estimate.

		 \vspace{1\baselineskip} \noindent
		 ($ii$) \smash{\underline{\textit{Case} $i = \ord$}}: starting from the function $u_{\zv, \cut}(\xv) := U(\xv\,\cut + (\zv, a))$ defined for $\xv > 0$ and for any $\zv = (\zi_1, \dots, \zi_{\ord-1})$ with $(\zv, a) \in \Sigma_{\ord, a}$, the proof uses the exact same arguments as above, except the inverse image under $\operatorname{T}$ becomes the whole half-space $\mathcal{Q}_{\ord, a} := \{\yv \in \R^\ord_+,\ \yi_\ord > a\}$.
	\end{dem}

	\vspace{1\baselineskip} \noindent
		 {The previous result does not hold in general for functions which are only $H^1_\cut$ in sub-domains of the half-space $\R^\ord_+$}. In particular when it comes to the half-cylinder $\Omega^\diese$, one is led to apply the one-dimensional trace theorem on segments that become smaller in the neighbourhood of the “corners”, \emph{i.e.} the intersections of two faces. To overcome this difficulty, let us consider the sets (see Figure \ref{fig:domains_3D_bis}) 
		\begin{equation}
			\label{eq:Sigma_ia_securite}
			\spforall 0 < b < 1/2, \quad \Sigma_{i, a}^{\diese, b} = \{\yv \in \Sigma^\diese_{i,a},\quad  \operatorname{dist}(\yv,\, \partial \Sigma^\diese_{i, a}) := \inf_{\zv \, \in\, \partial \Sigma^\diese_{i, a}} |\yv - \zv| > b\}.
		\end{equation}
		Using these domains, the traces on $\Sigma_{i,a}^\diese$ can be defined as locally integrable functions in the sense of the following proposition.
		\begin{figure}[H]
			\centering
			\begin{subfigure}[t]{0.9\textwidth}
				\centering
				\def\lenX{2.5}
				\def\lenY{1.0}
				\def\lenZ{5.0}
				\def\distancesecu{0.25}
				\def\hppXmax{0.75*\lenX}\def\hppXmin{-1.25}\def\hppXlen{\hppXmax-\hppXmin}
				\def\hppYmax{1.50*\lenY}\def\hppYmin{-0.5}\def\hppYlen{\hppYmax-\hppYmin}
				\def\hppZmax{0.85*\lenZ}\def\hppZmin{0}\def\hppZlen{\hppZmax-\hppZmin}
				\def\thetaX{0.15}
				\def\thetaY{0.35}
				\def\thetaZ{1}
				\pgfmathsetmacro\xhppZ{\hppZmax/\thetaZ}
				\pgfmathsetmacro\xlenZ{\lenZ/\thetaZ}

				\tdplotsetmaincoords{60}{58}

				\begin{tikzpicture}[tdplot_main_coords, scale=1.,]
					\begin{scope}
						\node (ptSigma21sharp) at (0.5, 2, 0.9*\lenZ) {\color{rougeTerre}$\Sigma^{\diese, b}_{2, 1}$};
						\node (ptSigma10sharp) at (0.5, -1.25, 3.5) {\color{rougeTerre}$\Sigma^{\diese, b}_{1, 0}$};
						\node (ptSigma30sharp) at (0.5, -1.25, -0.5) {\color{rougeTerre}$\Sigma^{\diese, b}_{3, 0}$};
						\draw[latex-, rougeTerre] (0.5, 1, 0.6*\lenZ) to[bend right=20] (ptSigma21sharp);
						\draw[latex-, rougeTerre] (0, 0.5, 0.6*\lenZ) to[bend right=50] (ptSigma10sharp);
						\draw[latex-, rougeTerre] (0.5, 0.5, 0) to[bend left=50] (ptSigma30sharp);
						\node (xaxis) at (1.1*\lenX, 0, 0) {$y_1$};
						\node (yaxis) at (0, 2.25*\lenY, 0) {$y_2$};
						\node (zaxis) at (0, 0, \lenZ) {$y_3$};

						\draw[-latex] (-0.75, 0, 0) -- (xaxis);
						\draw[-latex] (0, -0.75, 0) -- (yaxis);
						\draw[-latex] (0, 0, -0.3*\lenZ) -- (zaxis);
						\draw[dashed, black!80] (0, 1, 0) -- (0, 1, \hppZmax); \draw (1, 1, 0) -- (1, 1, \hppZmax);
						\draw[dashed, black!80] (0, 0, 0) -- (0, 1, 0); \draw (1, 0, 0) -- (1, 1, 0);
						\draw[dashed, black!80] (0, 1, 0) -- (1, 1, 0);
						%
						\fill[rougeTerre!50, fill opacity=0.8] (0, \distancesecu, \distancesecu) -- (0, 1-\distancesecu, \distancesecu) -- (0, 1-\distancesecu, \hppZmax) -- (0, \distancesecu, \hppZmax) -- cycle;
						\fill[rougeTerre!50, fill opacity=0.8] (\distancesecu, 1, \distancesecu) -- (1-\distancesecu, 1, \distancesecu) -- (1-\distancesecu, 1, \hppZmax) -- (\distancesecu, 1, \hppZmax) -- cycle;
						\fill[rougeTerre!50, fill opacity=0.8] (\distancesecu, \distancesecu, 0) -- (1-\distancesecu, \distancesecu, 0) -- (1-\distancesecu, 1-\distancesecu, 0) -- (\distancesecu, 1-\distancesecu, 0) -- cycle;
						%
						%
						\draw[dashed, black!80] (0, 1, \hppZmax) -- (0, 1, \lenZ); \draw[dashed, black!80] (1, 1, \hppZmax) -- (1, 1, \lenZ); \draw[dashed, black!80] (1, 0, \hppZmax) -- (1, 0, \lenZ);
						\draw[black!80] (0, 0, 0) -- (1, 0, 0);
						\draw[black!80] (1, 0, 0) -- (1, 1, 0);
						\draw[black!80] (0, 0, 0) -- (0, 0, \hppZmax);
						\draw[black!80] (1, 0, 0) -- (1, 0, \hppZmax);
						\draw[black!80] (1, 1, 0) -- (1, 1, \hppZmax);
						%
						%
						\draw[latex-, red, opacity=0.5, dashed] (0, 0.5, 0.6*\lenZ) to[bend right=50] (ptSigma10sharp);

						\draw[latex-, red, opacity=0.3, dashed] (0.5, 1, 0.6*\lenZ) to[bend right=20] (ptSigma21sharp);
						\draw[latex-, red, opacity=0.5, dashed] (0.5, 0.5, 0) to[bend left=50] (ptSigma30sharp);
						\draw[<->] (0, \distancesecu, 0.5*\hppZmax) -- (0, 0, 0.5*\hppZmax);
						\draw (-0.5, 0.5*\distancesecu, 0.5*\hppZmax) node {$b$};
					\end{scope}
				\end{tikzpicture}
				\hfill
				\begin{tikzpicture}[tdplot_main_coords, scale=1.,]
					\begin{scope}
						\node (ptOmegaSharp1) at (0.5, 0.5, 0.75*\hppZmax) {};
						\node (ptOmegaSharp2) at (0.5, 0.5, 1.25*\hppZmax) {$\Omega^\diese$};
						\node (ptTn) at (\lenX, 1, 0) {\color{rougeTerre}$T_n$};
						\draw[latex-] (ptOmegaSharp1) to (ptOmegaSharp2);
						\draw[latex-, rougeTerre!80] (1, 0.5, 0) to[bend right=20] (ptTn);
						\node (xaxis) at (1.1*\lenX, 0, 0) {$y_1$};
						\node (yaxis) at (0, 2.25*\lenY, 0) {$y_2$};
						\node (zaxis) at (0, 0, \lenZ) {$y_3$};

						\draw[-latex] (-0.75, 0, 0) -- (xaxis);
						\draw[-latex] (0, -0.75, 0) -- (yaxis);
						\draw[-latex] (0, 0, -0.3*\lenZ) -- (zaxis);
						\draw[dashed, rougeTerre!100] (0, 1, 0) -- (0, 1, \hppZmax); \draw (1, 1, 0) -- (1, 1, \hppZmax);
						\draw[dashed, rougeTerre!100] (0, 0, 0) -- (0, 1, 0); \draw (1, 0, 0) -- (1, 1, 0);
						\draw[dashed, rougeTerre!100] (0, 1, 0) -- (1, 1, 0);
						%
						\fill[black!30, fill opacity=0.8] (0, 0, 0) -- (0, 1, 0) -- (0, 1, \hppZmax) -- (0, 0, \hppZmax) -- cycle;
						\fill[black!30, fill opacity=0.8] (0, 1, 0) -- (1, 1, 0) -- (1, 1, \hppZmax) -- (0, 1, \hppZmax) -- cycle;
						\fill[black!30, fill opacity=0.8] (0, 0, 0) -- (1, 0, 0) -- (1, 0, \hppZmax) -- (0, 0, \hppZmax) -- cycle;
						\fill[rougeTerre!50, fill opacity=0.8] (0, 0, 0) -- (1, 0, 0) -- (1, 1, 0) -- (0, 1, 0) -- cycle;
						\fill[black!20, fill opacity=0.8] (0, 0, \hppZmax) -- (1, 0, \hppZmax) -- (1, 1, \hppZmax) -- (0, 1, \hppZmax) -- cycle;
						\fill[black!30, fill opacity=0.5] (1, 0, 0) -- (1, 1, 0) -- (1, 1, \hppZmax) -- (1, 0, \hppZmax) -- cycle;
						\fill[bleuDoux!30, fill opacity=0.4] (0, 0, 0) -- (0, 1, 0) -- (0, 1, 0.5*\hppZmax) -- (0, 0, 0.5*\hppZmax) -- cycle;
						\fill[bleuDoux!30, fill opacity=0.4] (0, 1, 0) -- (1, 1, 0) -- (1, 1, 0.5*\hppZmax) -- (0, 1, 0.5*\hppZmax) -- cycle;
						\fill[bleuDoux!30, fill opacity=0.4] (0, 0, 0) -- (1, 0, 0) -- (1, 0, 0.5*\hppZmax) -- (0, 0, 0.5*\hppZmax) -- cycle;
						\fill[bleuDoux!30, fill opacity=0.2] (1, 0, 0) -- (1, 1, 0) -- (1, 1, 0.5*\hppZmax) -- (1, 0, 0.5*\hppZmax) -- cycle;

						\fill[bleuDoux!50, fill opacity=0.8] (0, 0, 0.5*\hppZmax) -- (1, 0, 0.5*\hppZmax) -- (1, 1, 0.5*\hppZmax) -- (0, 1, 0.5*\hppZmax) -- cycle;

						\draw[dashed, rougeTerre] (0, 1, \hppZmax) -- (0, 1, \lenZ); \draw[dashed, rougeTerre] (1, 1, \hppZmax) -- (1, 1, \lenZ); \draw[dashed, rougeTerre] (1, 0, \hppZmax) -- (1, 0, \lenZ);
						\draw[rougeTerre, line width=0.5mm] (0, 0, 0) -- (1, 0, 0);
						\draw[rougeTerre, line width=0.5mm] (1, 0, 0) -- (1, 1, 0);
						\draw[black!50] (0, 0, \hppZmax) -- (1, 0, \hppZmax); \draw[black!50] (0, 1, \hppZmax) -- (1, 1, \hppZmax);
						\draw[rougeTerre, line width=0.5mm] (0, 0, 0) -- (0, 0, \hppZmax);
						\draw[rougeTerre, line width=0.5mm] (1, 0, 0) -- (1, 0, \hppZmax);
						\draw[rougeTerre, line width=0.5mm] (1, 1, 0) -- (1, 1, \hppZmax);
						\draw[black!50] (0, 0, \hppZmax) -- (0, 1, \hppZmax); \draw[black!50] (1, 0, \hppZmax) -- (1, 1, \hppZmax);

						%
						%
						\draw[latex-, opacity=0.5, dashed] (ptOmegaSharp1) to (ptOmegaSharp2);
						\draw (-0.25, -0.25, 0.5*\hppZmax) node {$a$};
						\node[bleuDoux] (ptOmegasharpa) at (0.5, 2, 0.4*\hppZmax) {$\Omega^\diese_{a, -}$};
						\draw[latex-, bleuDoux] (0.5, 1, 0.25*\hppZmax) to[bend right=20] (ptOmegasharpa);

					\end{scope}
				\end{tikzpicture}
				\hfill
				\begin{tikzpicture}[tdplot_main_coords, scale=1.,]
					\begin{scope}
						\coordinate (ptCut) at ({\xhppZ*\thetaX}, {\xhppZ*\thetaY}, {\xhppZ*\thetaZ});
						\coordinate (ptCut2) at ({\xlenZ*\thetaX}, {\xlenZ*\thetaY}, {\xlenZ*\thetaZ});

						\node (ptOmegaSharp1) at ($(1, 0.5, 0)+0.75*(ptCut)$) {};
						\node (ptOmegaSharp2) at ($(2.5, 0.5, 0)+1.25*(ptCut)$) {$\Omega^\diese_\cut$};


						\draw[latex-, bend right=40] (ptOmegaSharp1) to (ptOmegaSharp2);
						\node (xaxis) at (1.1*\lenX, 0, 0) {$y_1$};
						\node (yaxis) at (0, 2.25*\lenY, 0) {$y_2$};
						\node (zaxis) at (0, 0, \lenZ) {$y_3$};

						\draw[-latex] (-0.75, 0, 0) -- (xaxis);
						\draw[-latex] (0, -0.75, 0) -- (yaxis);
						\draw[-latex] (0, 0, -0.3*\lenZ) -- (zaxis);
						\draw[dashed, black!80] (0, 1, 0) -- ($(0, 1, 0) + (ptCut2)$);
						\draw (1, 1, 0) -- ($(1, 1, 0) + (ptCut)$); \draw[dashed] ($(1, 1, 0) + (ptCut)$) -- ($(1, 1, 0) + (ptCut2)$);
						\draw[dashed, black!80] (0, 0, 0) -- (0, 1, 0); \draw (1, 0, 0) -- (1, 1, 0);
						\draw[dashed, black!80] (0, 1, 0) -- (1, 1, 0);
						%
						\fill[black!30, fill opacity=0.8] (0, 0, 0) -- (0, 1, 0) -- ($(0, 1, 0) + (ptCut)$) -- (ptCut) -- cycle;
						\fill[black!30, fill opacity=0.8] (0, 1, 0) -- (1, 1, 0) -- ($(1, 1, 0) + (ptCut)$) -- ($(0, 1, 0) + (ptCut)$) -- cycle;
						\fill[black!30, fill opacity=0.8] (0, 0, 0) -- (1, 0, 0) -- ($(1, 0, 0) + (ptCut)$) -- (ptCut) -- cycle;
						\fill[black!30, fill opacity=0.8] (0, 0, 0) -- (1, 0, 0) -- (1, 1, 0) -- (0, 1, 0) -- cycle;
						\fill[black!20, fill opacity=0.8] (ptCut) -- ($(1, 0, 0) + (ptCut)$) -- ($(1, 1, 0) + (ptCut)$) -- ($(0, 1, 0) + (ptCut)$) -- cycle;
						\fill[black!30, fill opacity=0.5] (1, 0, 0) -- (1, 1, 0) -- ($(1, 1, 0) + (ptCut)$) -- ($(1, 0, 0) + (ptCut)$) -- cycle;
						%
						%
						\draw (0, 0, 0) -- (ptCut); \draw[dashed] (ptCut) -- (ptCut2);
						\draw (1, 0, 0) -- ($(1, 0, 0) + (ptCut)$); \draw[dashed] ($(1, 0, 0) + (ptCut)$) -- ($(1, 0, 0) + (ptCut2)$);
						\draw (0, 0, 0) -- (1, 0, 0);
						\draw (1, 0, 0) -- (1, 1, 0);
						\draw[black!50] (ptCut) -- ($(1, 0, 0) + (ptCut)$); \draw[black!50] ($(0, 1, 0) + (ptCut)$) -- ($(1, 1, 0) + (ptCut)$);
						\draw[black!50] (ptCut) -- ($(0, 1, 0) + (ptCut)$); \draw[black!50] ($(1, 0, 0) + (ptCut)$) -- ($(1, 1, 0) + (ptCut)$);
						\fill [red!50, opacity=0.2] (0, -0.0, 0) -- (0, \hppYmax, 0) -- (0, \hppYmax, \hppZmax) -- (0, -0.0, \hppZmax) -- cycle;
						\fill [red!50, opacity=0.2] (-0.0, 0, 0) -- (\hppXmax, 0, 0) -- (\hppXmax, 0, \hppZmax) -- (-0.0, 0, \hppZmax) -- cycle;
						%
						\draw[latex-, opacity=0.5, dashed, bend right=40] (ptOmegaSharp1) to (ptOmegaSharp2);
					\end{scope}
				\end{tikzpicture}
			\end{subfigure}
			\caption{
		  				 From left to right: $\Sigma^{\diese, b}_{i, a}$ \eqref{eq:Sigma_ia_securite}, $T_n$ \eqref{eq:definition_table}, $\Omega^\diese_{a, -}$ \eqref{eq:Omega_diese_a_plus}, and $\Omega^\diese_\cut$ \eqref{eq:Omegadiesetheta} represented for $\ord = 3$.%
							\label{fig:domains_3D_bis}}
		\end{figure}
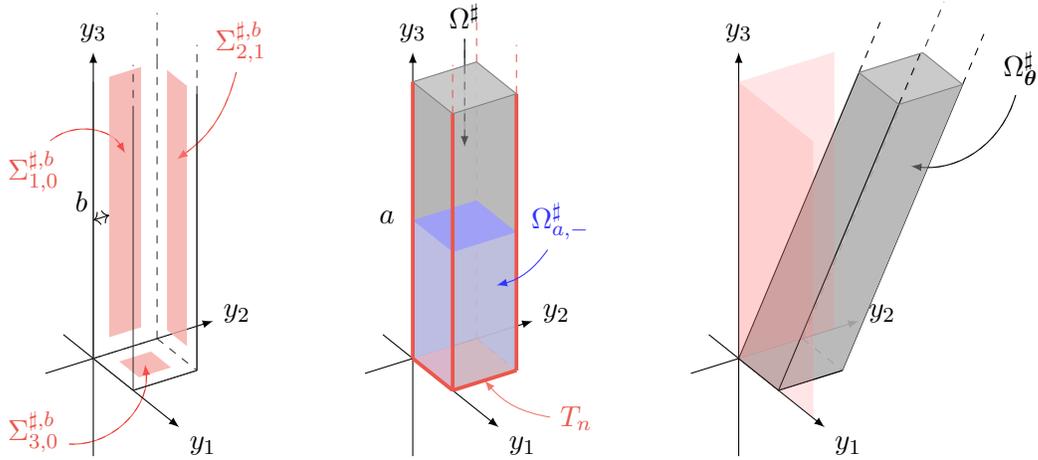
		\begin{prop}
		\label{prop:trace_L2loc_demi_cylindre}%
		Let $a \in \{0, 1\}$ and $i\in \llbracket 1, \ord\rrbracket $. The mapping $\gamma_{i,a}^\diese : \mathscr{C}^\infty_0(\overline{\Omega}^\diese) \to \mathscr{C}^\infty_0(\Sigma^\diese_{i,a})$ defined by $\gamma_{i, a}^\diese U = \restr{U}{\Sigma_{i, a}^\diese}$ extends by continuity to a linear mapping still denoted $\gamma_{i, a}^\diese $, from $H^1_{\cut}(\Omega^\diese)$ to $L^2_{\textit{loc}}(\Sigma_{i, a}^\diese)$, and which satisfies the estimate
	 	\begin{equation}
	 		\displaystyle
			\label{eq:trace_L2loc_demi_cylindre}
			\spforall  0 < b < 1/2,\quad \spexists C_b > 0, \quad \spforall U \in H^1_{\cut}(\Omega^\diese), \quad \|\gamma_{i, a}^\diese  U\|^2_{L^2(\Sigma_{i, a}^{\diese, b})} \leq \frac{C_b}{\cuti_i}\, \|U\|^2_{H^1_{\cut}(\Omega^\diese)}.
	 	\end{equation}
	 \end{prop}
	 \begin{dem}
		 Using the density result stated in Proposition \ref{prop:density_Cinfy_H1cut}, one only has to show \eqref{eq:trace_L2loc_demi_cylindre} for $U \in \mathscr{C}^\infty_0(\overline{\Omega}^\diese)$. Let us assume that $i = 1$ and $a = 0$, the arguments in the following extending without any difficulty to $i\in \llbracket 1, \ord \rrbracket$ and $a \in \{0, 1\}$. Define
		 \begin{equation}\displaystyle
		 		\Gamma^\diese_{1, 0} := \{\zv = (\zi_2, \dots, \zi_\ord),\quad (0, \zv) \in \Sigma^\diese_{1, 0}\} \equiv (0, 1)^{\ord-1} \times \R_+.
				\label{eq:GammaDiese_1a}
		 \end{equation}
		 We introduce the length function defined by 
		 \[\displaystyle
		 	 	\spforall \zv \in \Gamma^\diese_{1, 0}, \quad \lambda_{1, 0}(\zv) := \big| \{\cut\,\R + (0, \zv)\} \cap \Omega^\diese \big| = \sup\{\xv > 0, \;  \xv\,\cuti_1\leq 1,\ \xv\,\cuti_i + \zi_i \leq 1\ \ \spforall i\in\llbracket2, \ord-1\rrbracket  \}.
		 \]
		 We deduce easily that
		 \begin{equation}
			 \label{eq:expression_longueur}
				\lambda_{1, 0}(\zv) = \min\bigg\{ \frac{1}{\cuti_1};\; \min_{2 \leq j \leq \ord-1} \Big( \frac{1 - \zi_j}{\cuti_j} \Big) \bigg\}.
		 \end{equation}
		 %
		 For $U \in \mathscr{C}^\infty_0(\overline{\Omega}^\diese)$ and $\zv \in \Gamma^\diese_{1, 0}$, we define
		 \begin{equation}
			 \label{eq:preuve_trace_demi_cylindre_-1}
			 \displaystyle
			 \spforall 0 < \xv < \lambda_{1, 0}(\zv), \quad u_{\zv, \cut}(\xv) = U(\xv\,\cut + (0, \itbf{z})).
		 \end{equation}
		 Since $u_{\zv, \cut} \in \smash{H^1\big(0, \lambda_{1, 0}(\zv)\big)}$, {Lemma \ref{lem:1Dtrace_theorem}} and an integration with respect to $\zv$ give 
		 \begin{equation}
			 \displaystyle
			 \int_{\Gamma^\diese_{1, 0}} w_{1, 0}(\zv)\; |u_{\zv, \cut}(0)|^2\; d\zv \leq \int_{\Gamma^\diese_{1, 0}} \|u_{\zv, \cut}\|^2_{H^1(0, \gamma_{i, a}(\zv))} \; d\zv, \quad \textnormal{with} \  w_{1, 0}(\zv) = \tanh[ \lambda_{1, 0}(\zv)].
			 \label{eq:preuve_trace_demi_cylindre_0}
		 \end{equation}
		 On the other hand, consider the $\mathscr{C}^1$--diffeomorphism $\operatorname{T}$ given by \eqref{eq:preuve_trace_3}. The set $\mathcal{Q}^{\diese}_{1, 0} := \{\operatorname{T}^{-1}(\xv, \zv),\ \  0 < \xv < \lambda_{1, 0}(\zv),\ \zv \in \Gamma^{\diese}_{1, 0}\}$ is clearly included in $\Omega^\diese$. {Thus, by analogy with \eqref{eq:preuve_trace_4} in} the proof of Proposition \ref{prop:trace_H1cut_demi_espace}, we have from \eqref{eq:preuve_trace_demi_cylindre_-1}, the chain rule, and the change of variables $\yv \mapsto \operatorname{T}\yv$ that
		 \begin{equation}
			 \displaystyle
			 \int_{\Gamma^\diese_{1, 0}} w_{1, 0}(\zv)\; |U(0, \zv)|^2\; d\zv \leq \frac{1}{\cuti_1} \, \|U\|^2_{H^1_\cut(\Omega^\diese)}.
			 \label{eq:preuve_trace_demi_cylindre_1}
		 \end{equation}
		 More generally, we can show that $\gamma^\diese_{i, a}$ can be defined from $H^1_\cut(\Omega^\diese)$ to the weighted space $L^2(\Sigma^\diese_{i, a}, w_{i, a}\, d\zv)$, where the weight $w_{i, a}$ is given in \eqref{eq:preuve_trace_demi_cylindre_0} for $i = 1$ and $a = 0$. Now, the expression \eqref{eq:expression_longueur} of $\lambda_{1, 0}$ implies that $w_{1, 0}$ degenerates at the neighbourhood of the corners $\zi_j = 1$. However, the weight $w_{1, 0}$ is bounded from below on $\Sigma^{\diese, b}_{1, 0}$ with
		 \begin{equation}
			 \label{eq:preuve_trace_demi_cylindre_2}
			 \displaystyle
			  \inf_{(0, \zv) \in \Sigma^{\diese, b}_{1, 0}} w_{1, 0}(\zv) = \tanh \bigg[ \min\Big\{ \frac{1}{\cuti_1};\; b\min_{2 \leq j \leq \ord-1} \frac{1}{\cuti_j} \Big\} \bigg] > 0.
		 \end{equation}
		 If we set $\smash{C_b := [\inf_{(0, \zv) \in \Sigma^{\diese, b}_{1, 0}} w_{1, 0}(\zv)]^{-1} > 0}$, then \eqref{eq:trace_L2loc_demi_cylindre} follows directly from \eqref{eq:preuve_trace_demi_cylindre_1} by integrating with respect to $\{\zv,\ (0, \zv) \in \Sigma^{\diese, b}_{1, 0}\}$, instead of $\Gamma^\diese_{1, 0}$.
	 \end{dem}
	 \color{black}%
	\begin{rmk}\label{rem:traceH1theta}
		The best constant in the previous proposition necessarily blows up when $b$ tends to 0. The above proof shows that traces could be defined on the whole faces in appropriate weighted $L^2$-spaces. More details about traces in anisotropic spaces can be found in \cite{joly1992some}.
	\end{rmk}

	\subsubsection{Green's formulas}\label{sub:Green}
	{Let us now introduce} the set $H^1_{\cut, \textit{loc}}(\R^n_+)$ of functions which are \textcolor{surligneur}{$H^1_\cut$ in any half-cylinder $S \times \R_+$ where $S$ is a bounded open set in $\R^{\ord-1}$}. More rigorously, we define for any $\varphi \in \mathscr{C}^\infty_0(\R^{n-1})$ the $\ord$--dimensional function $\check{\varphi} \in \mathscr{C}^\infty(\R^{\ord})$
	such that
	\begin{equation}
		\label{eq:cut_off_etendu_Rn}
		\check{\varphi}(y_1,\ldots,y_{n-1},y_n)= \varphi(y_1,\ldots,y_{n-1}).
	\end{equation}
	Note that for any $U\in L^2_\textit{loc}(\R^n_+)$, the support of $\check{\varphi}\, U$ is bounded in the directions $y_j,\,j\neq n$. Starting from this remark, we define
	\begin{equation}
		\begin{array}{c}
			H^1_{\cut, \textit{loc}}(\R^n_+) := \Big\{U \in L^2_{\textit{loc}}(\R^n_+), \quad \check{\varphi}\, U\; \in H^1_\cut(\R^+_n)\ \ \spforall\!\varphi \in \mathscr{C}^\infty_0(\R^{n-1}) \Big\}.
		\end{array}
	\end{equation}
	%
	%

	%
	%
	%
	%
	\noindent
	%
	%
	%
	%
	%
	%
	%
	Let us introduce a 1D cut-off function $\chi\in \mathscr{C}^\infty_0(\R)$ such that $\chi=1$ on $(0,1)$, from which we define $\check{\chi}_\diese\in \mathscr{C}^\infty_0(\R^\ord)$ as
	\begin{equation}
		\label{eq:definition_cutoff}
		\check{\chi}_\diese(y_1,\ldots,y_{n-1}, y_\ord)=\chi(y_1)\ldots\chi(y_{n-1}).
	\end{equation}
	We deduce in particular that
	\begin{equation}\label{eq:restr_halfspace}
		\spforall U\in H^1_{\cut, \textit{loc}}(\R^{\ord}_+),\quad \restr{U}{\Omega^\diese}=\restr{(\check{\chi}_\diese\, U)}{\Omega^\diese}\;\in\;H^1_\cut({\Omega^\diese}).
	\end{equation}
	Moreover, by Proposition \ref{prop:trace_H1cut_demi_espace}, it is obvious that we can define without any ambiguity the trace map $\gamma_{i,a}^\diese$ to $H^1_{\cut, \textit{loc}}(\R^{\ord}_+)$ as follows
	\begin{equation}
		\label{eq:trace_H1loc}
		\spforall U\in H^1_{\cut, \textit{loc}}(\R^{\ord}_+),\quad \gamma_{i,a}^\diese U:=\restr{\gamma_{i,a} (\check{\chi}_\diese U)}{\Sigma_{i,a}^\diese}\;\in L^2(\Sigma_{i,a}^\diese).
	\end{equation}
	For simplicity, when considering traces on $\Sigma_{i,a}^\diese$, we  shall write $U$ instead of $\gamma_{i,a}^\diese U$.
	We can now state the following Green's formula.
	\begin{prop}\label{prop:Green_formula}
	For any $U, V \in H^1_{\cut, \textit{loc}}(\R^{\ord}_+)$, we have the Green's formula
	\begin{equation}
		\int_{\Omega^\diese} \left( \Dt{} U \; \overline{V} + U\; \Dt{} \overline{V} \right)\; d \yv = \frac{1}{\cuti_\ord} \int_{\Sigma^\diese_{n,0}} U \; \overline{V} \;d s + \sum_{i = 1}^{\ord-1} \frac{1}{\cuti_i} \Big( \int_{\Sigma^\diese_{i,1}} U \; \overline{V} \;d s - \int_{\Sigma^\diese_{i,0}} U \; \overline{V} \;d s \Big). \label{eq:Green_formula}
	\end{equation}
\end{prop}
	 \begin{dem}
		 Let $U, V \in H^1_{\cut, \textit{loc}}(\R^{\ord}_+)$. By definition, for any $\chi \in \mathscr{C}^\infty_0(\R)$ such that $\chi = 1$ on $(0, 1)$, the functions $\check{\chi}_\diese\, U$ and $\check{\chi}_\diese\, V$ belong to $H^1_\cut(\R^\ord_+)$, where $\check{\chi}_\diese$ is defined in \eqref{eq:definition_cutoff}.
		 Since Proposition \ref{prop:density_Cinfy_H1cut} ensures that $\mathscr{C}^\infty_0(\overline{\R^\ord_+})$ is dense in $H^1_\cut(\R^\ord_+)$, there exist two sequences $(U_k)_{k \in \N}, (V_k)_{k \in \N}$ of functions in $\mathscr{C}^\infty_0(\overline{\R^\ord_+})$, such that
		 \[\displaystyle
		 		U_k \to \check{\chi}_\diese\, U \quad \textnormal{and} \quad V_k \to \check{\chi}_\diese\, V \quad \textnormal{in} \quad H^1_\cut(\R^\ord_+), \quad k \to +\infty.
		 \]
		 It follows from Green's formula for smooth functions that $U_k$ and $V_k$ satisfy \eqref{eq:Green_formula} for any $k \in \N$. Passing to the limit and using the trace continuity result stated in Propsition \ref{prop:trace_H1cut_demi_espace} imply that \eqref{eq:Green_formula} is satisfied by $\check{\chi}_\diese\, U$ and $\check{\chi}_\diese\, V$, \emph{i.e.} by $U$ and $V$, since $\check{\chi}_\diese = 1$ in $\Omega^\diese$.
		 %
	 \end{dem}

	 \vspace{1\baselineskip} \noindent
	 We next focus on functions which are periodic with respect to their $(n-1)$ first variables. More precisely,
	 for any $U \in L^2(\Omega^\diese)$ and any $\varphi \in L^2(\Sigma^\diese_{\ord, 0})$, we introduce the respective periodic extensions $\widetilde{U} \in L^2_{\textit{loc}}(\R^\ord_+)$ and $\widetilde{\varphi} \in L^2_{\textit{loc}} (\Sigma_{\ord, 0})$ as defined {for any $i\in \llbracket 1, \ord-1\rrbracket $} by
	\begin{equation}\label{eq:per_extension}
		\left\{
		\begin{array}{l@{\quad}l@{\quad \textnormal{and} \quad}l}
			\aeforall \yv \in \R^n_+, & \widetilde{U}(\yv + \vec{\ev}_i) = \widetilde{U}(\yv) & \restr{\widetilde{U}}{\Omega^\diese}= U.
			\\[8pt]
			 \aeforall \sv \in \Sigma_{\ord, 0}, &  \widetilde{\varphi}(\sv + \vec{\ev}_i) = \widetilde{\varphi}(\sv) &  \restr{\widetilde{\varphi}}{\Sigma^\diese_{\ord, 0}}= \varphi.
		\end{array}
		\right.
	\end{equation}
	An appropriate functional framework is provided by the space
	\begin{equation}
		\label{eq:H1thetaper}
		\displaystyle
		H^1_{\cut, \textit{per}}(\Omega^\diese) = \Big\{ U \in L^2(\Omega^\diese),\  \widetilde{U} \in H^1_{\cut, \textit{loc}}(\R^n_+) \Big\} \ \subset H^1_{\cut}(\Omega^\diese),
	\end{equation}
	where the inclusion follows from \eqref{eq:restr_halfspace} and \eqref{eq:per_extension}.
	If $\mathscr{C}^\infty_{\textit{per}}(\Omega^\diese)$ denotes the set of smooth functions in $\mathscr{C}^\infty(\Omega^\diese)$ which are $1$--periodic with respect to their first $\ord-1$ variables, that is,
	\begin{equation}
		\displaystyle
		\label{eq:Cinfyper0}
		\mathscr{C}^\infty_{\textit{per}}({\Omega}^\diese) = \Big\{V \in \mathscr{C}^\infty(\Omega^\diese),\quad \widetilde{V} \in \mathscr{C}^\infty(\R^\ord_+) \Big\},
	\end{equation}
	then one can show the following result by adapting classical properties of $H^1$ functions.
	\begin{lem}
		\label{prop:density_Cinfyper_H1cutper}
		The space $\mathscr{C}^\infty_{\textit{per}}({\Omega}^\diese)$ is dense in $H^1_{\cut, \textit{per}}(\Omega^\diese)$.
	\end{lem}%

	\noindent
	Note that the traces of functions in $H^1_{\cut, \textit{per}}(\Omega^\diese)$ on $\Sigma_{i,a}^\diese$ are well-defined in $L^2$ by \eqref{eq:trace_H1loc}. Moreover, using the continuity estimate \eqref{eq:trace_H1cut_demi_espace} we have
\begin{equation}
	\label{eq:trace_func_per}
	\gamma_{i,a}^\diese\in\mathcal{L}(H^1_{\cut, \textit{per}}(\Omega^\diese),L^2(\Sigma_{i,a}^\diese)).
\end{equation}
	{One has the characterization}
	\begin{equation}\label{eq:charac_H1per}
		\displaystyle
		H^1_{\cut, \textit{per}}(\Omega^\diese) = \Big\{ U \in H^1_\cut(\Omega^\diese),\quad \gamma_{i,0}^\diese U =\gamma_{i,1}^\diese U\ \  \spforall i\in \llbracket 1, \ord-1\rrbracket   \Big\},
	\end{equation}
	where {the traces of functions in $ H^1_\cut(\Omega^\diese)$ are defined in Proposition \ref{prop:trace_L2loc_demi_cylindre}} and the equality of traces has to be understood up to the identification of functions on $\Sigma_{i,0}^\diese$ and	$\Sigma_{i,1}^\diese$. It is clear from \eqref{eq:charac_H1per} that $H^1_{\cut, \textit{per}}(\Omega^\diese)$ is a closed subspace of $H^1_{\cut}(\Omega^\diese)$, thus it is an Hilbert space when equipped with the norm of $H^1_{\cut}(\Omega^\diese)$. %
	From Proposition \ref{prop:Green_formula} and \eqref{eq:charac_H1per}, we deduce the Green's formula on $H^1_{\cut, \textit{per}}(\Omega^\diese)$.
	\begin{prop}\label{prop:Green_formula_H1per}
		For any $U, V \in H^1_{\cut, \textit{per}}(\Omega^\diese)$, we have the Green's formula
		\begin{equation}
			\int_{\Omega^\diese} \left( \Dt{} U \; \overline{V} + U\; \Dt{} \overline{V}\right) d \yv = \frac{1}{\cuti_\ord} \int_{\Sigma^\diese_{n,0}} U \; \overline{V} \;d s. \label{eq:Green_formula_H1per}
	   \end{equation}
	\end{prop}

	\noindent
	From the Green's formula \eqref{eq:Green_formula_H1per}, we can easily deduce the following result.
	\begin{cor}\label{cor:rest_H1thetaper}
		Let $a>0$, and define the sets with common boundary $\Sigma^\diese_{\ord, a}$ (see Figure \ref{fig:domains_3D_bis}):
		\begin{equation}\label{eq:Omega_diese_a_plus}
			\Omega^\diese_{a,+} := \Omega^\diese \cap \{\yi_\ord > a\} \quad \textnormal{and} \quad \Omega^\diese_{a,-} := \Omega^\diese \cap \{\yi_\ord < a\}.
		\end{equation}
		Consider a function $U \in L^2(\Omega^\diese)$ such that $U_\pm := U|_{\Omega^\diese_{a,\pm}} \in H^1_{\cut, \textit{per}}(\Omega^\diese_{a, \pm})$, where $H^1_{\cut, \textit{per}}(\Omega^\diese_{a,\pm})$ is defined as in \eqref{eq:charac_H1per}. Then
		\[
			U\in H^1_{\cut, \textit{per}}(\Omega^\diese)\quad \Longleftrightarrow \quad \gamma_{n,a}^\diese U_+ =\gamma_{n,a}^\diese U_-.
		\]
	\end{cor}
	 \noindent%
	 We finish this section with a more technical Green's formula, used in the proof of Proposition \ref{prop:half_guide_FV}, involving functions $U$ that only belong to $H^1_\cut(\Omega^\diese)$, provided that the test function $V$ vanishes in the neighborhood of the skeleton $T_\ord$ defined by
	 \begin{equation}
		 \label{eq:definition_table}
		 \displaystyle
		 T_2 = \overline{\Sigma}^\diese_{2, 0} \quad \textnormal{and} \quad T_\ord = \overline{\Sigma}^\diese_{\ord, 0} \cup \bigg[ \bigcup_{j = 1}^{\ord-1} \big(\partial \Sigma^\diese_{j, 0} \cup \partial \Sigma^\diese_{j, 1} \big) \bigg] \quad \textnormal{for $\ord \geq 3$}.
	 \end{equation}
	 This domain is represented in Figure \ref{fig:domains_3D_bis} for $\ord = 3$.
	 \begin{prop}\label{prop:formule_Green_demi_cylindre}
 			For $U \in H^1_{\cut}(\Omega^\diese)$ and $V \in \mathscr{C}^\infty_0(\overline{\Omega}^\diese \setminus T_\ord)$, the Green's formula \eqref{eq:Green_formula} still holds.
	\end{prop}
	\begin{dem}
		Consider $U \in H^1_{\cut}(\Omega^\diese)$ and $V\in \mathscr{C}^\infty_0(\overline{\Omega}^\diese \setminus T_\ord)$. Since by Proposition \ref{prop:density_Cinfy_H1cut}, $\mathscr{C}^\infty_0(\overline{\Omega}^\diese)$ is dense in $H^1_{\cut}(\Omega^\diese)$, there exists a sequence $(U_k)_{k \in \N}$ of functions in $\mathscr{C}^\infty_0(\overline{\Omega}^\diese)$ which tends to $U$. It follows from Green's formula in $\Omega^\diese$ for smooth functions that $U_k$ and $V$ satisfy \eqref{eq:Green_formula} for any $k \in \N$. For $0 < b < 1/2$, let $\Omega^{\diese, b}$ be the domain
		\begin{equation}
			\displaystyle
			\Omega^{\diese, b} = \{\yv \in \Omega^\diese,\quad  \operatorname{dist}(\yv,\, T_n) := \inf_{\zv \, \in\, T_n} |\yv - \zv| > b\}.
		\end{equation}
		Since $V \in \mathscr{C}^\infty_0(\overline{\Omega}^\diese \setminus T_\ord)$, there exists a real number $0 < b < 1/2$ such that $\restr{V}{\Omega^{\diese, b}} \in \mathscr{C}^\infty_0(\Omega^{\diese, b})$. Consequently, for any $i\in \llbracket 1, \ord-1\rrbracket $, the surface integral on $\Sigma^\diese_{i, a}$ is reduced to the set $\Sigma^{\diese, b}_{i, a}$ defined by \eqref{eq:Sigma_ia_securite}. When $k$ tends to $+\infty$, we can then use the trace continuity result stated in Proposition \ref{prop:trace_L2loc_demi_cylindre} on $\Sigma^{\diese, b}_{i, a}$, to deduce that \eqref{eq:Green_formula} is satisfied by $U$ and $V$.
	\end{dem}
	\color{black}
	%
	%

\subsubsection{An oblique change of variables}\label{sub:oblique_cov}
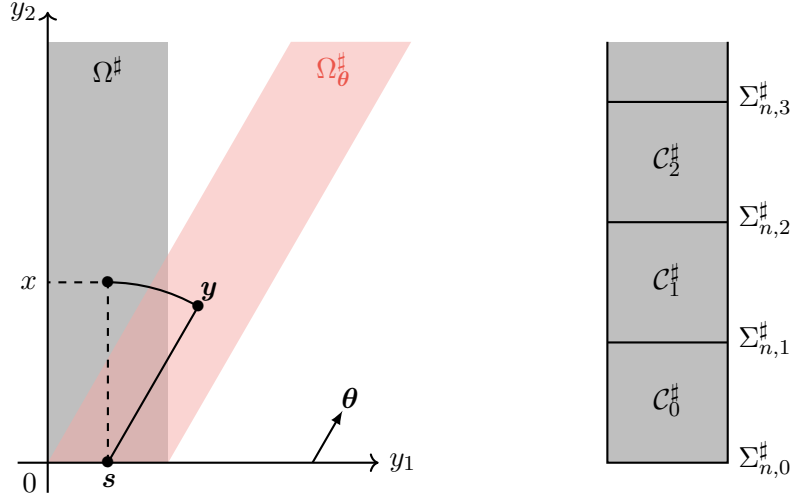
\begin{figure}
	\centering
	\setlength{\belowcaptionskip}{-10pt}
	\begin{tikzpicture}[thick, scale=0.8]
		\def\cotecellule{2}
		\def\nbcellules{3}
		%
		%
		\begin{scope}[xshift=0]
			\pgfmathsetmacro\lenX{2.5*\cotecellule+0.25*\cotecellule}
			\pgfmathsetmacro\bandeY{\nbcellules*\cotecellule+0.5*\cotecellule}
			\pgfmathsetmacro\lenY{\nbcellules*\cotecellule+0.75*\cotecellule}
			\def\anglecoupe{60}
			\pgfmathsetmacro\chgtvarS{0.5*\cotecellule}
			\def\chgtvarX{3.0}
			\fill[black!50, opacity=0.5] (\cotecellule , 0) -- (\cotecellule , \bandeY) -- (0 , \bandeY) -- (0 , 0) -- cycle;
			\fill[rougeTerre!50, opacity=0.5] (\cotecellule , 0) -- ({\cotecellule + \bandeY / tan(\anglecoupe)}, \bandeY) -- ({0 + \bandeY / tan(\anglecoupe)}, \bandeY) -- (0 , 0) -- cycle;
			\draw[->] (-0.25*\cotecellule, 0) -- (\lenX, 0) node[right] {$\yi_1$};
			\draw[->] (0, -0.25*\cotecellule) -- (0, \lenY) node[left]  {$\yi_2$};
			\draw (\chgtvarS, 0) node {$\bullet$} node[below] {$\sv$} -- ({\chgtvarS + \chgtvarX*cos(\anglecoupe)}, {\chgtvarX*sin(\anglecoupe)}) node {$\bullet$} node [pos=1.1] {$\yv$} arc(\anglecoupe:90:\chgtvarX) node {$\bullet$};
			\draw[dashed] (\chgtvarS, 0) -- (\chgtvarS, \chgtvarX) -- (0, \chgtvarX) node[left] {$\xv$};
			\draw[-latex] (0.8*\lenX, 0) -- +({1.0*cos(\anglecoupe)}, {1.0*sin(\anglecoupe)}) node [pos=1.25] {$\cut$};
			\draw (0.5*\cotecellule, {\bandeY-0.25*\cotecellule}) node {$\Omega^\diese$};
			\draw[rougeTerre] ({0.5*\cotecellule + (\bandeY-0.25*\cotecellule) / tan(\anglecoupe)}, {\bandeY-0.25*\cotecellule}) node {$\Omega^\diese_\cut$};
			\draw (0, 0) node[below left] {$0$};
		\end{scope}
		%
		%
		\begin{scope}[xshift=0.6\textwidth]
			\fill[fill=black!50, opacity=0.5] (0, \nbcellules*\cotecellule+0.5*\cotecellule) -- (0, 0) -- (\cotecellule, 0) -- (\cotecellule, \nbcellules*\cotecellule+0.5*\cotecellule);
			\draw (0, \nbcellules*\cotecellule+0.5*\cotecellule) -- (0, 0) -- (\cotecellule, 0) -- (\cotecellule, \nbcellules*\cotecellule+0.5*\cotecellule);
			\draw (\cotecellule, 0) node[right] {$\Sigma^\diese_{\ord, 0}$};
			\foreach [evaluate=\idI as \idJ using {int(\idI-1)}] \idI in {1,...,\nbcellules} {
        \draw (0, \idI*\cotecellule) -- (\cotecellule, \idI*\cotecellule);
        \draw (\cotecellule, \idI*\cotecellule) node[right] {$\Sigma^\diese_{n, \idI}$};
        \draw (0.5*\cotecellule, \idI*\cotecellule-0.5*\cotecellule) node {$\mathcal{C}^\diese_\idJ$};
      }
		\end{scope}
	\end{tikzpicture}
	\caption{The half-cylinders $\Omega^\diese$ and $\Omega^\diese_\cut$ (left), and the domains $\mathcal{C}^\diese_\ell$ and $\Sigma^\diese_{\ord, k}$ (right) for $\ord = 2$ \label{fig:domains}}
\end{figure}

Before stating \textcolor{surligneur}{Proposition \ref{prop:trace_lifting} which is} the main result of this section, let us introduce the change of variables in $\R^\ord_+$:
\begin{equation}\label{eq:changt_vars}
	(\sv, x) \in \R^n_+ \mapsto \yv = (\sv, 0) + \xv\,\cut \in \R^\ord_+,
\end{equation}
and denote by $\Omega^\diese_\cut$ the image of $\Omega^\diese$ by the above transformation:
\begin{equation}\label{eq:Omegadiesetheta}
	\Omega^\diese_\cut := \{ (\sv, 0) + \xv\, \cut, \ \ \sv \in (0, 1)^{\ord-1},\ \xv > 0\}.
\end{equation}
This is illustrated in Figure \ref{fig:domains_3D_bis} for $n = 3$ and in Figure \ref{fig:domains} for $\ord = 2$ and $|\cut| = 1$. The following simple lemma will be used in the sequel. %
%
%
\begin{lem}\label{lem:integrale}
	For any $V \in L^1(\Omega^\diese)$, we have
	\begin{equation}
		\displaystyle\label{eq:integrale}
		\int_{\Omega^\diese_\cut} \widetilde{V}(\yv)\; d\yv = \int_{\Omega^\diese} \widetilde{V}(\yv)\; d\yv,
	\end{equation}
	where $\widetilde{V} \in L^1_{\textit{loc}}(\R^\ord_+)$ denotes the periodic extension of $V$, defined by \eqref{eq:per_extension}.
\end{lem}
\begin{dem}
	We will use the notation $\kv = (\ki_1, \dots, \ki_d) \in \Z^d$ for a vector of integers. For any set $\mathcal{O} \subset \R^\ord$, let $\smash{\mathbbm{1}_{\mathcal{O}}}$ be the indicator function of $\mathcal{O}$, that is, the function which equals $1$ in $\mathcal{O}$ and $0$ elsewhere.
By density, it suffices to prove \eqref{eq:integrale} for $V \in \mathscr{C}^\infty_0(\Omega^\diese)$. By additivity of integration,
	\[
		\int_{\Omega^\diese_\cut} \widetilde{V}(\yv)\; d\yv = \int_{\R^\ord_+} \mathbbm{1}_{\Omega^\diese_\cut}(\yv) \; \widetilde{V}(\yv)\; d\yv = \sum_{\kv \in \Z^{\ord-1}} \int_{\Omega^\diese + (\kv, 0)} \mathbbm{1}_{\Omega^\diese_\cut}(\yv) \; \widetilde{V}(\yv)\; d\yv,
	\]
	where the sum over $\kv \in \Z^{\ord-1}$ is finite because of $\mathbbm{1}_{\Omega^\diese_\cut}$ and because $V$ is compactly supported. We then use the change of variables $\zv \mapsto \zv + (\kv, 0)$ which leads to
	\begin{align}
		\displaystyle \int_{\Omega^\diese_\cut} \widetilde{V}(\yv)\; d\yv &= \sum_{\kv \in \Z^{\ord-1}} \int_{\Omega^\diese} \mathbbm{1}_{\Omega^\diese_\cut}(\zv + (\kv, 0)) \; \widetilde{V}(\zv)\; d\zv \quad \textnormal{because $\widetilde{V}$ is periodic} \nonumber
		\\
		&= \int_{\Omega^\diese} \Big[ \sum_{\kv \in \Z^{\ord-1}} \mathbbm{1}_{\Omega^\diese_\cut - (\kv, 0)} (\zv) \Big]\; \widetilde{V}(\zv)\; d\zv \quad \textnormal{by linearity.}\label{eq:preuve_lemme_integrales_1}
	\end{align}
	Furthermore, by noticing that the collection of sets $\{\Omega^\diese_\cut - (\kv, 0),\ \ \kv \in \Z^{\ord-1}\}$ forms a partition of $\R^\ord_+$, it follows that
	\begin{equation}
		\displaystyle
		\spforall \zv \in \Omega^\diese, \quad \sum_{\kv \in \Z^{\ord-1}} \mathbbm{1}_{\Omega^\diese_\cut - (\kv, 0)} (\zv) = \mathbbm{1}_{\R^\ord_+}(\zv) = 1.
		\label{eq:preuve_lemme_integrales_2}
	\end{equation}
	Combining \eqref{eq:preuve_lemme_integrales_1} and \eqref{eq:preuve_lemme_integrales_2} implies that \eqref{eq:integrale} is satisfied for $V \in \mathscr{C}^\infty_0(\Omega^\diese)$.
\end{dem}

\vspace{1\baselineskip}\noindent
The inversion of the change of variables \eqref{eq:changt_vars} leads us to introduce:
\begin{equation}
	\displaystyle \label{eq:transverse_coo}
	\spforall \yv \in \R^\ord, \quad \sv_\cut(\yv) := \hat{\yv\;} - (\yi_\ord/\cuti_\ord)\, \hat{\cut\;} \in \R^{\ord-1},
\end{equation}
so that,
\begin{equation}
 \yv = (\sv, 0) + \xv\,\cut \quad \Longleftrightarrow \quad \sv = \sv_\cut(\yv) \quad \textnormal{and}\quad \xv = \yi_\ord/\cuti_\ord.
\end{equation}
%
%
%
%
%
%
\noindent
The next proposition emphasizes the fact that through the change of variables \eqref{eq:changt_vars}, the differential operator $\Dt{}$ simply becomes the partial derivative with respect to $\yi_\ord$ (which is obvious for smooth functions).
%
%
%
\begin{prop}\label{prop:trace_lifting}
	Let $\Psi \in L^2(\Omega^\diese)$. Then the periodic function $\Psi_\cut$ defined as
	\begin{equation}\label{eq:fonction}
		\displaystyle
			\aeforall \yv \in \R^\ord_+, \quad \widetilde{\Psi}_\cut(\yv) := \widetilde{\Psi}(\sv_\cut(\yv),\yi_\ord/\cuti_\ord),
	\end{equation}
 	(where $\widetilde{\Psi}$ is the periodic extension of $\Psi$) belongs to  $L^2(\Omega^\diese)$ and
	 \begin{equation}
		\label{eq:norme_fonction}
		\begin{array}{r@{\ =\ }l}
			\displaystyle\|{\Psi}_\cut\|_{L^2(\Omega^\diese)} &\displaystyle \sqrt{\cuti_\ord}\; \|{\Psi}\|_{L^2(\Omega^\diese)}.
		\end{array}
	\end{equation}
	Moreover, if $\partial_{\yi_\ord}\Psi\in L^2(\Omega^\diese)$, then $\Psi_\cut$ belongs to $H^1_{\cut, \textit{per}}(\Omega^\diese)$ with  directional derivative
	\begin{equation}
		\label{eq:derivee_fonction}
		\displaystyle
		\aeforall \yv \in \R^\ord_+, \quad D_\cut \widetilde{\Psi}_\cut (\yv)  = \frac{\partial \widetilde{\Psi}}{\partial \yi_\ord}(\sv_\cut(\yv),\yi_\ord/\cuti_\ord). 
	\end{equation}
\end{prop}
\begin{dem}
	The map $(\sv, x) \mapsto (\sv, 0) + x\, \cut$ from $\Sigma^\diese_{\ord, 0} \times \R_+$ to $\Omega^\diese_\cut$ defines a $\mathscr{C}^1$--diffeomorphism with a non-vanishing Jacobian $\cuti_\ord \neq 0$. Therefore, by using the definition \eqref{eq:Omegadiesetheta} of $\Omega^\diese_\cut$, a change of variables as well as the property $s_\cut((\sv, 0) + x\, \cut) = \sv$, we obtain that
	\[
		\setlength{\abovedisplayskip}{4pt}
		\setlength{\belowdisplayskip}{4pt}
		\displaystyle
		\int_{\Omega^\diese_\cut} |\widetilde{\Psi}_\cut(\yv)|^2\; d\yv = \cuti_\ord\; \int_{\Sigma^\diese_{\ord, 0}} \int_0^{+\infty} |\widetilde{\Psi}_\cut((\sv, 0) + x\, \cut)|^2\; d\xv\, d\sv = \cuti_\ord\; \int_{\Sigma^\diese_{\ord, 0}} \int_0^{+\infty} |\widetilde{\Psi}(\sv, x)|^2\; d\xv\, d\sv.
	\]
	We deduce from Lemma \ref{lem:integrale} that $\Psi_\cut \in L^2(\Omega^\diese)$, and that \eqref{eq:norme_fonction} holds.

	\vspace{1\baselineskip} \noindent
	Now in order to derive the expression of $D_\cut \widetilde{\Psi}_\cut$ in the sense of distributions, consider a test function $\Phi \in \mathscr{C}^\infty_0(\R^\ord_+)$. The change of variables $(\sv, \xv) \mapsto (\sv, 0) + x\, \cut$ combined with Fubini's theorem for integrable functions leads to
	\begin{equation}
		\int_{\R^\ord_+} \widetilde{\Psi}_\cut(\yv)\; D_\cut \Phi(\yv) \; d\yv = \cuti_\ord \int_{\R^{\ord-1}} \int_0^{+\infty} \widetilde{\Psi}(\sv, x)\; D_\cut \Phi((\sv, 0) + \xv\, \cut) \; d\xv d\sv.
		\label{eq:preuve_fonction_var_sep2}
	\end{equation}
	Furthermore the 1D function $\phi_{\sv, \cut}$ defined by $\phi_{\sv, \cut}(\xv) :=  \Phi((\sv, 0) + \xv\, \cut)$ belongs to $\mathscr{C}^\infty_0(\R_+)$ and we have $[d\phi_{\sv, \cut}/d\xv](\xv) = D_\cut \Phi((\sv, 0) + \xv\, \cut)$ from the chain rule. Since $\partial_{\yi_\ord}\Psi$ is in $L^2$, we can integrate by parts the inner integral in \eqref{eq:preuve_fonction_var_sep2} to obtain
	\begin{align}
		\int_{\R^\ord_+} \widetilde{\Psi}_\cut(\yv)\; D_\cut \Phi(\yv) \; d\yv  &= -\cuti_\ord \int_{\R^{\ord-1}} \int_0^{+\infty}  \frac{\partial \Psi}{\partial \yi_\ord}(\sv,x)\; \phi_{\sv, \cut}(\xv) \; d\xv d\sv \nonumber
		\\
		&= - \int_{\R^\ord_+}  \frac{\partial \Psi}{\partial \yi_\ord}(\sv_\cut(\yv),y_\ord/{\cuti_\ord})\; \Phi(\yv) \; d\yv,
	\end{align}
	where the last equality comes from the change of variables $\yv \mapsto (\sv_\cut(\yv), \yi_\ord/\cuti_\ord)$. This gives the expression of $D_\cut \widetilde{\Psi}_\cut$ in \eqref{eq:derivee_fonction}.
\end{dem}
\begin{rmk}
	\setlength{\abovedisplayskip}{1pt}
	\setlength{\belowdisplayskip}{1pt}
	It will be often useful to use \eqref{eq:derivee_fonction} in the form
	\begin{equation}
		\aeforall (\sv, \xv) \in \R^\ord_+, \quad D_\cut \widetilde{\Psi}_\cut ((\sv, 0) + x\, \cut)  = \frac{\partial \widetilde{\Psi}}{\partial \yi_\ord}(\sv, \xv).\label{eq:derivee_fonction_2}
	\end{equation}
\end{rmk}
\noindent
The previous proposition allows in particular to deduce the surjectivity of the trace operator from $H^1_{\cut, \textit{per}}(\Omega^\diese)$ to $L^2(\Sigma^\diese_{\ord, 0})$.
\begin{cor}\label{cor:trace_lifting}
	Let $\varphi \in L^2(\Sigma^\diese_{\ord, 0})$, and $\psi \in H^1(\R_+)$ such that $\psi(0)=1$. Then the periodic function defined by
	\begin{equation}\label{eq:fonction_var_sep}
		\displaystyle
			\aeforall \yv \in \R^\ord_+, \quad \mathcal{R}\varphi \,(\yv) := \widetilde{\varphi}(\sv_\cut(\yv))\; \psi(\yi_\ord/\cuti_\ord)
	\end{equation}
 belongs to $H^1_{\cut, \textit{per}}(\Omega^\diese)$, and its trace is $\restr{\mathcal{R}\varphi}{\Sigma^\diese_{\ord, 0}} = \varphi$. Moreover, $\mathcal{R}$ defines a continuous map from $L^2(\Sigma^\diese_{\ord, 0})$ to $H^1_{\cut, \textit{per}}(\Omega^\diese)$.
\end{cor}

\subsection{Link with a periodic half-guide problem}
\label{sec:link_with_a_periodic_half_guide_problem}
For any boundary data $\varphi \in L^2(\Sigma^\diese_{\ord, 0})$, we can now introduce $U^+_\cut$ as the solution in $H^1_{\cut}(\Omega^\diese)$ of the half-guide problem
\begin{equation}
	\label{eq:half_guide_problem}
	\left|
	\begin{array}{r@{\ }c@{\ }ll}
	\displaystyle - \Dt{} \big( \mu_p \; \Dt{} U^+_\cut \big) - \rho_p \; \omega^2 \; U^+_\cut &=& 0, \quad \textnormal{in}\quad \Omega^\diese,
	\\[8pt]
	\displaystyle \restr{U^+_\cut}{\Sigma^\diese_{\ord, 0}} &=& \varphi,
	\\[8pt]
	\displaystyle \restr{U^+_\cut}{\Sigma^\diese_{i, 0}} &=& \restr{U^+_\cut}{\Sigma^\diese_{i, 1}}  &\spforall i \in \llbracket 1, \ord-1 \rrbracket,
	\\[8pt]
	\restr{\mu_p\; \Dt{} U^+_\cut}{\Sigma^\diese_{i, 0}} &=& \restr{\mu_p\; \Dt{} U^+_\cut}{\Sigma^\diese_{i, 1}} &\spforall i \in \llbracket 1, \ord-1 \rrbracket.
	\end{array}
	\right.
\end{equation}
Note that the third equation above implies that $U^+_\cut \in H^1_{\cut, \textit{per}}(\Omega^\diese)$, the first one implies that $\smash{\mu_p\; \Dt{} U^+_\cut \in H^1_{\cut}(\Omega^\diese)}$, and finally the fourth one implies that $\smash{\mu_p\; \Dt{} U^+_\cut \in H^1_{\cut, \textit{per}}(\Omega^\diese)}$. The space of the boundary data \textcolor{surligneur}{can} seem surprising {compared to the Helmholtz equation with an elliptic principal part}, but recall from Corollary \ref{cor:trace_lifting} that the trace mapping on $\Sigma^\diese_{\ord, 0}$ is surjective from $H^1_{\cut, \textit{per}}(\Omega^\diese)$ to $L^2(\Sigma^\diese_{\ord, 0})$.

\vspace{1\baselineskip} \noindent
With the functional framework introduced in the previous section, we can now show that Problem \eqref{eq:half_guide_problem} is well-posed.
\begin{prop}\label{prop:half_guide_FV}
	For any $\varphi \in L^2(\Sigma^\diese_{\ord, 0})$, Problem \eqref{eq:half_guide_problem} is equivalent to the variational formulation
	\begin{equation}
		\label{eq:half_guide_FV}
		\left|
		\begin{array}{l}
		\textnormal{\textit{Find $U^+_\cut \in H^1_{\cut, \textit{per}}(\Omega^\diese)$ such that $\restr{U^+_\cut}{\Sigma^\diese_{\ord, 0}} = \varphi$ and}}
		\\[8pt]
		\displaystyle \spforall V \in H^1_{\cut, \textit{per}}(\Omega^\diese) \textnormal{ such that }\,\restr{V}{\Sigma^\diese_{\ord, 0}}=0, \quad \int_{\Omega^\diese} \left( \mu_p \; \Dt{} U^+_\cut \; \Dt{} \overline{V} - \rho_p \; \omega^2 \; U^+_\cut \; \overline{V} \right) = 0,
		\end{array}
		\right.
	\end{equation}
	for which Lax-Milgram's theorem applies.
\end{prop}

\begin{dem}
	The variational formulation is obtained {by multiplying the first equation of \eqref{eq:half_guide_problem} by $V \in H^1_{\cut, \textit{per}}(\Omega^\diese)$, and by using} Green's formula \eqref{eq:Green_formula_H1per}. The application of the Lax-Milgram's theorem in $\{V \in H^1_{\cut, \textit{per}}(\Omega^\diese),\ \gamma_{\ord, 0} V = 0\}$, thanks to Corollary \ref{cor:trace_lifting}, is direct.

	\vspace{1\baselineskip}\noindent
	For the equivalence, as usual, one picks test functions $V\in \mathscr{C}^\infty_0(\Omega^\diese)$ to deduce that the solution $U^+_\cut \in H^1_{\cut, \textit{per}}(\Omega^\diese)$ of \eqref{eq:half_guide_FV} satisfies the first equation of \eqref{eq:half_guide_problem}. This implies that  $\mu_p\; \Dt{} U^+_\cut \in H^1_{\cut}(\Omega^\diese)$. The real difficulty is to show that $U^+_\cut$ satisfies the fourth equation in \eqref{eq:half_guide_problem} or equivalently that $\mu_p\; \Dt{} U^+_\cut \in H^1_{\cut, \textit{per}}(\Omega^\diese)$. According to Proposition \ref{prop:trace_L2loc_demi_cylindre}, we have
	\[\displaystyle
		\spforall 1 \leq i \leq \ord-1, \quad \restr{\mu_p\; \Dt{} U^+_\cut}{\Sigma^\diese_{i, a}} \in L^2_\textit{loc}(\Sigma^\diese_{i, a}).
	\]
	Therefore, Proposition \ref{prop:formule_Green_demi_cylindre} allows us to use Green's formula \eqref{eq:Green_formula} for $U = \mu_p\; \Dt{} U^+_\cut $ and for $V \in \mathscr{C}^\infty_0(\overline{\Omega}^\diese \setminus T_\ord) \cap H^1_{\cut, \textit{per}}({\Omega}^\diese)$, where $T_\ord$ is the skeleton defined in \eqref{eq:definition_table}. By combining this with the fact that $U^+_\cut$ solves \eqref{eq:half_guide_FV} and the first equation of \eqref{eq:half_guide_problem}, one obtains that for any integer $i \in \llbracket 1, \ord-1 \rrbracket$,
	\[
		 \spforall V\in \mathscr{C}^\infty_0(\overline{\Omega}^\diese \setminus T_\ord)\cap H^1_{\cut, \textit{per}}({\Omega}^\diese),\quad \Big( \int_{\Sigma^\diese_{i,1}} \mu_p\; \Dt{} U^+_\cut \; \overline{V} \;d s - \int_{\Sigma^\diese_{i,0}} \mu_p\; \Dt{} U^+_\cut \; \overline{V} \;d s \Big)=0.
	\]
	Furthermore, $\mathscr{C}^\infty_0(\Sigma^\diese_{i, 0})$ is included in $\{\restr{V}{\Sigma^\diese_{i, 0}},\ V \in \mathscr{C}^\infty_0(\overline{\Omega}^\diese \setminus T_\ord)\cap H^1_{\cut, \textit{per}}({\Omega}^\diese) \}$. In fact, any $\psi \in \mathscr{C}^\infty_0(\Sigma^\diese_{i, 0})$ admits the extension $\Psi: \yv \in \Omega^\diese \mapsto \psi(\yi_1,\dots,\yi_{i-1},\yi_{i+1},\dots,\yi_\ord)$, which belongs to $\mathscr{C}^\infty_0(\overline{\Omega}^\diese \setminus T_\ord)\cap H^1_{\cut, \textit{per}}({\Omega}^\diese)$. Finally, since $\mathscr{C}^\infty_0(\Sigma_{i,0}^\diese)$ is dense in $L^2(\Sigma^\diese_{i,0})$, it is easy to show that the fourth equation of \eqref{eq:half_guide_problem} holds and that $\restr{\mu_p\; \Dt{} U^+_\cut}{\Sigma^\diese_{i, 1}} \in L^2(\Sigma^\diese_{i, 1})$ for any $i \in \llbracket 1, \ord-1 \rrbracket$.
\end{dem}
%
%
%
%

\vspace{0\baselineskip}\noindent
We now make the link between $U^+_\cut(\varphi)$ and the solution of the half-line problem \eqref{eq:half_line_problem} that fully justifies the introduction of the half-guide problem \eqref{eq:half_guide_problem}. 

\vspace{1\baselineskip}\noindent
To do so, first, let us introduce the quasiperiodic coefficients defined for any $\sv \in \R^{\ord-1}$ by
	\begin{equation}
	\displaystyle
	\spforall \xv \in \R, \quad \mu_{\sv, \cut}(\xv) := \mu_p\big((\sv, 0) + \xv\, \cut\big) \quad \textnormal{and} \quad \rho_{\sv, \cut}(\xv) := \rho_p\big((\sv, 0) + \xv\, \cut\big),
	\end{equation}
	as well as the one-dimensional problems
	\begin{equation}
	\left|
	\begin{array}{r@{\ }c@{\ }l@{\quad}l}
	\displaystyle - \frac{d}{d \xv} \Big( \mu_{\sv, \cut} \; \frac{d u^+_{\sv, \cut}}{d \xv} \Big) - \rho_{\sv, \cut} \; \omega^2 \; u^+_{\sv, \cut} &=& 0, \quad \textnormal{in} & \R_+,
	\\[8pt]
	\displaystyle u^+_{\sv, \cut}(0) &=& 1.
	\end{array}
	\right.
	\label{eq:half_line_problems}
	\end{equation}
	Note that \eqref{eq:half_line_problem} corresponds to \eqref{eq:half_line_problems} taken with $\sv = 0$.

	\vspace{1\baselineskip} \noindent
	Under the assumptions \eqref{eq:coef_ellipt} and \eqref{eq:dissipation}, Problem \eqref{eq:half_line_problems} admits a unique solution $u^+_{\sv, \cut}$ in $H^1(\R_+)$ for any $\sv \in \R^{\ord-1}$. Moreover, $u^+_{\sv, \cut}$ decays exponentially at infinity, uniformly with respect to $\sv$, that is, there exist constants $\alpha, c > 0$ depending only on $\mu_\pm, \rho_\pm$ such that
	\begin{equation}
		\label{eq:exp_decay_halfline_s}
		\quad \spforall \sv \in \R^{\ord-1}, \quad \big\|\euler^{-\alpha \Imag \omega\, \xv}\, u^+_{\sv, \cut}\big\|_{H^1(\R^+)} \leq c.
	\end{equation}
	Furthermore, thanks to the continuity of $\mu_p$ and $\rho_p$, we can show that $u^+_{\sv, \cut}$ is continuous with respect to $\sv$, as stated in the next proposition.
	%
	%
	\begin{prop}
		\label{prop:properties_us}
		The mapping $\sv \in \R^{\ord-1} \mapsto u^+_{\sv, \cut}$, which associates with a real vector $\sv$ the solution in $H^1(\R_+)$ of the problem \eqref{eq:half_line_problems}, defines a uniformly continuous function which is periodic of period $1$ in each direction. 
	\end{prop}
	\begin{dem}
		%

		To show that $\sv \mapsto u^+_{\sv, \cut}$ is $1$--periodic in each direction, one simply has to note that since $\mu_{\sv,\cut}$ and $\rho_{\sv, \cut}$ are $1$--periodic with respect to each $\si_i$, both $u^+_{\sv, \cut}$ and $u^+_{\sv + \vec{\ev}_i, \cut}$ satisfy the same half-line problem \eqref{eq:half_line_problems}. Thus, by well-posedness of \eqref{eq:half_line_problems}, $u^+_{\sv, \cut} = u^+_{\sv + \vec{\ev}_i, \cut}$.

		\vspace{1\baselineskip} \noindent
		Now let us prove the regularity of $\sv \mapsto u^+_{\sv, \cut}$. For any $\sv_1, \sv_2 \in \R^{\ord-1}$, by writing the variational formulations satisfied by $u^+_{\sv_1, \cut}$ and $u^+_{\sv_2, \cut}$, and by substracting one from the other, we obtain
		\begin{multline*}\displaystyle
			\spforall v  \in H^1_0(\R_+), \quad \int_{\R_+} \Big[ \mu_{\sv_1, \cut}\; \frac{d}{d \xv} (u^+_{\sv_1, \cut} - u^+_{\sv_2, \cut}) \; \overline{\frac{d v}{d \xv}} - \rho_{\sv_1, \cut} \; \omega^2 \; (u^+_{\sv_1, \cut} - u^+_{\sv_2, \cut}) \; \overline{v} \Big] =
			\\
			\int_{\R_+} \Big[ (\mu_{\sv_2, \cut} - \mu_{\sv_1, \cut})\; \frac{d u^+_{\sv_2, \cut}}{d \xv} \; \overline{\frac{d v}{d \xv}} - (\rho_{\sv_1, \cut} - \rho_{\sv_2, \cut}) \; \omega^2 \; u^+_{\sv_2, \cut} \Big].
		\end{multline*}
		Now choose $v = u^+_{\sv_1, \cut} - u^+_{\sv_2, \cut} \in H^1_0(\R_+)$ in the above equality. The well-posedness of \eqref{eq:half_line_problems}, a Cauchy-Schwarz inequality applied to the right-hand side and \eqref{eq:exp_decay_halfline_s} imply that there exists a real number $c > 0$ independent of $\sv$ and $\cut$ such that
		\begin{equation}\label{eq:cont}\displaystyle
		\big\|u^+_{\sv_1, \cut} - u^+_{\sv_2, \cut} \big\|_{H^1(\R_+)} \leq c\;  \Big( \|\mu_{\sv_2, \cut} - \mu_{\sv_1, \cut}\|_\infty + \|\rho_{\sv_2, \cut} - \rho_{\sv_1, \cut}\|_\infty \Big).
		\end{equation}
		The functions $ \mu_p$ and $\rho_p$ are continuous and $1$--periodic in each direction: from Heine-Cantor theorem, they are uniformly continuous. Let us define the modulus of uniform continuity
		\[
			\spforall \mu \in \mathscr{C}^0(\R^\ord),\ \spforall \eps>0, \quad \delta(\mu,\eps)=\sup_{\itbf{y},\itbf{z}}\{|\mu(\itbf{y})-\mu(\itbf{z})|,\;|\itbf{y}-\itbf{z}|<\eps\}
		\]
		A function $\mu$ is uniformly continuous if $\delta(\mu,\eps)$ tends to $0$ as $\eps$ tends to $0$.
		It follows from \eqref{eq:cont} that
		\[
			\big\|u^+_{\sv_1, \cut} - u^+_{\sv_2, \cut} \big\|_{H^1(\R_+)} \leq c\;  \Big(\delta(\mu_p,|\sv_1-\sv_2|) + \delta(\rho_p,|s_1-s_2|) \Big).
		\]
		Therefore, $\sv \mapsto u^+_{\sv, \cut}$ is continuous from $\R^{n-1}$ in $H^1(\R^+)$.
	\end{dem}

\begin{prop}\label{prop:structure_half_guide}
	Let $\sv_\cut$ be the mapping defined by \eqref{eq:transverse_coo}, and $\widetilde{U}^+_\cut$ (resp. $\widetilde{\varphi}$) be the periodic extension of $U^+_\cut$ (resp. $\varphi$) the solution of \eqref{eq:half_guide_problem}. Then, we have
	\begin{equation}
		\label{eq:concatenation_half_guide}
		\aeforall \yv \in \R^\ord_+, \quad \widetilde{U}^+_\cut(\widetilde{\varphi})(\yv) = \widetilde{\varphi}\big(\sv_\cut(\yv)\big)\; u^+_{\sv_\cut(\yv), \cut} (\yi_\ord/\cuti_\ord),
	\end{equation}
	{\color{surligneur}or equivalently
	\begin{equation}\displaystyle\label{eq:preuve_lien_2D_1D_1}
		\aeforall (\sv, \xv) \in \R^{\ord-1} \times \R_+, \quad \widetilde{U}^+_\cut(\widetilde{\varphi})((\sv, 0) + \cut\, \xv) = \widetilde{\varphi}(\sv)\; u^+_{\sv, \cut} (\xv).
	\end{equation}}
	Moreover if $\widetilde{\varphi}$ is continuous in the neighbourhood of $0$ and satisfies $\widetilde{\varphi}(0) = 1$, then
	\begin{equation}\label{eq:lien_2D_1D}
		\aeforall \xv \in \R, \quad u^+_\cut(\xv) = \widetilde{U}^+_\cut(\widetilde{\varphi})(\xv\, \cut)
	\end{equation}
\end{prop}
\begin{dem}
	We begin by proving \eqref{eq:concatenation_half_guide}. Let us denote for $\aeforall \yv \in \R^\ord_+$, $U^{}_1(\yv)$ the right-hand side of \eqref{eq:concatenation_half_guide}. Note that $\Psi:(\sv,x)\mapsto \widetilde{\varphi}(\sv)\; u^+_{\sv, \cut} (x)$ is $1$--periodic with respect to $\sv$ (thanks to Proposition \ref{prop:properties_us}), and belongs to $L^2(\Omega^\diese)$ since
	\[
	\|\Psi\|^2_{L^2(\Omega^\diese)} = \int_{\Sigma^\diese_{\ord, 0}} |\varphi(\sv)|^2\ \|u^+_{\sv, \cut}\|^2_{L^2(\R_+)} \; d \sv \leq \cuti_\ord\, c^2 \; \|\varphi\|^2_{L^2(\Sigma^\diese_{\ord, 0})},\ \ \textnormal{with}\ \ c = \sup_{\sv} \|u^+_{\sv, \cut}\|_{L^2(\R_+)}.
	\]
	Moreover, since for all $\sv$, $u^+_{\sv, \cut} \in H^1(\R^+)$,  $\partial_{y_n}\Psi$ is also in $L^2(\Omega^\diese)$ (using similar inequalities to the above). By Proposition \ref{prop:trace_lifting}, $U^{}_1$ belongs to $H^1_{\cut,\textit{per}}(\Omega^\diese)$  with
	\[
		\displaystyle
		\aeforall \yv \in \R^\ord_+, \quad D_\cut\, \widetilde{U_1} (\yv)  = \widetilde{\varphi}\big(\sv_\cut(\yv)\big)\; \frac{d u^+_{\sv_\cut(\yv), \cut}}{d x} (\yi_\ord/\cuti_\ord).
	\]
	Finally, since $u^+_{\sv, \cut}(0) = 1$, it is clear that $\restr{U^{}_1}{\Sigma^\diese_{\ord, 0}} = \varphi$. By repeating the same argument, we can show that $\mu_p D_\cut\, U^{}_1$ belongs to $H^1_{\cut,\textit{per}}(\Omega^\diese)$ with
	\[
		\displaystyle
		\aeforall \yv \in \R^\ord_+, \quad D_\cut\, [\mu_p\,D_\cut\, \widetilde{U_1}] (\yv)  = \widetilde{\varphi}\big(\sv_\cut(\yv)\big)\; \frac{d}{d x}\Big(\mu_{\sv_\cut(\yv),\cut} \frac{d u^+_{\sv_\cut(\yv), \cut}}{d x}\Big) (\yi_\ord/\cuti_\ord).
	\]
	Since  $u^+_{\sv, \cut}$ satisfies \eqref{eq:half_line_problems}, it is clear that $U_1$ satisfies \eqref{eq:half_guide_problem}. By well-posedness of \eqref{eq:half_guide_problem}, we have $U^{}_1 = U^+_\cut$.

	\vspace{1\baselineskip}\noindent
	\textcolor{surligneur}{The equivalence between \eqref{eq:concatenation_half_guide} and \eqref{eq:preuve_lien_2D_1D_1} is directly obtained using the change of variables $(\sv, \xv) \mapsto ((\sv, 0) + \cut\, \xv)$. Moreover, %
	w}e have from Proposition \ref{prop:properties_us} that $\sv \mapsto u^+_{\sv, \cut}$ is continuous. If in addition to that, $\widetilde{\varphi}$ is continuous in a neighbourhood of $0$, then \eqref{eq:preuve_lien_2D_1D_1} becomes true \emph{for any} $\sv$ in that neighbourhood. In particular, \eqref{eq:preuve_lien_2D_1D_1} can be written for $\sv = 0$, thus leading to \eqref{eq:lien_2D_1D}.
\end{dem}
\noindent
In particular, we deduce from the above proprosition that
\begin{equation}
	\label{eq:link_DU_du}
	\aeforall \yv \in \R^\ord_+, \quad \Dt{} \widetilde{U}^+_\cut(\widetilde{\varphi})(\yv) = \widetilde{\varphi}\big(\sv_\cut(\yv)\big)\; \frac{d u^+_{\sv_\cut(\yv), \cut}}{d\xv} (\yi_\ord/\cuti_\ord).
\end{equation}
%
%
\begin{rmk}
	The half-guide solution $U^+_\cut$ depends on $\varphi$ whereas $u^+_{\sv, \cut}$ does not. In this sense, the relation \eqref{eq:concatenation_half_guide} \textcolor{surligneur}{can} seem surprising at first sight. Numerical results presented in Section \ref{sec:numerical_results} will illustrate this property.
	%
	%
\end{rmk}

\section{Resolution of the half-guide problem}
\label{sec:resolution_half_guide_problem}
The advantage of the lifting process lies in the periodic nature of \eqref{eq:half_guide_problem}, which allows us to exploit tools that are well-suited for periodic waveguides. %
In this paper, we use a DtN-based method \cite{flissthese, jolyLiFliss}, developed for the elliptic\footnote{By \textit{elliptic} Helmholtz equation, we refer to the Helmholtz equation with an \textit{elliptic principal part}.} Helmholtz equation $- \nabla \cdot (\mu_p\; \nabla U) - \rho_p\; \omega^2\; U = 0$ in unbounded periodic guides. This method does not rely on decay properties, and therefore remains robust when the absorption tends to $0$. As we essentially transpose this method to our directional Helmholtz equation, we will see below that the framework remains exactly the same, although the analysis has to be adapted. Let us mention the \textit{recursive doubling method} \cite{yuan2007recursive, ehrhardt2008numerical}, suited for bounded periodic waveguides, and a method \cite{zhang2021numerical} based on the Floquet-Bloch transform, although its extension to our non-elliptic equation seems unclear.


\vspace{1\baselineskip} \noindent
In what follows, $\mathcal{C}^\diese_\ell$ is the cell defined for every $\ell \in \N$ by
\begin{equation}
\mathcal{C}^\diese_0 = (0, 1)^\ord \quad \textnormal{and} \quad \mathcal{C}^\diese_\ell = \mathcal{C}^\diese_0 + \ell\, \vec{\ev}_\ord, \quad \textnormal{so that} \quad \Omega^\diese = \bigcup_{\ell \in \N} \mathcal{C}^\diese_\ell.
\end{equation}
For $\ell > 0$, we call $\Sigma^\diese_{\ord, \ell}$ the interface between the cells $\mathcal{C}^\diese_\ell$ and $\mathcal{C}^\diese_{\ell+1}$, that is, $\Sigma^\diese_{\ord, \ell} = \Sigma^\diese_{\ord, 0} + \ell\, \vec{\ev}_\ord$. By periodicity, each cell $\mathcal{C}^\diese_\ell$ can be identified to $\mathcal{C}^\diese_0$. Similarly, each interface $\Sigma^\diese_{\ord, \ell}$ can be identified to $\Sigma^\diese_{\ord, 0}$. The cells and interfaces are represented in Figure \ref{fig:domains}.

\subsection{Structure of the solution}
\label{sec:structure_of_the_solution}
The solution $U^+_\cut(\varphi)$ of \eqref{eq:half_guide_problem} has a particular structure that we explain in this section. Denote by $\mathcal{P} \in \mathcal{L}\big(L^2(\Sigma^\diese_{\ord, 0})\big)$ the operator
\begin{equation}
	\label{eq:definition_propagation_operator}
	\spforall \varphi \in L^2(\Sigma^\diese_{\ord, 0}), \quad \mathcal{P}\varphi := \restr{U^+_\cut(\varphi)}{\Sigma^\diese_{\ord, 1}},
\end{equation}
where $L^2(\Sigma^\diese_{\ord, 1})$ and $L^2(\Sigma^\diese_{\ord, 0})$ have been identified to each other in an obvious manner. This identification will be used systematically in what follows, even if not mentioned. Note that the operator $\mathcal{P}$ is well-defined, due to the continuity of the trace operator on $\Sigma^\diese_{i, a}$ \eqref{eq:trace_func_per}.
\begin{prop}
	\label{prop:periodic_structure_solution}
  For any $\varphi$ in $L^2(\Sigma^\diese_{\ord, 0})$, we have
	\begin{equation}
		\label{eq:periodic_structure_solution}
		\spforall \ell \in \N,\;\aeforall \yv \in \Omega^\diese, \quad U^+_\cut(\varphi)(\yv + \ell\, \vec{\ev}_\ord) = U^+_\cut(\mathcal{P}^\ell \varphi)(\yv).
	\end{equation}
	Moreover, the spectral radius of $\mathcal{P}$ is strictly less than one.
\end{prop}
\begin{dem}
	We only present the outline of the proof, which is quite similar to the one in \cite{flissthese, jolyLiFliss}. Given $\varphi \in L^2(\Sigma^\diese_{\ord, 0})$, consider the function $U^{}_1$ defined in $\Omega^\diese$ by $U^{}_1(\yv) = U^+_\cut(\varphi)(\yv + \vec{\ev}_\ord)$ for almost any $\yv \in \Omega^\diese$. Since the coefficients $\mu_p$ and $\rho_p$ are periodic, one deduces that $U^{}_1$ satisfies the volume equation as well as the periodicity condition in \eqref{eq:half_guide_problem}. Furthermore,
	\[\displaystyle
		\restr{U^{}_1}{\Sigma^\diese_{\ord, 0}} = \restr{U^+_\cut(\varphi)}{\Sigma^\diese_{\ord, 1}} = \mathcal{P}\varphi.
	\]
	Thus, by well-posedness of \eqref{eq:half_guide_problem}, we have \eqref{eq:periodic_structure_solution} for $\ell = 1$. The result \eqref{eq:periodic_structure_solution} for $\ell \geq 2$ is proved by induction.

	\vspace{1\baselineskip} \noindent
	It remains to show that the spectral radius is strictly less than $1$. To this end, by analogy with \eqref{eq:exp_decay_halfline_s}, one can show the existence of constants $\alpha, c > 0$ such that
	\begin{equation}\label{eq:exp_decay_halfguide}\displaystyle
	\spforall \varphi \in L^2(\Sigma^\diese_{\ord, 0}), \quad \big\|e^{\alpha\Imag \omega\, \yi_\ord/\cuti_\ord}\, U^+_\cut \big\|_{H^1_\cut(\Omega^\diese)} \leq c\; \|\varphi\|_{L^2(\Sigma^\diese_{\ord, 0})}.
	\end{equation}
	Since $\mathcal{P}^\ell \varphi = U^+_\cut(\varphi)(\cdot, \ell)$, the estimate above implies that $\|\mathcal{P}^\ell\| \leq c\; e^{-\alpha\Imag \omega\, \ell/\cuti_\ord}$.	Hence, using Gelfand's formula \cite[\S 10.3]{rudin1991functional}, the spectral radius can be estimated as follows:
	\[
	\rho(\mathcal{P}) = \lim_{\ell \to +\infty}\|\mathcal{P}^\ell\|^{1/\ell} \leq 
	e^{-\beta\Imag \omega/\cuti_\ord} < 1.
	\]%
\end{dem}
%

\noindent
The operator $\mathcal{P}$ is called the \emph{propagation operator}, as it describes how the solution of \eqref{eq:half_guide_problem} evolves from one interface to another. Provided that $\mathcal{P}$ is known, the solution $U^+_\cut(\varphi)$ may be constructed using \emph{local cell problems}. Let us first introduce the appropriate functional framework in a periodicity cell
\begin{equation}
	\label{eq:H1thetaper_celle}
	\displaystyle
	H^1_{\cut, \textit{per}}(\mathcal{C}_0^\diese) := \Big\{ U \in H^1_{\cut}(\mathcal{C}_0^\diese),\  \widetilde{U} \in H^1_{\cut, \textit{loc}}(\mathcal{B}_0) \Big\} ,
\end{equation}
where $\mathcal{B}_0:=\R^\ord_+ \cap\{ 0<y_\ord<1 \}$. Similarly to Section \ref{sub:trace}, one can show that any function of $H^1_{\cut, \textit{per}}(\mathcal{C}_0^\diese)$ has a $L^2$ trace on the boundary of $\mathcal{C}_0^\diese$. We can prove in particular that
	\[
	H^1_{\cut, \textit{per}}(\mathcal{C}^\diese_0) = \Big\{U \in H^1_\cut(\mathcal{C}^\diese_0)\ /\ \restr{U}{\yi_i = 0} = \restr{U}{\yi_i = 1},\ \spforall i \in \llbracket 1, \ord-1 \rrbracket \Big\}.
	\]
	We can now introduce the local cell problems: %
for all $\varphi\in L^2(\Sigma^\diese_{\ord, 0})$, for $j\in\{0,1\}$, let $E^j(\varphi)\in  H^1_{\cut, \textit{per}}(\mathcal{C}^\diese_0) $ satisfy
\begin{equation}
	\label{eq:local_cell_problem}
	\left|
	\begin{array}{r@{\ }c@{\ }l@{\quad}l}
	\displaystyle - \Dt{} \big( \mu_p \; \Dt{} E^j \big) - \rho_p \; \omega^2 \; E^j &=& 0, \quad \textnormal{in} & \mathcal{C}^\diese_0,
	\\[4pt]
	\multicolumn{4}{l}{ 
	\restr{\mu_p\; \Dt{} E^j}{\yi_i = 0} =\restr{\mu_p\; \Dt{} E^j}{\yi_i = 1} \quad \spforall i \in \llbracket 1, \ord-1 \rrbracket,}
	\end{array}
	\right.
\end{equation}
defined for $j = 0, 1$, with the boundary conditions
\begin{equation}
	\label{eq:local_cell_BC}
	\left|
	\begin{array}{c@{\quad}c@{\quad}c}
		\restr{E^0}{\Sigma^\diese_{\ord, 0}} = \varphi &\textnormal{and}& \restr{E^0}{\Sigma^\diese_{\ord, 1}} = 0,
		\\[4pt]
		\restr{E^1}{\Sigma^\diese_{\ord, 0}} = 0 &\textnormal{and}& \restr{E^1}{\Sigma^\diese_{\ord, 1}} = \varphi.
	\end{array}
	\right.
\end{equation}
A variational formulation can be derived as in Proposition \ref{prop:half_guide_FV}, and
the well-posedness follows once again from a lifting argument (see Proposition \ref{prop:trace_lifting}) combined with Lax-Milgram's theorem in $H^1_{\cut, \textit{per}}(\mathcal{C}^\diese_0)$.

\vspace{1\baselineskip} \noindent
Proposition \ref{prop:periodic_structure_solution} implies that $\restr{U^+_\cut(\varphi)(\cdot + \ell\, \vec{\ev}_\ord)}{\Sigma^\diese_{\ord, 0}} = \mathcal{P}^\ell \varphi$. Hence, if the propagation operator $\mathcal{P}$ is known, by linearity, the solution of the half-guide problem can be entirely constructed cell by cell as follows:
\begin{equation}
	\displaystyle
	\spforall \ell \in \N, \quad \restr{U^+_\cut(\varphi)(\cdot + \ell\, \vec{\ev}_\ord)}{\mathcal{C}^\diese_0} = E^0(\mathcal{P}^{\ell} \varphi) + E^1(\mathcal{P}^{\ell+1} \varphi).
	\label{eq:UfromEis}
\end{equation}
\subsection{Characterization of the propagation operator: the Riccati equation}\label{sec:Riccati}
%
In the sequel, $\langle \cdot, \cdot \rangle$ denotes the canonical $L^2$ scalar product on $\Sigma^\diese_{\ord, 0}$ (or equivalently on $\Sigma^\diese_{\ord, 1}$).

\vspace{1\baselineskip} \noindent
In order to characterize the propagation operator $\mathcal{P}$, it is useful to introduce the \emph{local DtN operators} $\mathcal{T}^{jk} \in \mathcal{L}(L^2(\Sigma^\diese_{\ord, 0}))$, defined for $j, k = 0, 1$ by
\begin{equation}
	\label{eq:DtNloc}
	\spforall \varphi \in L^2(\Sigma^\diese_{\ord, 0}), \quad \mathcal{T}^{jk} \varphi = (-1)^{k+1}\; \cuti_\ord\; \restr{\left[\mu_p\; \Dt{} E^j(\varphi)\right]}{\Sigma^\diese_{\ord, k}}.
\end{equation}
where $E^j(\varphi)$ satisfies \eqref{eq:local_cell_problem}-\eqref{eq:local_cell_BC}.
By Green's formula \eqref{eq:Green_formula}, note that for all $j, k = 0, 1$ and for $(\varphi, \psi) \in L^2(\Sigma^\diese_{\ord, 0})^2$, these operators satisfy
\begin{equation}\label{eq:DtNloc_weak_form}
	\displaystyle
	\Big\langle \mathcal{T}^{j k} \varphi,\;{\psi} \Big\rangle = \int_{\mathcal{C}^\diese_0} \left[ \mu_p\; \Dt{} E^j(\varphi) \; \Dt{} \overline{E^k(\psi)} - \rho_p\; \omega^2\; E^j(\varphi)\; \overline{E^k(\psi)}\, \right].
\end{equation}
Before deriving other useful properties of the local DtN operators, we need to introduce some additional notations. For any closed operator $\mathcal{A} \in \mathcal{L}(L^2(\Sigma^\diese_{\ord, 0}))$, we denote $\mathcal{A}^*$ the adjoint of $\mathcal{A}$, and $\overline{\mathcal{A}}$ its «\,\textit{complex conjugate}\,», that is,
 	\[
 	\spforall \varphi \in L^2(\Sigma^\diese_{\ord, 0}), \quad \overline{\mathcal{A}}\varphi = \overline{\mathcal{A} \overline{\varphi}}.
	\]
%
\noindent
It is not difficult to see that $\overline{\mathcal{A}^*} = \overline{\mathcal{A}}^*$, and $\overline{\overline{\mathcal{A}}} = \mathcal{A}$.

\begin{prop}
	\label{prop:pties_local_DtN}
	The local DtN operators $\mathcal{T}^{jk}$ satisfy
	\begin{equation}
		\label{eq:adjoints_local_DtN}
		\left[\mathcal{T}^{00}\right]^* = \overline{\mathcal{T}^{00}}, \quad \left[\mathcal{T}^{11}\right]^* = \overline{\mathcal{T}^{11}}, \quad \left[\mathcal{T}^{01}\right]^* = \overline{\mathcal{T}^{10}}, \quad \left[\mathcal{T}^{10}\right]^* = \overline{\mathcal{T}^{01}}.
	\end{equation}
	Furthermore, the operators $\mathcal{T}^{00}$, $\mathcal{T}^{11}$, and $\mathcal{T}^{00} + \mathcal{T}^{11}$ are invertible.
\end{prop}
\begin{dem}
	The property \eqref{eq:adjoints_local_DtN} follows from Green's formula applied to $E^j(\varphi)$ and $\overline{E^k(\overline{\psi})}$, see for instance \cite[Proposition 2.2.4]{flissthese} in the case of the Helmholtz equation.

	\vspace{1\baselineskip} \noindent
	The operators $\mathcal{T}^{00}$, $\mathcal{T}^{11}$, and $\mathcal{T}^{00} + \mathcal{T}^{11}$ are bounded. We are going to show that they are also coercive. Their invertibility will then follow from Lax-Milgram's theorem.

	From the expression \eqref{eq:DtNloc_weak_form}, one has the existence of a constant $c \equiv c(\mu_-, \rho_-, |\omega|) > 0$ such that
	\[\displaystyle
		- |\omega|\; \Imag\left[ \frac{1}{\omega} \Big\langle\mathcal{T}^{kk}  \varphi,\; {\varphi} \Big\rangle \right] \geq c\; \Imag \omega\; \|E^k(\varphi)\|^2_{H^1_\cut(\mathcal{C}^\diese_0)} \geq \tilde{c}\, \Imag \omega\, \|\varphi\|^2_{L^2(\Sigma^\diese_{\ord, 0})},
	\]
	since from \eqref{eq:trace_func_per}, the trace application from $H^1_{\cut, \textit{per}}(\mathcal{C}^\diese_0)$ to $L^2(\Sigma^\diese_{\ord, 0})$ is continuous. %
	It follows that the operators $\mathcal{T}^{00}$ and $\mathcal{T}^{11}$ are coercive, and therefore invertible. The inequalities above summed for $k = 0, 1$ imply the coercivity and hence the invertibility of $\mathcal{T}^{00} + \mathcal{T}^{11}$ as well.
\end{dem}

\vspace{1\baselineskip} \noindent
As seen earlier, the solution of the half-guide problem \eqref{eq:half_guide_problem} is given by \eqref{eq:UfromEis}. Now let us use the characterization of $H^1_{\textit{per}, \cut}(\Omega^\diese)$, namely, Corollary \ref{cor:rest_H1thetaper} with $a = 1$, so that $\Omega^\diese_{a, -} = \mathcal{C}^\diese_0$ and $\Omega^\diese_{a, +} = \Omega^\diese \setminus \mathcal{C}^\diese_0$. Since $\mu_p\; \Dt{}U^+_\cut(\varphi)$ belongs to $H^1_{\cut, \textit{per}}(\Omega^\diese)$, the directional derivative of $U^+_\cut(\varphi)$ is continuous across the interface $\Sigma^\diese_{\ord, 1}$, \emph{i.e.}
\begin{equation}
	\restr{\left[\mu_p\; {\Dt{} U^+_\cut(\varphi)}\right]}{\Sigma^\diese_{\ord, 1}} = \restr{\left[\mu_p\; {\Dt{} U^+_\cut(\varphi)((\cdot + \vec{\ev}_\ord)}\right]}{\Sigma^\diese_{\ord, 0}},
\end{equation}
or equivalently,
\begin{equation}
	\begin{array}{l}
		\restr{\left[\mu_p\; \Dt{} E^0(\varphi)\right]}{\Sigma^\diese_{\ord, 1}} + \; \restr{\left[\mu_p\; \Dt{} E^1(\mathcal{P} \varphi)\right]}{\Sigma^\diese_{\ord, 1}}
		\\[15pt]
		\ \;\qquad\qquad\qquad\qquad=\; \restr{\left[\mu_p\; \Dt{} E^0(\mathcal{P} \varphi)\right]}{\Sigma^\diese_{\ord, 0}} +\; \restr{\left[\mu_p\; \Dt{} E^1(\mathcal{P}^2 \varphi)\right]}{\Sigma^\diese_{\ord, 0}}.
	\end{array}
\label{eq:cty_across_Sigma}
\end{equation}
By using the definition of the local DtN operators $\mathcal{T}^{jk}$, \eqref{eq:cty_across_Sigma} leads to the following characterization.
\begin{prop}
	The propagation operator $\mathcal{P}$ \textcolor{surligneur}{defined by \eqref{eq:definition_propagation_operator}} is the unique solution of the constrained Riccati equation
	\begin{equation}
		\label{eq:Riccati}
		\left|
		\begin{array}{l}
			\textnormal{\textit{Find $\mathcal{P} \in \mathcal{L}(L^2(\Sigma^\diese_{\ord, 0}))$ such that $\rho(\mathcal{P}) < 1$ and}}
			\\[12pt]
			\multicolumn{1}{c}{\displaystyle \mathcal{T}^{10}\mathcal{P}^2 + (\mathcal{T}^{00} + \mathcal{T}^{11})\, \mathcal{P} + \mathcal{T}^{01} = 0.}
		\end{array}
		\right.
	\end{equation}
\end{prop}
\begin{dem}
	The proof is identical to the one for the elliptic Helmholtz equation \cite[Theorem 4.1]{jolyLiFliss}. We know from Proposition \ref{prop:periodic_structure_solution} that $\mathcal{P}$ has a spectral radius which is strictly less than $1$. Moreover \eqref{eq:cty_across_Sigma} ensures that $\mathcal{P}$ satisfies the Riccati equation.

	\vspace{1\baselineskip} \noindent
	In order to prove the uniqueness, let us consider an operator $\mathcal{P}_1$ which satisfies \eqref{eq:Riccati}. The function defined cell by cell by
	\[\displaystyle
		\spforall \varphi \in L^2(\Sigma^\diese_{\ord, 0}), \quad \spforall \ell \in \N^*, \quad \restr{U_1(\varphi)(\cdot + \ell\, \vec{\ev}_\ord)}{\mathcal{C}^\diese_0} = E^0(\mathcal{P}_1^{\ell} \varphi) + E^1(\mathcal{P}_1^{\ell+1} \varphi),
	\]
	solves \eqref{eq:half_guide_problem} in each cell $\mathcal{C}_\ell$ and is continuous across each interface $\Sigma^\diese_{\ord, \ell}$, by definition \eqref{eq:local_cell_problem}, \eqref{eq:local_cell_BC} of $E^0$ and $E^1$. By Corollary \ref{cor:rest_H1thetaper}, $U_1$ is locally $H^1_\cut$ in $\Omega^\diese$.

	\vspace{1\baselineskip} \noindent
	Moreover, since $\mathcal{P}_1$ satisfies \eqref{eq:Riccati}, the directional derivative $\mu_p \Dt{} U_1$ is continuous across each interface. Thus, using Corollary \ref{cor:rest_H1thetaper}, we deduce that $U_1$ satisfies \eqref{eq:half_guide_problem} in $\Omega^\diese$.

	\vspace{1\baselineskip}\noindent
	Furthermore, given that $\rho(\mathcal{P}_1) < 1$, Gelfand's formula and the well-posedness of the cell problems ensure that there exist positive constants $c, \rho_*$, with $\rho_* < 1$ such that, for $\ell \in \N$ large enough,
	\[
	 \|U_1(\varphi)\|_{H^1_\cut(\mathcal{C}^\diese_\ell)} \leq c\; \rho_*^\ell\; \|\varphi\|_{L^2(\Sigma^\diese_{\ord, 0})}.
	\]
	Hence $U_1(\varphi)$ belongs to $H^1_{\cut, \textit{per}}(\Omega^\diese)$ and satisfies the half-guide problem \eqref{eq:half_guide_problem}. By well-posedness of \eqref{eq:half_guide_problem}, $U_1(\varphi)$ and $U^+_\cut(\varphi)$ coincide, and thus have the same trace on $\Sigma^\diese_{\ord, 1}$, that is $\mathcal{P}_1 \varphi = \mathcal{P} \varphi$ for any $\varphi \in L^2(\Sigma^\diese_{\ord, 0})$.
\end{dem}

\vspace{1\baselineskip} \noindent
As a consequence, the propagation operator can be obtained by solving the Riccati equation in \eqref{eq:Riccati}, and by choosing the unique solution whose spectral radius is strictly less than $1$. One important thing to retain from the above is that both the propagation operator and the solution of the half-guide problem only require the computation of $E^0$, $E^1$, and the operators $\mathcal{T}^{00}$, $\mathcal{T}^{10}$, $\mathcal{T}^{01}$, and $\mathcal{T}^{11}$, which involve problems defined on a periodicity cell. However, the resolution of the constrained Riccati equation \eqref{eq:Riccati} is not obvious at all. The properties of this equation are investigated in further details in Section \ref{sec:about_Riccati_equation}.

\subsection{The DtN operator and the DtN coefficient}
\label{sec:the_DtN_operator_and_the_DtN_coefficient}
The goal of this part is to see how the half-guide problem and the local cell problems can be used to compute the DtN coefficient $\lambda^+$. We recall that
\[\displaystyle
	\lambda^+ = - \mu_\cut(0) \; \frac{d u^+_\cut}{d \xv}(0).
\]
Therefore, it is natural to introduce the DtN operator $\Lambda \in \mathcal{L}(L^2(\Sigma^\diese_{\ord, 0}))$ defined by
\begin{equation}
	\label{eq:def_DtN_operator}
	\displaystyle
	\spforall \varphi \in L^2(\Sigma^\diese_{\ord, 0}), \quad \Lambda \varphi := - \cuti_\ord\; \restr{\left[\mu_p\; \Dt{} U^+_\cut(\varphi)\right]}{\Sigma^\diese_{\ord, 0}}.
\end{equation}
This operator also has the following properties, whose proof is exactly identical to the one of Proposition \ref{prop:pties_local_DtN}.
\begin{prop}
	\label{prop:pties_DtN}
	One has $\Lambda^* = \overline{\Lambda}$. Moreover, $\Lambda$ and $\Lambda + \mathcal{T}^{11}$ are invertible operators.
\end{prop}
\vspace{1\baselineskip} \noindent
\noindent
 Taking the directional derivative of \eqref{eq:UfromEis} (for $\ell=0$) on $\Sigma^\diese_{\ord, 0}$
and using the definition \eqref{eq:DtNloc} of the local DtN operators $\mathcal{T}^{00}$ and $\mathcal{T}^{10}$ leads to
\begin{equation}
	\label{eq:DtN_operator}
	\displaystyle
	\Lambda = \mathcal{T}^{00} + \mathcal{T}^{10} \mathcal{P}.
\end{equation}
Besides, by writing the formula \eqref{eq:link_DU_du} after multiplication by $\mu_p$, and by evaluating it for $\yv = (\sv, 0)$, so that $\sv_\cut(\yv) =\sv$, we obtain
%
%
\begin{equation}
	\label{eq:expression_Lambda}
	\Lambda \varphi(\sv) = \cuti_\ord\; \lambda^{}_\cut(\sv)\; \varphi(\sv), \quad \textnormal{with} \quad \lambda^{}_\cut(\sv) = - \Big[\mu_{\sv, \cut} \; \frac{d u^+_{\sv, \cut}}{d \xv}\Big](0),
\end{equation}
namely, $\Lambda$ is a multiplication operator. We deduce from \eqref{eq:expression_Lambda} the DtN coefficient $\lambda^+$.
\begin{prop}\label{prop:DtNcoeff_from_DtNoperator}
	The function $\lambda^{}_\cut : \R^{\ord-1} \to \C$ defined by \eqref{eq:expression_Lambda} is continuous. Moreover, if $\varphi \in \mathscr{C}_{\textit{per}}(\R^{\ord-1})$ is a given function which satisfies $\varphi(0) = 1$, then we have
	\begin{equation}
		\lambda^+ = \lambda^{}_\cut(0) = \frac{1}{\cuti_\ord}\; (\Lambda \varphi)(0).\label{eq:DtNcoeff_from_DtNoperator}
	\end{equation}
\end{prop}
\begin{dem}
	Using Green's formula, we have that for all $\sv \in \R^{\ord-1}$
	\[\displaystyle
		\lambda^{}_\cut(\sv) = a_\sv(u^+_{\sv, \cut}, u^+_{\sv, \cut}), \quad \textnormal{with} \quad a_\sv(u, v) = \int_{\R_+} \Big( \mu_{\sv, \cut} \; \frac{d u}{d \xv} \; \overline{\frac{d v}{d \xv}} - \rho_{\sv, \cut} \; \omega^2 \; u \; \overline{v} \Big).
	\]
	The continuity of $u \mapsto a_\sv(u, u)$ results directly from the properties of the coefficients $\mu_p$ and $\rho_p$. Moreover, Proposition \ref{prop:properties_us} ensures that the function $\sv \mapsto u^+_{\sv, \cut}$ is continuous. Therefore, as the composition of these two functions, $\lambda^{}_\cut$ is also continuous.

	If in addition $\varphi$ is continuous, then $\Lambda \varphi$ is also continuous. Hence, $(\Lambda \varphi) (0) = \cuti_\ord\; \lambda^{}_\cut(0) \varphi(0)$ which yields the desired result.
\end{dem}

\subsection{Spectral properties of the Riccati equation}
\label{sec:about_Riccati_equation}
We now present some properties regarding Equation \eqref{eq:Riccati}. \textcolor{surligneur}{These properties will be exploited for the numerical resolution of the Riccati equation, by constructing the operator $\mathcal{P}$ from its eigenpairs (this will be done in Section \ref{sec:discrete_Riccati_equation} after space discretization). For this reason, it is worhwhile to reformulate a spectral version (Proposition \ref{prop:spectrum_P_riccati}) of the Riccati equation that would characterize these eigenpairs, while taking into account the spectral radius constraint. This is precisely the purpose of this section.}

\vspace{1\baselineskip} \noindent
Recall that $\mathcal{T}(\mathcal{P}) = 0$, where $\mathcal{T}$ is the bounded operator defined by
\begin{equation}
	\spforall X \in \mathcal{L}\big(L^2(\Sigma^\diese_{\ord, 0})\big), \quad \mathcal{T}(X) = \mathcal{T}^{10} X^2 + (\mathcal{T}^{00} + \mathcal{T}^{11}) X + \mathcal{T}^{01}.
\end{equation}
In the sequel, we will write $\mathcal{T}(\lambda)$ for $\mathcal{T}(\lambda I)$. We begin with the following factorization lemma.
\begin{lem}
	\label{lme:factorization_Riccati_operator}
	Let $\mathcal{P}$ be the propagation operator defined by \eqref{eq:definition_propagation_operator}. For any number $\lambda \in \C$, 
	\begin{equation}
		\mathcal{T}(\lambda) = (\lambda \overline{\mathcal{P}^*} - I)\; (\Lambda + \mathcal{T}^{11})\; (\mathcal{P} - \lambda),
	\end{equation}
	where $\mathcal{T}^{11}$ is defined by \eqref{eq:DtNloc} and $\Lambda$ is defined by \eqref{eq:def_DtN_operator}.
\end{lem}
\begin{dem}
	Let $\lambda \in \C$. Since the propagation operator satisfies $\mathcal{T}(\mathcal{P}) = 0$, one obtains that
	\begin{align}
		\mathcal{T}(\lambda) &= \mathcal{T}(\lambda) - \mathcal{T}(\mathcal{P}) \nonumber\\
		&= \left[ \mathcal{T}^{10} (\lambda + \mathcal{P}) + \mathcal{T}^{00} + \mathcal{T}^{11} \right]\; (\lambda - \mathcal{P}) \nonumber\\
		&= (\lambda \mathcal{T}^{10} + \Lambda + \mathcal{T}^{11}) \; (\lambda - \mathcal{P}), \quad \textnormal{from \eqref{eq:DtN_operator}.}\label{eq:T(lambda I)}
	\end{align}
	We use once again the fact that $\mathcal{T}(\mathcal{P}) = 0$ which, by the expression \eqref{eq:DtN_operator}, is equivalent to $\mathcal{T}^{01} = -(\Lambda + \mathcal{T}^{11})\; \mathcal{P}$. By transposing this equation, and by taking the complex conjugate, one obtains that $\overline{\left[\mathcal{T}^{01}\right]^*} = -\overline{\mathcal{P}^* \vphantom{\cramped{\mathcal{T}^{11}}}}\; \overline{(\Lambda + \mathcal{T}^{11})^*}$. Since $\left[\mathcal{T}^{11}\right]^* = \overline{\mathcal{T}^{11}}$ and $\left[\mathcal{T}^{01}\right]^* = \overline{\mathcal{T}^{10}}$ as ensured by Proposition \ref{prop:pties_local_DtN}, and since $\Lambda^* = \overline{\Lambda}$ from Proposition \ref{prop:pties_DtN}, it follows that
	\[
		\mathcal{T}^{10} = -\overline{\mathcal{P}^*}\; (\Lambda + \mathcal{T}^{11}).
	\]
	Inserting this expression of $\mathcal{T}^{10}$ in \eqref{eq:T(lambda I)} therefore leads to
	\[
		\mathcal{T}(\lambda) = \left[ - \lambda \overline{\mathcal{P}^*}\; (\Lambda + \mathcal{T}^{11}) + \Lambda + \mathcal{T}^{11} \right] \; (\lambda - \mathcal{P}) = (I - \lambda \overline{\mathcal{P}^*})\;  (\Lambda + \mathcal{T}^{11})\; (\lambda - \mathcal{P}).
	\]
	which is the desired result.
\end{dem}
\vspace{1\baselineskip} \noindent
The previous factorization lemma allows one to characterize the spectrum of the propagation operator as follows.
\begin{prop}\label{prop:spectrum_P_riccati}
	For any complex number $\lambda$, one has
	\begin{equation}
		\label{eq:spectrum_P_riccati}
		\displaystyle
		\lambda \in \sigma(\mathcal{P})\quad \Longleftrightarrow \quad 0 \in \sigma\big[\mathcal{T}(\lambda)\big]\ \ \text{and}\ \ |\lambda| < 1.
	\end{equation}
\end{prop}
\begin{dem}
	Proving \eqref{eq:spectrum_P_riccati} amounts to showing that for any $\lambda \in \C$ such that $|\lambda| < 1$, $\mathcal{P} - \lambda$ is invertible if and only if $\mathcal{T}(\lambda)$ is invertible. To this end, using Lemma \ref{lme:factorization_Riccati_operator}, it is sufficient to prove that $(\lambda \overline{\mathcal{P}^*} - I)\; (\Lambda + \mathcal{T}^{11})$ is an invertible operator. Proposition \ref{prop:pties_DtN} ensures the invertibility of $\Lambda + \mathcal{T}^{11}$ already. It thus remains to show that $\lambda \overline{\mathcal{P}^*} - I$ is invertible, which is true when $|\lambda| < 1$.

	Indeed, if $\lambda = 0$, then $\lambda \overline{\mathcal{P}^*} - I = - I$ is obviously invertible. Otherwise, it is not difficult to see that $\mathcal{P}$ and $\overline{\mathcal{P}^*}$ have the same spectrum. Hence, given that $|1/\lambda| > 1 > \rho(\overline{\mathcal{P}^*})$, it follows that $1/\lambda$ does not belong to $\sigma(\overline{\mathcal{P}^*})$. In other words, $\overline{\mathcal{P}^*} - (1/\lambda)\,I$ is an invertible operator.
\end{dem}

\begin{rmk}
	Note that the property \eqref{eq:spectrum_P_riccati} can be proved easily (and without Lemma \ref{lme:factorization_Riccati_operator}) for the point spectrum:
	\begin{equation}
		\label{eq:point_spectrum_P_riccati}
		\displaystyle
		\lambda \in \sigma_p(\mathcal{P})\quad \Longleftrightarrow \quad 0 \in \sigma_p\big[\mathcal{T}(\lambda)\big]\ \ \text{and}\ \ |\lambda| < 1.
	\end{equation}
 	This property was already proved in \cite{jolyLiFliss} for the Helmholtz equation. In this context, this was sufficient since the operator $\mathcal{P}$ was compact, which is no longer the case here.
\end{rmk}
\vspace{1\baselineskip} \noindent
Finally, it is worth noting that the values $\lambda \neq 0$ for which $0 \in \sigma\big[\mathcal{T}(\lambda)\big]$ can be paired in the following way.
\begin{prop}\label{prop:pairs_Riccati}
	For any complex number $\lambda \neq 0$, one has the following equivalence:
	\begin{equation}
		\displaystyle
		0 \in \sigma\big[\mathcal{T}(\lambda)\big] \quad \Longleftrightarrow \quad 0 \in \sigma\big[\mathcal{T}(1/\lambda)\big].
	\end{equation}
\end{prop}
\begin{dem}
	Let $\lambda \in \C^*$. From the properties of the local DtN operators (see Proposition \ref{prop:pties_local_DtN}), we deduce that
	\begin{equation}
		\overline{\left[\mathcal{T}(\lambda)\right]^*} = \lambda^2\, \mathcal{T}^{01} + \lambda (\mathcal{T}^{00} + \mathcal{T}^{11}) + \mathcal{T}^{10} = \lambda^2\, \mathcal{T}(1/\lambda).
	\end{equation}
	The operators $\mathcal{T}(\lambda)$ and $\overline{\left[\mathcal{T}(\lambda)\right]^*}$ have the same spectrum, hence the result.
\end{dem}
%
\begin{rmk}
	\label{rmk:pairs_Riccati}
	As Proposition \ref{prop:pairs_Riccati} shows, the values $\lambda \neq 0$ for which
	\[
		0 \in \sigma\big[\mathcal{T}(\lambda)\big]
	\]
	come by pairs $(\lambda,\lambda^{-1})$. From a numerical point of view, it suffices to choose $\lambda$ such that $|\lambda|<1$ and discard $\lambda^{-1}$.
\end{rmk}
\subsection{Spectral properties of the propagation operator}
\label{sec:propagation_operator}
\textcolor{surligneur}{This section, contrary to Section \ref{sec:about_Riccati_equation} is not related to the construction of our numerical method; it is of theoretical interest. On one hand, the result of this section, that is Proposition \ref{prop:properties_P}, is useful for interpreting some of the numerical results in Section \ref{sec:results:absorption}. On the other hand, it emphasizes the differences between the spectral properties of $\mathcal{P}$, and the ones of the corresponding operator for classical waveguide problems.} For the elliptic Helmholtz equation, $\mathcal{P}$ is compact (see \cite[Theorem 3.1]{jolyLiFliss}) and its spectrum hence consists only in \textcolor{surligneur}{isolated eigenvalues which accumulate to $0$. However, the picture is completely different in this case, because the spectrum has no isolated points}.

\vspace{1\baselineskip} \noindent
One useful way to study the properties of the propagation operator (especially its spectrum) is through an analytic formula: according to \eqref{eq:concatenation_half_guide}, $\mathcal{P}$ can be expressed for all $\varphi$ in $L^2(\Sigma^\diese_{\ord, 0})$ and for $\sv \in \R^{\ord-1}$ as
\begin{equation}
	\label{eq:expression_P}
	\mathcal{P} \varphi(\sv) = p_\cut(\sv)\; \widetilde{\varphi}\big(\sv - \cutslope\big), \quad \textnormal{with} \quad p_\cut(\sv) = u^+_{\sv - \cutslope, \cut} (1/\cuti_\ord) \quad \textnormal{and} \quad \cutslope = \hat{\cut\,}/\cuti_\ord \in \R^{\ord-1}.
\end{equation}
Note that since $\cut$ is an irrational vector, $\cutslope$ is also an irrational vector.

The properties of the mapping $\sv \mapsto u^+_{\sv, \cut}$ stated in Proposition \ref{prop:properties_us} imply that the fonction $p_\cut$ is continuous and $1$-periodic in each direction.

\vspace{1\baselineskip} \noindent
Operators that can be written under the form \eqref{eq:expression_P} are known as \emph{weighted shift operators}, and have been studied for instance in \cite{antonevich2012linear}. In particular, the spectral properties of $\mathcal{P}$ are given by the following result.
{%
\begin{prop}
	\label{prop:properties_P}
	Let $p_\cut: \Sigma^\diese_{\ord, 0} \to \C$ be the function defined in \eqref{eq:expression_P}. Then, $p_\cut(\sv) \neq 0$ for all $s$ in $\Sigma^\diese_{\ord, 0}$, and the spectral radius of $\mathcal{P}$ is given by
	\begin{equation}
		\label{eq:radius_spectrum}
		\rho(\mathcal{P})= \exp \left(\int_{\Sigma^\diese_{\ord, 0}} \log |p_\cut(\sv)|\; d \sv \right).
	\end{equation}
	Moreover, the spectrum of $\mathcal{P}$ is a circle of radius $\rho(\mathcal{P})$.
\end{prop}
\noindent
	 This result can be found in \cite[Theorem 2.1]{antonevich2012linear} for $\ord = 2$. We give below the proof for $\ord > 2$, which requires the following lemma (see Theorem 6.1 and Example 6.1 of \cite{kuipers}), known as a particular case of Birkhoff's ergodic theorem for continuous functions.
	\begin{lem}
		\label{lem:mean_value_discrete}
		Let $\psi : \Sigma^\diese_{\ord, 0} \to \C$ be continuous and $1$--periodic in each direction. Let $\alpha \in \R^{\ord-1}$ be an irrational vector. Then, we have the following uniform convergence:
		\[\displaystyle
			\lim_{\ell \to +\infty} \Big\|\frac{1}{\ell} \sum_{m = 0}^{\ell-1} \psi(\cdot - m \boldsymbol{\alpha}) - \int_{\Sigma^\diese_{\ord, 0}} \psi \Big\|_{\infty} = 0.
		\]
	\end{lem}
}
\begin{dem}[of Proposition \ref{prop:properties_P}]
	Let us first show by contradiction that $p_\cut$ or equivalently the function $\sv \mapsto u^+_{\sv, \cut}(1/\cuti_\ord)$ is nowhere vanishing. To do so, we use an argument of unique continuation. In fact, assume that there exists $s \in \Sigma^\diese_{\ord, 0}$ such that $u^+_{\sv, \cut}(1/\cuti_n) = 0$. Then $u^+_{\sv, \cut}$ satisfies the problem
	\[
	- \frac{d}{d \xv} \left( \mu_{\sv, \cut} \; \frac{d u^+_{\sv, \cut}}{d \xv} \right) - \rho_{\sv, \cut} \; \omega^2 \; u^+_{\sv, \cut} = 0, \ \textnormal{in} \ (1/\cuti_n,+\infty), \quad \textnormal{and} \quad u^+_{\sv, \cut}(1/\cuti_n) = 0.
	\]
	From the well-posedness of this problem, it follows that $u^+_{\sv, \cut} = 0$ in $(1/\cuti_n,+\infty)$. Therefore, by unique continuation, one deduces that $u^+_{\sv, \cut} = 0$ in $\R_+$, which contradicts the boundary condition $u^+_{\sv, \cut}(0) = 1$.

	\vspace{1\baselineskip}\noindent{
	We now establish the expression of the spectral radius $\rho(\mathcal{P})$. One has $\smash{\displaystyle \rho(\mathcal{P}) = \lim_{\ell \to +\infty}\|\mathcal{P}^\ell\|^{1/\ell}}$ from Gelfand's formula, and by induction, $\mathcal{P}^\ell$ can be expressed under the form
	\[\displaystyle
	\mathcal{P}^\ell \varphi(\sv) = p^{(\ell)}_\cut(\sv)\; \varphi(\sv - \ell\cutslope), \quad \textnormal{with} \quad p^{(\ell)}_\cut(\sv) = \prod_{m = 0}^{\ell-1} p_\cut(\sv - m \cutslope).
	\]
	Since the translation operator $\varphi \mapsto \varphi(\cdot - \ell\cutslope)$ is isometric and bijective, the norm of $\mathcal{P}^\ell$ is equal to the norm of the multiplication operator $\varphi \mapsto p^{(\ell)}_\cut\, \varphi$, that is $\|p^{(\ell)}_\cut\|_\infty$. Hence, given that $p_\cut(\sv) \neq 0$ for all $\sv$, one has
	\[
	\rho(\mathcal{P}) = \lim_{\ell \to +\infty} \Big\|\prod_{m = 0}^{\ell-1} p_\cut(\cdot - m \cutslope) \Big\|^{1/\ell}_\infty = \lim_{\ell \to +\infty} \exp  \Big\|\frac{1}{\ell} \sum_{m = 0}^{\ell-1} \log\big(|p_\cut(\cdot - m \cutslope)|\big)  \Big\|_\infty
	\]
	Since $\cut$ is an irrational vector, $\cutslope = \hat{\cut} / \cuti_\ord$ is also an irrational vector. Therefore, Lemma \ref{lem:mean_value_discrete} can be applied with $\boldsymbol{\alpha} = \cutslope$, and $\psi : \sv \mapsto \log |p_\cut(\sv )|$, which is well-defined and continuous. Hence the spectral radius is given by
	\[
		\rho(\mathcal{P}) = M_{\log}(p_\cut) := \exp \left(\int_{\Sigma^\diese_{\ord, 0}} \log |p_\cut(\sv)|\; d \sv \right).
	\]
	%
	 Let us now characterize the spectrum. To begin, note that the inverse of $\mathcal{P}$ is well-defined, since $p_\cut$ vanishes nowhere: for all $\varphi\in L^2(\Sigma^\diese_{\ord, 0}),\; \mathcal{P}^{-1} \varphi(\sv) := p_\cut(\sv)^{-1}\; \widetilde{\varphi}\big(\sv + \cutslope\big)$. Therefore, all the computations above can be applied to $\mathcal{P}^{-1}$, thus yielding
		\[
		\rho(\mathcal{P}^{-1}) = M_{\log}(p^{-1}_\cut) = \frac{1}{M_{\log}(p_\cut)} = \frac{1}{\rho(\mathcal{P})}
		\]
		Since the spectrum of $\mathcal{P}$ is always included in the annulus $\smash{\displaystyle \rho(\mathcal{P}^{-1})^{-1} \leq |z| \leq \rho(\mathcal{P})}$, it follows that $\sigma(\mathcal{P})$ is included in the circle $|z| = \rho(\mathcal{P}) = M_{\log}(p_\cut)$.

		\vspace{1\baselineskip} \noindent
		Conversely, for $\kv \in \Z^{\ord-1}$, let $S_\kv$ be the multiplication operator by $\sv \in \R^{\ord-1} \mapsto \exp(2\icplx\pi\, \kv\cdot \sv)$. From the expression \eqref{eq:expression_P} of the propagation operator, we obtain that
		\[
			S_\kv\; \mathcal{P}\; S^{-1}_\kv = e^{2\icplx\pi\, \kv\, \cdot\, \cutslope}\, \mathcal{P}.
		\]
		The operators $\mathcal{P}$ and $e^{2\icplx\pi \kv\, \cdot\, \cutslope}\, \mathcal{P}$ are similar, and thus have the same spectrum. Now consider an element $\lambda_0$ of $\sigma(\mathcal{P})$. Then, $|\lambda_0| =  M_{\log}(p_\cut)$, and $\lambda_\kv := e^{2\icplx\pi \kv\, \cdot\, \cutslope}\, \lambda_0$ also belongs to $\sigma(\mathcal{P})$ for all $\kv \in \Z^{\ord-1}$. Since $\cutslope$ is irrational, we have from Kronecker's theorem (Theorem \ref{thm:kronecker}) that the set $\{\lambda_\kv,\ \kv \in \Z^{\ord-1}\}$ is dense in the circle $|z| = |\lambda_0| = M_{\log}(p_\cut)$. Consequently, this whole circle is included in the spectrum, since the latter is a closed set.
	}
\end{dem}


%

%
%

\section{Resolution algorithm and discretization issues for \texorpdfstring{$\ord = 2$}{n = 2}} \label{sec:resolution_algorithm}
In order to compute the solution of Equation \eqref{eq:whole_line_problem}, the previous sections provide an algorithm which sums up as follows:
%

\begin{enumerate}[label=\arabic*., ref = \arabic*]
	\item\label{item:step1} \textit{\underline{Compute the solution $u^+_\cut$ of \eqref{eq:half_line_problems_0} and the DtN coefficient $\lambda^+$ defined by \eqref{eq:DtN_coefficients}}} by using the following procedure:
	\begin{enumerate}[label=(\textit{\alph*})., ref = \textit{1.\alph*}]
		\item\label{item:stepa} for any boundary data $\varphi \in L^2(\Sigma^\diese_{\ord, 0})$, compute the solutions $E^0(\varphi)$, $E^1(\varphi)$ of the local cell problems \eqref{eq:local_cell_problem};
		\item\label{item:stepb} compute the local DtN operators $(\mathcal{T}^{00}, \mathcal{T}^{01}, \mathcal{T}^{10}, \mathcal{T}^{11})$, defined by \eqref{eq:DtNloc}--\eqref{eq:DtNloc_weak_form};
		\item\label{item:stepc} compute the propagation operator $\mathcal{P}$ as the unique solution of the constrained Riccati equation \eqref{eq:Riccati};
		\item\label{item:stepd} using an arbitrarily chosen boundary data $\varphi \in \mathscr{C}_{\textit{per}}(\R^{\ord-1})$ which satisfies $\varphi(0) = 1$,
		\begin{itemize}
			\item from \eqref{eq:UfromEis}, construct the solution $U^+_\cut$ of the half-guide problem cell by cell;
			\item deduce the half-line solution $u^+_\cut$ via the formula \eqref{eq:lien_2D_1D};
		 \end{itemize}
		 \item\label{item:stepe} compute the DtN operator $\Lambda$ defined by \eqref{eq:DtN_operator}, and deduce $\lambda^+$ from \eqref{eq:DtNcoeff_from_DtNoperator}.
	\end{enumerate}
	\item\label{item:step2} \textit{\underline{Compute the solution $u^-_\cut$ of \eqref{eq:half_line_problems_0} and the DtN coefficient $\lambda^-$ defined by \eqref{eq:DtN_coefficients}}} by \textcolor{surligneur}{using exactly the same procedure as in} Step \ref{item:step1} (but independently from Step \ref{item:step1}).
	\item Finally, \textit{\underline{solve the interior problem \eqref{eq:interior_problem} in $(-a, a)$, and extend the solution everywhere}} by using \eqref{eq:solution_of_whole_line_problem}, as well as Step \ref{item:step1} and Step \ref{item:step2}.
\end{enumerate}
%

\noindent
For convenience sake, the quasiperiodicity order is set to $\ord = 2$. The most original aspects of the algorithm are the steps \eqref{item:stepa}--\eqref{item:stepd}, and the rest of this section focuses on the discretization of these four steps. We present in Sections \ref{sec:methode_2D} and \ref{sec:methode_quasi1D} two different methods that are linked to a choice of discretization of the step \eqref{item:stepa}, which influences the implementation of the steps \eqref{item:stepb} and \eqref{item:stepd}. The treatment of the step \eqref{item:stepc} is independent of this choice, and will be presented in Section \ref{sec:discrete_Riccati_equation}.

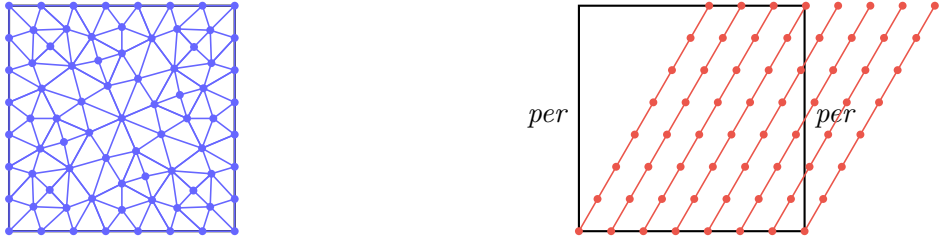
\begin{figure}[H]
	\centering
	\begin{subfigure}[b]{0.48\textwidth}
		\centering
		\begin{tikzpicture}[scale=3]
				\draw[thick] (0, 0) rectangle (1, 1);
				\ifthenelse{\boolean{afficherGraphes}}{
					\plotmsh{figures/cell.msh}{bleuDoux!75}{bleuDoux!75, semithick}
				}{}
		\end{tikzpicture}
	\end{subfigure}
	\hfill
	\begin{subfigure}[b]{0.48\textwidth}
		\centering
		\begin{tikzpicture}[scale=3]
			\def\cutangle{60}
			\xdef\nbnoeuds{8}
			\draw[thick] (0, 0) rectangle (1, 1);
			\draw (0, 0.5) node[left]  {\textit{per}};
			\draw (1, 0.5) node[right] {\textit{per}};
			\foreach [evaluate=\idS as \xS using {(\idS-1)/(\nbnoeuds-1)}] \idS in {1,...,\nbnoeuds} {
				\draw[rougeTerre, semithick] (\xS, 0) -- +({1/tan(\cutangle)}, 1);
				\foreach [evaluate=\idT as \xT using {\xS+(1/tan(\cutangle))*(\idT-1)/(\nbnoeuds-1)},%
									evaluate=\idT as \yT using {(\idT-1)/(\nbnoeuds-1)}] \idT in {1,...,\nbnoeuds} {%
					\fill[rougeTerre] (\xT, \yT) circle [radius=0.5pt];
				}
			}
		\end{tikzpicture}
	\end{subfigure}
	\caption{Two-dimensional mesh for the $2$D method (left), and family of one-dimensional meshes for the quasi-$1$D method (right)\label{fig:meshes_2D_quasi1D}}
\end{figure}

\subsection{A fully two-dimensional method}
\label{sec:methode_2D}
The first method is inspired from the resolution of the elliptic Helmholtz equation (see \cite{flissthese} for instance), and consists in solving directly the local cell problems on an unstructured mesh of the periodicity cell $\mathcal{C}^\diese_0 = (0, 1)^2$ (see Figure \ref{fig:meshes_2D_quasi1D}). 

\vspace{1\baselineskip} \noindent
We start from a triangular mesh $\mathscr{T}_h(\mathcal{C}^\diese_0)$ of $\mathcal{C}^\diese_0 = (0, 1)^2$ with a mesh step $h > 0$. We assume that this mesh is \emph{periodic}, in the sense that one can identify the mesh nodes on the boundary $\yi_i = 0$ with those on $\yi_i = 1$, for $1 \leq i \leq 2$. In particular for $i = 1$, this condition allows us to handle the periodic boundary conditions.

%

\vspace{1\baselineskip} \noindent
Now let $\mathcal{V}_h(\mathcal{C}^\diese_0)$ be the usual $H^1$--conforming approximation by Lagrange finite elements of order $d > 0$. We also introduce
\[\displaystyle
\mathcal{V}_{h, \textit{per}}(\mathcal{C}^\diese_0) := \big\{V \in \mathcal{V}_h(\mathcal{C}^\diese_0) \ /\ \restr{V}{\yi_1 = 0} = \restr{V}{\yi_1 = 1}\big\}
\]
as an internal approximation of $H^1_{\cut, \textit{per}}(\mathcal{C}^\diese_0)$. Finally, to approximate $L^2(\Sigma^\diese_{2, 0})$ and $L^2(\Sigma^\diese_{2, 1})$, we consider the following subspaces:
\[\displaystyle
	\spforall a \in \{0, 1\}, \quad \mathcal{V}_{h, \textit{per}}(\Sigma^\diese_{2, a}) = \big\{ \restr{V_h}{\Sigma^\diese_{2, a}}\ / \ V_h \in \mathcal{V}_{h, \textit{per}}(\mathcal{C}^\diese_0) \big\}.
\]
Since the mesh nodes on $\Sigma^\diese_{2, 0}$ and $\Sigma^\diese_{2, 1}$ can be identified to each other by periodicity of $\mathscr{T}_h(\mathcal{C}^\diese_0)$, we can also make the identification $\mathcal{V}_{h, \textit{per}}(\Sigma^\diese_{2, 0}) \equiv \mathcal{V}_{h, \textit{per}}(\Sigma^\diese_{2, 1}) \equiv \mathcal{V}_{h, \textit{per}}(0, 1)$, as in the continuous case. Set $N := \dim \mathcal{V}_{h, \textit{per}}(0, 1)$, and consider a basis $\smash{(\varphi_p)_{1 \leq p \leq N}}$.



\vspace{1\baselineskip} \noindent
For any data $\varphi_h \in \mathcal{V}_{h, \textit{per}}(0, 1)$, we denote by $E^0_h(\varphi_h), E^1_h(\varphi_h) \in \mathcal{V}_{h, \textit{per}}(\mathcal{C}^\diese_0)$ the solutions of the discrete counterpart of the local cell problems \eqref{eq:local_cell_problem}--\eqref{eq:local_cell_BC} defined in a weak sense.
%
In practice, one has to compute $\smash{E^j_h(\varphi_p)}$, where $\smash{(\varphi_p)_{1 \leq p \leq N}}$ is a basis of $\mathcal{V}_{h, \textit{per}}(0, 1)$.

\vspace{1\baselineskip} \noindent
Similarly to the weak expression \eqref{eq:DtNloc_weak_form} of the continuous local DtN operators, the discrete local DtN operators $\mathcal{T}^{j k}_h \in \mathcal{L}(\mathcal{V}_{h, \textit{per}}(0, 1))$, $j, k = 0, 1$, are defined for any $\varphi_h, \psi_h \in \mathcal{V}_{h, \textit{per}}(0, 1)$ as follows:
\[\displaystyle
\Big\langle \mathcal{T}^{j k}_h \varphi_h, \; \psi_h \Big\rangle = \int_{\mathcal{C}^\diese_0} \left[ \mu_p\; \Dt{} E^j_h(\varphi_h) \; \Dt{} \overline{E^k_h(\psi_h)} - \rho_p\; \omega^2\; E^j_h(\varphi_h)\; \overline{E^k_h(\psi_h)}\, \right].
\]
In practice, these operators are represented as $N \times N$ matrices $\mathbb{T}^{jk}$ whose components are given by $\smash{\mathbb{T}^{jk}_{pq} = \big\langle \mathcal{T}^{j k}_h \varphi_q, \; \varphi_p \big\rangle}$, for $p, q \in \llbracket 1, N \rrbracket$.

\vspace{1\baselineskip} \noindent
Let $\varphi_h \in \mathcal{V}_{h, \textit{per}}(0, 1) \subset \mathscr{C}_{\textit{per}}(\R)$ such that $\varphi_h(0) = 1$. The computation of the propagation operator $\mathcal{P}_h \in \mathcal{L}(\mathcal{V}_{h, \textit{per}}(0, 1))$ is presented in Subsection \ref{sec:discrete_Riccati_equation}. Once this operator is determined, the solution of the half-guide problem \eqref{eq:half_guide_problem} can be approximated with the function defined cell by cell by
\[\displaystyle
\spforall \ell \in \N, \quad \restr{U^+_{\cut, h}(\varphi_h)(\cdot + \ell\, \vec{\ev}_\ord)}{\mathcal{C}^\diese_0} = E^0_h(\mathcal{P}^{\ell}_h\, \varphi_h) + E^1_h(\mathcal{P}^{\ell+1}_h \, \varphi_h).
\]

%
%

\noindent
Finally, a suitable approximation of the solution of the half-line problem \ref{eq:half_line_problem} is provided by
\[\displaystyle
\spforall \xv \in \R, \quad u^+_{\cut, h}(\xv) = U^+_{\cut, h}(\varphi)(\cut\, \xv).
\]

\subsection{A quasi one-dimensional method}
\label{sec:methode_quasi1D}
Though easy to implement, the two-dimensional approach described in the previous section does not exploit the fibered properties of the directional derivative $\Dt{}$. However, the periodic half-guide problem can be seen as a concatenation in a certain sense of one-dimensional half-line problems. This fibered structure is the core of the method presented in this section.

\subsubsection{Presentation}
\label{sec:a_quasi_one_dimensional_method_presentation}
For any $\si \in \R$, we consider the one-dimensional cell problems
\begin{equation}
\left|
\begin{array}{r@{\ }c@{\ }l@{\quad}l}
\displaystyle - \frac{d}{d \xv} \Big( \mu_{\si, \cut} \; \frac{d e^j_{\si, \cut}}{d \xv} \Big) - \rho_{\si, \cut} \; \omega^2 \; e^j_{\si, \cut} &=& 0, \quad \textnormal{in} & (0, 1/\cuti_2) := I_\cut,
\\[6pt]
\displaystyle e^0_{\si, \cut}(0) = 1\ \quad \textnormal{and}\ \quad e^0_{\si, \cut}(1/\cuti_2) &=& 0,
\\[6pt]
\displaystyle e^1_{\si, \cut}(0) = 0\ \quad \textnormal{and}\ \quad e^1_{\si, \cut}(1/\cuti_2) &=& 1.
\end{array}
\right.
\label{eq:QP_local_cell_problems}
\end{equation}
Then, by analogy with Proposition \ref{prop:structure_half_guide}, one easily shows that the local cell problems are concatenations of one-dimensional cell problems, in the following sense.
\begin{prop}\label{prop:structure_local_cell_problem}
	For any boundary data $\varphi$ in $L^2(0, 1)$, the solutions $E^0(\varphi)$ and $E^1(\varphi)$ of the local cell problems \eqref{eq:local_cell_problem} are given by
	\begin{equation}
		\label{eq:concatenation_local_cell_problems}
		\aeforall \yv \in \mathcal{C}^\diese_0, \quad E^j(\varphi)(\yv) = \widetilde{\varphi}\big( \sv_\cut(\yv) + j\, \cuti_1/\cuti_2 \big)\; e^j_{\sv_\cut(\yv), \cut} \bigg( \frac{y_2}{\cuti_2} \bigg),
	\end{equation}
	where $e^j_{\si, \cut}$ denotes the solution of the cell problems \eqref{eq:QP_local_cell_problems}.
\end{prop}

\noindent
Proposition \ref{prop:structure_local_cell_problem} also highlights the structure of the local DtN operators. To see this, let us introduce the \emph{local DtN functions} $t^{jk}_\cut$ defined for $j, k = 0, 1$, by
\begin{equation}\displaystyle
\spforall \si \in \R, \quad t^{jk}_\cut(\si) = (-1)^{k+1} \cuti_2\; \bigg[\mu_{\si, \cut}\; \frac{d e^j_{\si, \cut}}{d \xv}\bigg] \bigg( \frac{j}{\cuti_2} \bigg).
\label{eq:local_DtN_functions}
\end{equation}
Note that by periodicity of $\mu_p$ and $\rho_p$, the maps $s \mapsto e^j_{s, \cut}$ and $t^{jk}_\cut$ are $1$--periodic.

\vspace{1\baselineskip} \noindent
By applying the directional derivative operator $\Dt{}$ to \eqref{eq:concatenation_local_cell_problems}, and by using the relationship between $\Dt{} E^j(\varphi)$ and $d e^j_{s, \cut} / dx$ given by \eqref{eq:derivee_fonction_2}, it follows that the local DtN operators  defined by \eqref{eq:DtNloc} are weighted translation operators, similarly to the propagation operator.
\begin{prop}
The operators $\mathcal{T}^{jk}$ can be written for $\varphi \in L^2(0, 1)$ and $\si \in (0, 1)$ as
\begin{equation}
	\label{eq:DtNloc_wto}
	\begin{array}{l@{\quad}c@{\quad}l}
		\mathcal{T}^{00} \varphi(\si) = t^{00}_\cut(\si)\; \widetilde{\varphi}(\si) &\textnormal{and}& \mathcal{T}^{10} \varphi(\si) = t^{10}_\cut(\si)\; \widetilde{\varphi}(\si + \cuti_1/\cuti_2),
		\\[8pt]
		\mathcal{T}^{11} \varphi(\si) = t^{11}_\cut(\si - \cuti_1/\cuti_2)\; \widetilde{\varphi}(\si) &\textnormal{and}& \mathcal{T}^{01} \varphi(\si) = t^{01}_\cut(\si - \cuti_1/\cuti_2)\; \widetilde{\varphi}(\si - \cuti_1/\cuti_2),
	\end{array}
\end{equation}
where we recall that $\widetilde{\varphi}$ denotes the periodic extension of $\varphi$ on $\R$, defined by \eqref{eq:per_extension}.
\end{prop}

\vspace{1\baselineskip} \noindent
Finally, the solution $u^+_\cut$ of the half-line problem \eqref{eq:half_line_problem} can be computed directly from the functions $e^j_{\sv, \cut}$ and from the propagation operator. In fact, given a function $\varphi \in \mathscr{C}_{\textit{per}}(\Sigma^\diese_{\ord, 0})$ such that $\varphi(0) = 1$, taking formally the trace along $\cut\, \R$ in \eqref{eq:UfromEis} leads to
\begin{equation}
	\spforall \ell \in \N, \quad \restr{u^+_\cut(\cdot + \ell/\cuti_2)}{I_\cut} = (\widetilde{\mathcal{P}^{\ell} \varphi})(\ell\, \cuti_1/\cuti_2)\; e^0_{\ell \cuti_1/\cuti_2, \cut} + (\widetilde{\mathcal{P}^{\ell+1} \varphi})((\ell+1)\, \cuti_1/\cuti_2)\; e^1_{\ell \cuti_1/\cuti_2, \cut}.
	\label{eq:ufromQPeis}
\end{equation}
The proof of this result is similar to those of \eqref{eq:UfromEis} and Proposition \ref{prop:periodic_structure_solution}.

\vspace{1\baselineskip} \noindent
Expressions \eqref{eq:concatenation_local_cell_problems},  \eqref{eq:DtNloc_wto}, and \eqref{eq:ufromQPeis} form the basis of the \emph{quasi one-dimensional} or \emph{quasi-1D} strategy, which consists in approximating the solutions $e^j_{\sv, \cut}$ as well as the functions $t^{jk}_\cut$ and finally the local DtN operators $\mathcal{T}^{jk}$. Then once the propagation operator is computed by solving the constrained Riccati equation \eqref{eq:Riccati}, the solution $u^+_\cut$ may be constructed directly cell by cell using \eqref{eq:ufromQPeis}.
\subsubsection{Discretization}
The quasi-1D approach requires two distinct approximate spaces associated to the transverse and the $\cut$--oriented directions (see Figure \ref{fig:meshes_2D_quasi1D}).

\paragraph{Transverse direction.} We begin with a one-dimensional mesh $\mathscr{T}_h(0, 1)$ of $\Sigma^\diese_{2, 0} \equiv (0, 1)$ with a mesh step $h > 0$. Let $\mathcal{V}_h(0, 1)$ be the approximation space of $H^1(0, 1)$ by Lagrange finite elements of order $d > 0$. We denote by $\smash{(\varphi_p)_{0 \leq p \leq N}}$ the usual nodal basis, which satisfies in particular $\smash{\varphi_p(s_q) = \delta_{p,q}}$, where $(s_p)_{0 \leq p \leq N}$ are points (including the mesh vertices) such that $0 = s_0 < \dots < s_N = 1$. %
Then an internal approximation of $L^2(0, 1)$ is
\[\displaystyle
\mathcal{V}_{h, \textit{per}}(0, 1) := \vect\{ \varphi_0 + \varphi_N, \varphi_1, \dots , \varphi_{N-1}\},
\]
which is chosen so that $\mathcal{V}_{h, \textit{per}}(0, 1) \subset \mathscr{C}_{\textit{per}}(0, 1)$. In particular, from the definition of the basis functions $\varphi_i$, one has the following decomposition
\begin{equation}\displaystyle
\spforall\! \varphi_h \in \mathcal{V}_{h, \textit{per}}(0, 1), \quad \varphi_h = \sum_{p = 0}^{N} \varphi_h(s_p)\, \varphi_p, \quad \textnormal{with} \quad \varphi_h(s_0) = \varphi_h(s_N).
\label{eq:decomposition_Sigma0}
\end{equation}

\paragraph{$\cut$--oriented direction.} Let $\mathscr{T}_{h_\cut}(I_\cut)$ denote a mesh of the line segment $I_\cut$ with a mesh step $h_\cut > 0$. Set $\mathcal{V}_{h_\cut}(I_\cut)$ as the approximation space of $H^1(I_\cut)$ by Lagrange finite elements of order $d_\cut > 0$ and define $\mathcal{V}_{h_\cut, 0}(I_\cut) := \mathcal{V}_{h_\cut}(I_\cut) \cap H^1_0(I_\cut)$.
%
%
%


\vspace{1\baselineskip} \noindent
The approximation of $e^0_{s, \cut}$ and $e^1_{s, \cut}$ can be seen as a two-step process. First, for any $s \in \R$, consider the solution $e^j_{s, \cut, h_\cut}$ of the discrete variational formulation associated to \eqref{eq:QP_local_cell_problems}. 

\vspace{1\baselineskip} \noindent
In practice, the solution $e^j_{s, \cut, h_\cut}$ can only be computed for a finite number of $s \in (0, 1)$. This is where the discretization in the transverse direction comes into play: given $\xv \in I_\cut$, the function $\smash{\sv \mapsto e^j_{s, \cut, h_\cut}(\xv)}$ may be interpolated in $\mathcal{V}_{h, \textit{per}}(0, 1)$.

The interpolation process requires to compute the discrete solution $e^j_{s, \cut, h_\cut}$ only for $s = s_p$, $p \in \llbracket 0 , N-1 \rrbracket$. Then, using the decomposition formula \eqref{eq:decomposition_Sigma0}, $e^j_{\sv, \cut}$ shall be approximated by
\begin{equation}
	\displaystyle
	\spforall (s, \xv) \in (0, 1) \times I_\cut, \quad e^j_{s, \cut, \hv}(\xv) = \sum_{p = 0}^N e^j_{s_p, \cut, h_\cut}(\xv)\; \varphi_p(s), \quad \textnormal{with} \quad \hv = (h, h_\cut).
\end{equation}
where $e^j_{0, \cut, h_\cut} = e^j_{1, \cut, h_\cut}$ (because $e^j_{s, \cut}$ is $1$--periodic with respect to $s$).
%

From the solutions $e^j_{s, \cut, \hv}$, we introduce the discrete local DtN functions
\[\displaystyle
\displaystyle
\spforall s \in (0, 1), \quad	t^{jk}_{\cut, \hv}(s) = \cuti_\ord\; \int_0^{1/\cuti_\ord} \Big( \mu_{s, \cut} \; \frac{d e^j_{s, \cut, \hv}}{d \xv} \; \overline{\frac{d e^k_{s, \cut, h}}{d \xv}} - \rho_{s, \cut} \; \omega^2 \; e^j_{s, \cut, h} \; \overline{e^k_{s, \cut, \hv}} \; \Big),
\]
which are inspired from the weak expression \eqref{eq:local_DtN_functions} of the local DtN functions $t^{jk}_\cut$. Then, by analogy with \eqref{eq:DtNloc_wto}, we define the discrete DtN operators $\mathcal{T}^{jk}_\hv \in \mathcal{L}(\mathcal{V}_{h, \textit{per}}(0, 1))$ for any $\varphi_h$, $\psi_h \in \mathcal{V}_{h, \textit{per}}(0, 1)$ as follows:
\begin{equation}
	\displaystyle \Big\langle \mathcal{T}^{jk}_\hv \varphi_h,\; \psi_h \Big\rangle = \int_0^1 t^{jk}_{\cut, \hv}(s - k\, \cuti_1/\cuti_2)\; \varphi_h(s + (j - k)\,\cuti_1/\cuti_2) \; \overline{\psi_h(s)}\; d s.
	\label{eq:DtNloc_wto_weak_form}
\end{equation}
%
%
These discrete DtN operators, when computed for $\varphi_h$, $\psi_h$ being the basis functions of $\mathcal{V}_{h, \textit{per}}(0, 1)$, are represented as $N \times N$ matrices, where $N = \dim \mathcal{V}_{h, \textit{per}}(0, 1)$. The integrals which appear in \eqref{eq:DtNloc_wto_weak_form} are evaluated in practice using a specifically designed quadrature rule whose description is omitted here.

%
\vspace{1\baselineskip} \noindent
Finally, let $\varphi_h \in \mathcal{V}_{h, \textit{per}}(0, 1) \subset \mathscr{C}_{\textit{per}}(\R)$ such that $\varphi_h(0) = 1$. Then using \eqref{eq:ufromQPeis}, the solution of the half-line problem \eqref{eq:half_line_problem} can be approximated with the function defined cell by cell by
\[\displaystyle
\spforall \ell \in \N, \quad \restr{u^+_{\cut, \hv}(\cdot + \ell/\cuti_2)}{I_\cut} = (\mathcal{P}^{\ell}_\hv \varphi_h)(\ell\, \cuti_1/\cuti_2)\; e^0_{\ell \cuti_1/\cuti_2, \cut, \hv} + (\mathcal{P}^{\ell+1}_\hv \varphi_h)((\ell+1)\, \cuti_1/\cuti_2)\; e^1_{\ell \cuti_1/\cuti_2, \cut, \hv}.
\]
where $\mathcal{P}_\hv \in \mathcal{L}(\mathcal{V}_{h, \textit{per}}(0, 1))$ corresponds to a suitable discrete $\R^{N \times N}$ approximation of $\mathcal{P}$. The computation of such an operator is the subject of the next subsection.

\subsection{Approximation of the propagation operator}
\label{sec:discrete_Riccati_equation}
In order to find a suitable approximation $\mathcal{P}_h \in \mathcal{L}(\mathcal{V}_{h, \textit{per}}(0, 1))$ of the propagation operator $\mathcal{P}$, it is natural to introduce the discrete constrained Riccati equation
\begin{equation}\label{eq:Riccati_discrete}
\left|
\begin{array}{l}
	\textnormal{\textit{Find $\mathcal{P}_h \in \mathcal{L}(\mathcal{V}_{h, \textit{per}}(0, 1))$ such that $\rho(\mathcal{P}_h) < 1$ and $\displaystyle \mathcal{T}_h(\mathcal{P}_h) = 0$, where}}
	\\[12pt]
	\multicolumn{1}{c}{\displaystyle \mathcal{T}_h(\mathcal{P}_h) := \mathcal{T}^{10}_h\mathcal{P}^2_h + (\mathcal{T}^{00}_h + \mathcal{T}^{11}_h)\, \mathcal{P}_h + \mathcal{T}^{01}_h,}
\end{array}
\right.
\end{equation}
and where $(\mathcal{T}^{00}_h, \mathcal{T}^{01}_h, \mathcal{T}^{10}_h, \mathcal{T}^{11}_h)$ are obtained via one of the methods described in Sections \ref{sec:methode_2D} and \ref{sec:methode_quasi1D}. Using the same arguments as for the elliptic Helmholtz equation \cite{flissthese}, it can be proved that this discrete equation admits a unique solution.

\vspace{1\baselineskip} \noindent
In order to solve \eqref{eq:Riccati_discrete}, two methods have been proposed in \cite{jolyLiFliss}: a \emph{spectral decomposition method}, and a \emph{modified Newton method}. Here, we only describe the spectral approach.

\vspace{1\baselineskip} \noindent
The spectral decomposition method consists in characterizing $\mathcal{P}_h$ by means of its eigenpairs $(\lambda_i, \psi_i)$ of $\mathcal{P}_h$. Doing so however raises an important question: is $\mathcal{P}_h$ completely defined by its eigenpairs? This is equivalent to wondering if $\mathcal{P}_h$ is diagonalizable or not. The diagonalizability of $\mathcal{P}_h$ is an open question, but for the sake of simplicity, we will assume in the sequel that this is the case, namely
\[\displaystyle
\textnormal{The family of eigenfunctions $(\psi_i)_{1 \leq i \leq N}$ forms a basis of $\mathcal{V}_{h, \textit{per}}(0, 1)$.}
\]
In practice, this is the situation that we always have encountered. Moreover, in the case where this assumption fails to be true, one can still adapt the method, and recover $\mathcal{P}_h$ through a Jordan decomposition. (See \cite[Section 2.3.2.3]{flissthese} for more details.)

\vspace{1\baselineskip} \noindent
The spectral approach relies on the results presented in Section \ref{sec:about_Riccati_equation}, which remain true for the discrete equation. In particular, by analogy with Proposition \ref{prop:spectrum_P_riccati}, $(\lambda_h, \psi_h) \in \C \times \mathcal{V}_{h, \textit{per}}(0, 1)$ is an eigenpair of $\mathcal{P}_h$ if and only if it satisfies
\[\displaystyle
\mathcal{T}_h(\lambda_h)\, \psi_h = 0, \quad \textnormal{with} \quad \psi_h \neq 0 \quad \textnormal{and} \quad |\lambda_h| < 1.
\]
Solving the Riccati equation hence comes down to solving a quadratic eigenvalue problem:
\begin{equation}\label{eq:quadratic_EVP}
\left|
\begin{array}{l}
	\textnormal{\textit{Find $(\lambda_h, \psi_h) \in \C \times \mathcal{V}_{h, \textit{per}}(0, 1)$ such that $\psi_h \neq 0$, $|\lambda_h| < 1$ and}}
	\\[12pt]
	\multicolumn{1}{c}{\displaystyle \lambda^2_h\, \mathcal{T}^{10}_h\psi_h + \lambda_h\, (\mathcal{T}^{00}_h + \mathcal{T}^{11}_h) \psi_h + \mathcal{T}^{01}_h \psi_h = 0.}
\end{array}
\right.
\end{equation}
If one sets $N = \dim \mathcal{V}_{h, \textit{per}}(0, 1)$, then \eqref{eq:quadratic_EVP} can be reduced to a $2N \times 2N$ linear eigenvalue problem, thus yielding $2N$ eigenvalues. In order to pick the $N$ eigenvalues of the propagation operator, we need a criterion. To do so, note that with the 2D or the quasi-1D method, the properties of the local DtN operators (Proposition \ref{prop:pties_local_DtN}) remain preserved for the discrete operators $\mathcal{T}^{jk}_h$. Hence Proposition \ref{prop:pairs_Riccati} admits the following discrete version:
\[
\Ker \mathcal{T}_h(\lambda_h) \neq \{0\}  \quad \Longleftrightarrow \quad \Ker \mathcal{T}_h(1/\lambda_h) \neq \{0\}.
\]
Therefore, as already expected with Remark \ref{rmk:pairs_Riccati}, the solutions of \eqref{eq:quadratic_EVP} can be grouped into pairs $(\lambda_h, 1/\lambda_h)$, where $0 < |\lambda_h| < 1$. Consequently, in order to compute $\mathcal{P}_h$, one can solve \eqref{eq:quadratic_EVP} (using for instance linearization techniques), and choose the $N$ eigenpairs $(\lambda_h, \psi_h)$ which satisfy $|\lambda_h| < 1$.

\subsection{The DtN coefficient}\label{sec:discrete_DtN_coefficient}
Finally, consider a function $\varphi_h \in \mathcal{V}_{h, \textit{per}}(0, 1) \subset \mathscr{C}_{\textit{per}}(\R)$ which satisfies $\varphi_h(0) = 1$. Then by analogy with \eqref{eq:DtN_operator}, and in the spirit of Proposition \ref{prop:DtNcoeff_from_DtNoperator}, we define the discrete DtN operator and the discrete DtN coefficient as follows:
\[\displaystyle
 \Lambda_h = \mathcal{T}^{10}_h \mathcal{P}_h + \mathcal{T}^{00}_h \quad \textnormal{and} \quad \lambda^+_h = \frac{(\Lambda_h \varphi_h)(0)}{\cuti_2},
\]
where $\mathcal{T}^{10}_h$ and $\mathcal{T}^{00}_h$ are computed using one of the methods presented in Sections \ref{sec:methode_2D} and \ref{sec:methode_quasi1D}, and where $\mathcal{P}_h$ is the solution of the discrete Riccati equation \eqref{eq:Riccati_discrete}.


\subsection{Numerical results}
\label{sec:numerical_results}
\begin{figure}[H]
	\centering
	\begin{tikzpicture}
		\begin{axis}[
				width=0.99\textwidth,
				height=3cm,
				xmin=-6, xmax=6,
				ymin=0, ymax=2.75,
				xtick = {-6, -4, -2, 0, 2, 4, 6},
				ytick = {1, 2},
				enlargelimits=false,
				axis on top,
				legend style = {legend pos=south west, fill=none, draw=none},
			]
			\addplot+[domain=-6:-1, smooth, no marks, samples=200, color=bleuDoux!75, thick] {1.5 + cos(deg(pi*x)) * cos(deg(pi*sqrt(3)*x))};
			\draw[thick, color=bleuDoux!75] (axis cs:-1, 0.8339) -- (axis cs:-1/3, 0.8339) -- (axis cs:-1/3, 2) -- (axis cs:1/3, 2) -- (axis cs:1/3, 0.8339) -- (axis cs:1, 0.8339);
			\addplot+[domain=1:6, smooth, no marks, samples=200, color=bleuDoux!75, thick] {1.5 + cos(deg(pi*x)) * cos(deg(pi*sqrt(3)*x))};
			\draw[dashed, thick] (axis cs:-1, 0) -- (axis cs:-1, 2.75);
			\draw[dashed, thick] (axis cs:+1, 0) -- (axis cs:+1, 2.75);
			\legend{$\mu$}
		\end{axis}
	\end{tikzpicture}

	\begin{tikzpicture}
		\begin{axis}[
				width=0.99\textwidth,
				height=3cm,
				xmin=-6, xmax=6,
				ymin=0, ymax=2.75,
				xtick = {-6, -4, -2, 0, 2, 4, 6},
				ytick = {1, 2},
				enlargelimits=false,
				axis on top,
				legend style = {legend pos=south west, fill=none, draw=none},
			]
			\addplot+[domain=-6:-1, smooth, no marks, samples=200, color=rougeTerre, thick] {1.5 + 0.5*sin(deg(pi*x)) + 0.5*sin(deg(pi*sqrt(3)*x))};
			%
			\draw[thick, color=rougeTerre] (axis cs:-1, 1.8729) -- (axis cs:0, 1.8729) -- (axis cs:0, 1.1271) -- (axis cs:1, 1.1271);
			\addplot+[domain=1:6, smooth, no marks, samples=200, color=rougeTerre, thick] {1.5 + 0.5*sin(deg(pi*x)) + 0.5*sin(deg(pi*sqrt(3)*x))};
			\draw[dashed, thick] (axis cs:-1, 0) -- (axis cs:-1, 2.75);
			\draw[dashed, thick] (axis cs:+1, 0) -- (axis cs:+1, 2.75);
			\legend{$\rho$}
		\end{axis}
	\end{tikzpicture}
	\begin{tikzpicture}
		\begin{axis}[
				width=0.99\textwidth,
				height=3cm,
				xmin=-6, xmax=6,
				xtick = {-6, -4, -2, 0, 2, 4, 6},
				ytick = {0.5},
				enlargelimits=false,
				axis on top,
				legend style = {legend pos=south west, fill=none, draw=none},
			]
			\addplot+[domain=-0.99:0.99, smooth, no marks, samples=200, color=teal, thick] {exp(100 - 100/((1+x)*(1-x)))};
			\draw [teal, thick] (axis cs:-6,0)--(axis cs:-0.99,0);
			\draw [teal, thick] (axis cs:6,0)--(axis cs:0.99,0);
			\legend{$f$}
			\draw[dashed, thick] (axis cs:-1, 0) -- (axis cs:-1, 1);
			\draw[dashed, thick] (axis cs:+1, 0) -- (axis cs:+1, 1);
		\end{axis}
	\end{tikzpicture}
	\caption{The locally perturbed quasiperiodic coefficients $\mu$ and $\rho$, and the source term $f$. \label{fig:coefficients_mu_rho}}
\end{figure}
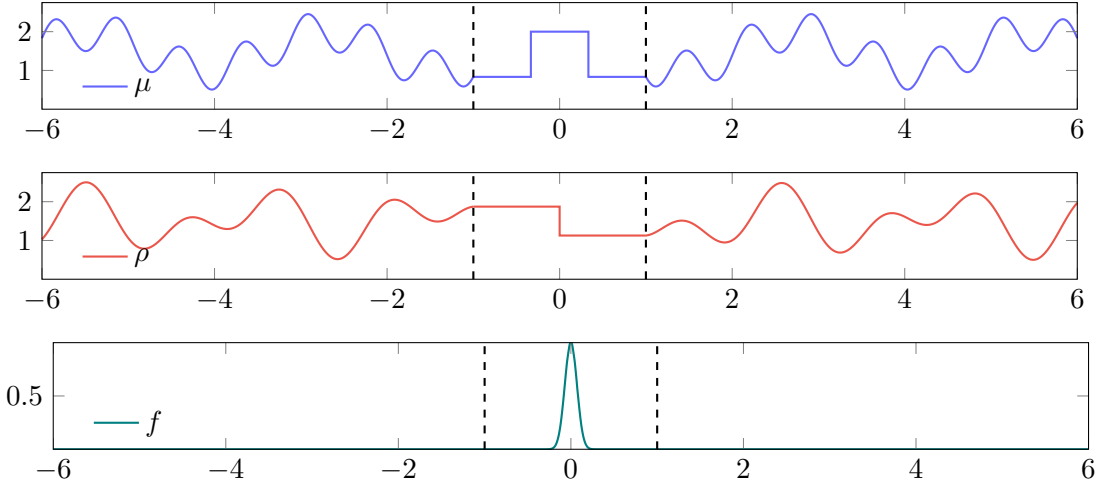

We present some numerical results to validate the method, to illustrate its efficiency, and to compare the multi-dimensional and the quasi one-dimensional methods in the case where the order of quasiperiodicity is set to $\ord = 2$. Simulations will be carried out with the periodic coefficients $\mu_p$ and $\rho_p$, defined for $\yv = (\yi_1, \yi_2) \in \R^2$ by
\[\displaystyle
\mu_p(\yv) = 1.5 + \cos (2\pi\yi_1)\, \cos (2\pi\yi_2) \quad \textnormal{and} \quad \rho_p(\yv) = 1.5 + 0.5\, \sin (2\pi\yi_1) + 0.5\, \sin (2\pi\yi_2).
\]
We set $\cut = (\cos\pi/3, \sin\pi/3)$. As the ratio $\cuti_2 / \cuti_1 = \sqrt{3}$ is irrational, $\cut$ is an irrational vector. For $a = 1$, the source term $f$ is the cut-off function
\[\displaystyle
\spforall \xv \in \R, \quad	f(\xv) = \exp \left( 100\, \big(1 - 1 / (1 - \xv^2)\big) \right)\; \chi_{(-1, 1)},
\]
and the local perturbations $\mu_i$ and $\rho_i$ are defined as piecewise constants, so that the coefficients $\mu$ and $\rho$ of the model problem \eqref{eq:whole_line_problem} are represented in Figure \ref{fig:coefficients_mu_rho}.
\noindent
\subsubsection{The half-line and the half-guide solutions}\label{sec:results:half_line_guide}
The model problem \eqref{eq:whole_line_problem} is solved by computing the solutions of the half-line problems \eqref{eq:half_line_problems_0}, as well as the DtN coefficients $\lambda^\pm$. In this part, only results regarding the numerical resolution of the problem \eqref{eq:half_line_problem} are going to be presented, as the problem set on $(-\infty, -a)$ provides the same overall results.

\paragraph{Error analysis}
In order to validate the method, we introduce for $L > 0$ the unique function $u^+_{\cut, L}$ in $H^1(0, L)$ that satisfies Problem \eqref{eq:half_line_problem} on the truncated domain $(0, L)$, with $u^+_{\cut, L}(L) = 0$. Similarly, define $\Omega_L := (0, 1)^{\ord-1} \times (0, L)$, and for any $\varphi \in L^2(\Sigma^\diese_{\ord, 0})$, let $U^+_{\cut, L}(\varphi) \in H^1_{\cut}(\Omega_L)$ denote the unique function that satisfies \eqref{eq:half_guide_problem} on $\Omega_L$, with $\restr{U^+_{\cut, L}(\varphi)}{y_2 = L} = 0$.

\vspace{1\baselineskip} \noindent
In presence of absorption, the solutions $u^+_\cut$ and $U^+_\cut(\varphi)$ decay exponentially at infinity (see \eqref{eq:exp_decay_halfline_s} and \eqref{eq:exp_decay_halfguide}), and by studying the problems satisfied by $u^+_{\cut, L} - u^+_\cut$ and $U^+_{\cut, L}(\varphi) - U^+_\cut(\varphi)$, it can be proved as in \cite{fliss_giovangigli} that there exist constants $\alpha, c > 0$ such that for any $L > 0$,
\begin{equation}
	\label{eq:estimate_reference_solution}
	\begin{array}{c}
		\|u^+_{\cut, L} - u^+_\cut\|_{H^1(0, L)} \leq c \,e^{-\alpha \Imag \omega L}\; \|u^+_\cut\|_{H^1(0, L)}
		\\[10pt]
		\|U^+_{\cut, L}(\varphi) - U^+_\cut(\varphi)\|_{H^1_\cut(\Omega_L)} \leq c \, e^{-\alpha \Imag \omega L}\; \|U^+_\cut(\varphi)\|_{H^1_\cut(\Omega_L)}.
	\end{array}
\end{equation}
with $\alpha = \sqrt{\rho_- / \mu_+}$. In particular, if $L$ is chosen large enough, then $u^+_{\cut, L}$ and $U^+_{\cut, L}(\varphi)$ can be viewed as suitable approximations of $u^+_\cut$ and $U^+_\cut(\varphi)$, and thus can serve as reference solutions. In the upcoming results, to make the truncation errors in \eqref{eq:estimate_reference_solution} negligible with respect to the errors induced by the numerical method, we choose $L$ so that
\begin{equation}
	\label{eq:criterion_length_truncated_domain}
	\exp\big(-\sqrt{\rho_- / \mu_+}\, \Imag \omega\, L\big) \leq 10^{-10}.
\end{equation}
The corresponding solutions $u^+_{\cut, L}$ and $U^+_{\cut, L}(\varphi)$, which will be denoted by $u^+_{\textit{ref}}$ and $U^+_{\textit{ref}}(\varphi)$ respectively, are computed via $\mathbb{P}^1$ Lagrange finite elements, with a mesh step $h = 5 \times 10^{-4}$.

\begin{figure}[ht!]
	\centering
	\begin{tikzpicture}
		\begin{groupplot}[
			group style={
				group name=myplot,
		    group size=2 by 1,
		    horizontal sep=0.85cm, 
				ylabels at=edge left,
		    yticklabels at=edge left,
			},
			width=0.4\textwidth, height=4.5cm,
			scale only axis,
			ymin = 2.46594e-04, ymax=1.07008e+00,
			ylabel={Relative errors},
			xmode=log, ymode=log,
			xtick={32, 64, 128, 256, 512},
			xticklabels = {32, 64, 128, 256, 512},
			enlargelimits=false,
			axis on top,
			legend style = {legend pos = south west, fill=none, legend columns = 1},
			]
			\nextgroupplot[legend to name={CommonLegend}, legend style={legend columns=2}]
			%
			\addlegendimage{mark=*, bleuDoux}
			\addlegendimage{mark=square*, rougeTerre}
			\addlegendentry{2D method}
			\addlegendentry{Quasi-1D method}
			\plotcurveandregression{\erreurs}{nDOFs}{quasiD_omega1}{mark=square*, rougeTerre}{dashed, rougeTerre}{0.5}{|-}{pos=0.25, left}{\slopeQuasiDa}
			\plotcurveandregression{\erreurs}{nDOFs}{multiD_omega1}{mark=*, bleuDoux}{dashed, bleuDoux}{0.25}{-|}{pos=0.75, right}{\slopeMultiDa}
			\nextgroupplot[]
			\plotcurveandregression{\erreurs}{nDOFs}{quasiD_omega2}{mark=square*, rougeTerre}{dashed, rougeTerre}{0.5}{|-}{pos=0.25, left}{\slopeQuasiDb}
			\plotcurveandregression{\erreurs}{nDOFs}{multiD_omega2}{mark=*, bleuDoux}{dashed, bleuDoux}{0.25}{-|}{pos=0.75, right}{\slopeMultiDb}
			%
			%
			%
			%
		\end{groupplot}

		\draw ($0.5*(myplot c1r1.south)+0.5*(myplot c2r1.south) - (0, 0.5)$) node[below] {Discretization parameter $1/h$};
		\draw ($0.5*(myplot c1r1.south)+0.5*(myplot c2r1.south) - (0, 1.1)$) node[below] {\ref{CommonLegend}};
		\node [text width=3cm, above=6pt] at (myplot c1r1.north) {\subcaption{$\omega = 8 + 0.25\,\icplx$\label{fig:erreurs:a}}};
		\node [text width=3cm, above=6pt] at (myplot c2r1.north) {\subcaption{$\omega = 20 + 0.25\,\icplx$\label{fig:erreurs:b}}};
	\end{tikzpicture}
	\caption{Relative error in $H^1$ norm of the half-line solution for different values of $\omega$.\label{fig:erreurs}}
\end{figure}
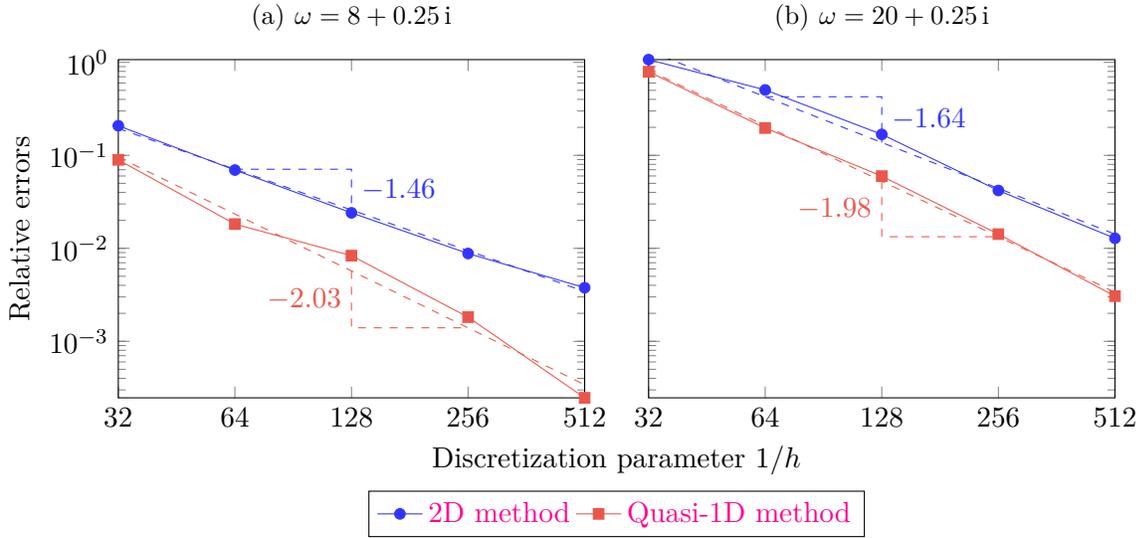

\vspace{1\baselineskip}\noindent
In the following, the boundary data is fixed to $\varphi = 1$, and is omitted in the notation of $U^+_\cut$ and $U^+_{\textit{ref}}$. Also, we only plot relative errors corresponding to the $1$D solution, as the errors for the $2$D solution behave similarly. In Figure \ref{fig:erreurs}, the relative error
\begin{equation}
	\displaystyle \label{eq:relative_error}
	\varepsilon(u^+_\cut) := \frac{\|u^+_{\cut, h} - u^+_{\textit{ref}}\,\|_{H^1(0, 4/\cuti_2)}}{\|u^+_{\textit{ref}}\,\|_{H^1(0, 4/\cuti_2)}}
\end{equation}
is represented with respect to the mesh step $h$, and for both the $2$D and the quasi-$1$D method (with $h_\cut = h$ for the quasi-1D method). The solutions are computed using Lagrange finite elements of degree $1$. 

\vspace{1\baselineskip} \noindent
One sees that the errors tend to $0$ as $h$ at least, as expected for $\mathbb{P}^1$ Lagrange finite elements. With the quasi-$1$D method however, $\varepsilon(u^+_\cut)$ behaves as $h^2$. This is a special superconvergence phenomenon, which is probably due to the fact that the problems solved in practice with the quasi-1D method are one-dimensional. %
%
%
Note also that in general, the quasi-1D method appears to be more accurate than the 2D method. 


\vspace{1\baselineskip} \noindent
For a fixed mesh step, the relative error increases with the real frequency $\Real \omega$. This is a well-known particularity of the Helmholtz equation: since $\Real \omega$ represents the spatial frequency of the time-harmonic waves, the discretization parameter $h$ has to be adapted in order to take their oscillations into account.

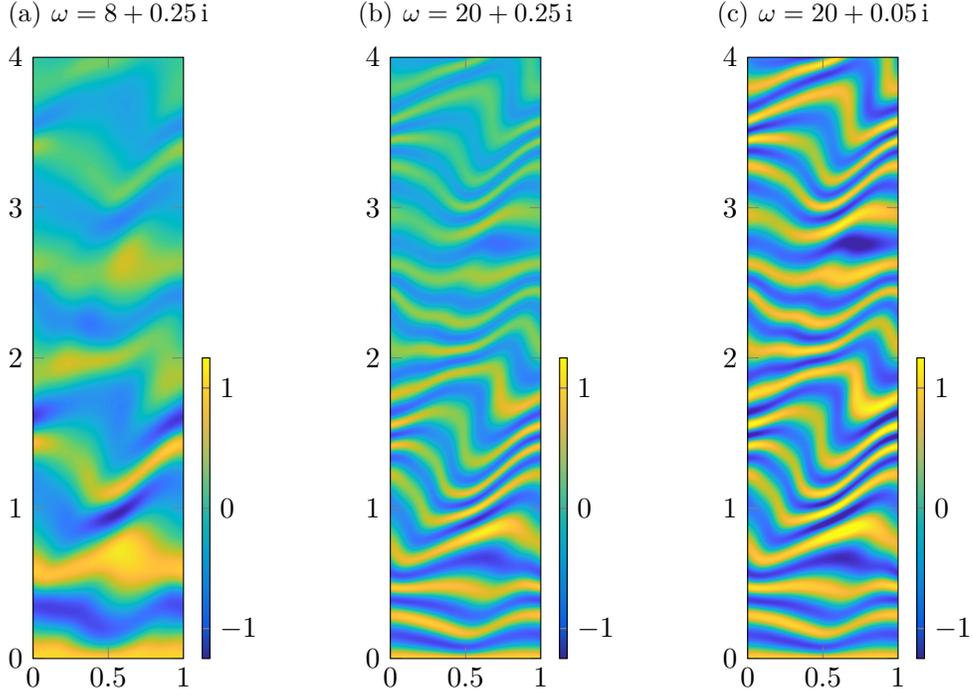
\begin{figure}[H]
	\centering
	\begin{tikzpicture}
		\begin{groupplot}[
			group style={
				group name=half_guide_omega_plot,
				group size=3 by 1,
				horizontal sep=2.75cm,
			},
			enlargelimits=false,
			axis on top,
			width=2cm, height=8cm,
			scale only axis,
			xtick = {0, 0.5, 1}, ytick = {0, 1, 2, 3, 4},
			ticklabel style = {scale=1},
			disabledatascaling,
			colorbar,
			point meta min=-1.25, point meta max=1.25,
			colorbar style={
			yticklabel style={
			/pgf/number format/precision=2,
			/pgf/number format/fixed,
			},
			ytick={-1, 0, 1},
			ticklabel style = {scale=1},
			major tick length=0.05*\pgfkeysvalueof{/pgfplots/parent axis width},
			scaled y ticks=false,
			at = {(2.25cm, 0)},
			anchor=south west,
			width  = 0.05*\pgfkeysvalueof{/pgfplots/parent axis  width},
			height = 0.50*\pgfkeysvalueof{/pgfplots/parent axis height},
			},
			]
			\nextgroupplot[]
			\addplot graphics [xmin=0, xmax=1, ymin=0, ymax=4]{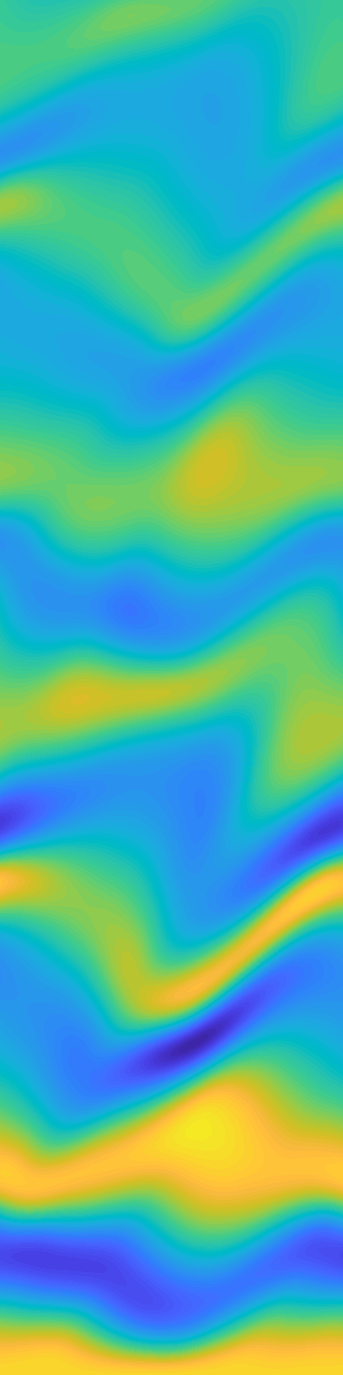};
			\nextgroupplot[]
			\addplot graphics [xmin=0, xmax=1, ymin=0, ymax=4]{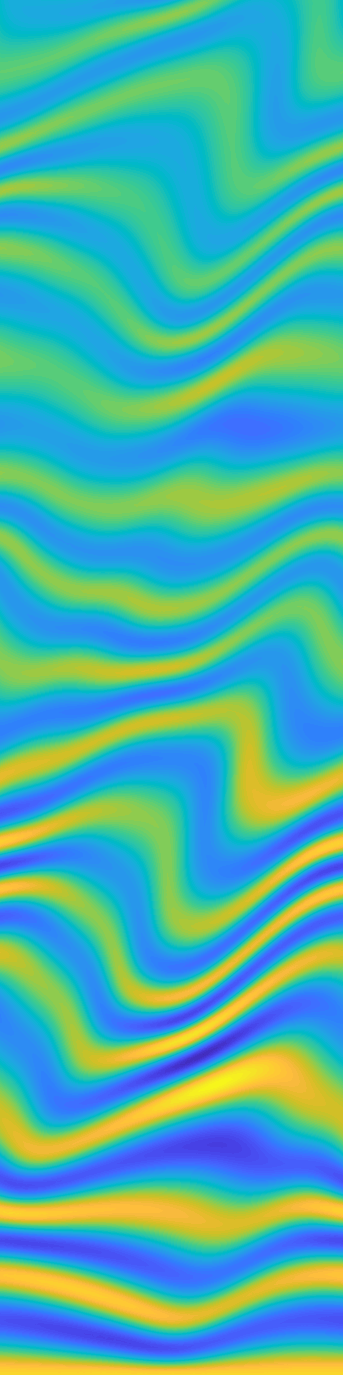};
			\nextgroupplot[]
			\addplot graphics [xmin=0, xmax=1, ymin=0, ymax=4]{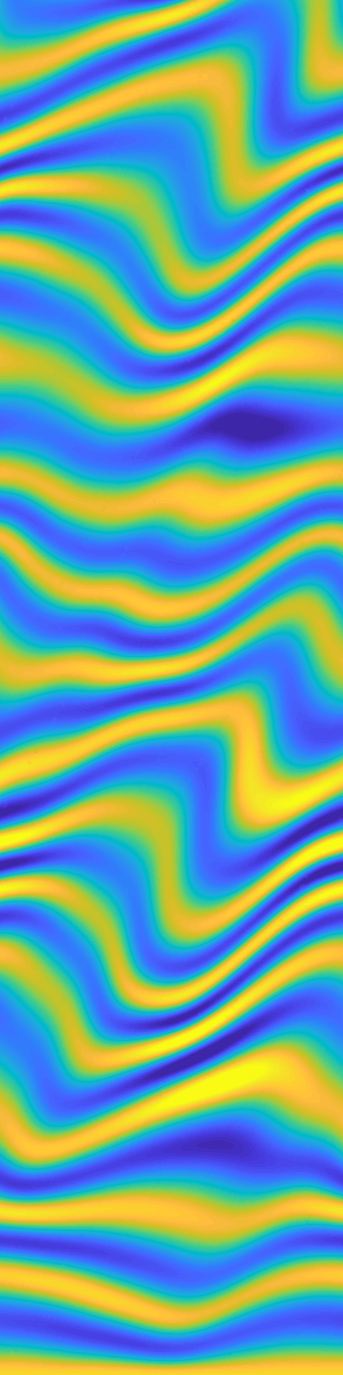};
		\end{groupplot}
		\node [text width=3cm, above=6pt] at (half_guide_omega_plot c1r1.north) {\subcaption{$\omega = 8 + 0.25\,\icplx$\label{fig:half_guide_solution_omega:a}}};
		\node [text width=3cm, above=6pt] at (half_guide_omega_plot c2r1.north) {\subcaption{$\omega = 20 + 0.25\,\icplx$\label{fig:half_guide_solution_omega:b}}};
		\node [text width=3cm, above=6pt] at (half_guide_omega_plot c3r1.north) {\subcaption{$\omega = 20 + 0.05\,\icplx$\label{fig:half_guide_solution_omega:c}}};
	\end{tikzpicture}
	\caption{Real part of the half-guide solution computed using the quasi-1D approach, with $\mathbb{P}^1$ Lagrange finite elements and $h = 2 \times 10^{-3}$, and for different values of $\omega$.\label{fig:half_guide_omega}}
\end{figure}

\paragraph{Representation of the half-guide solution} %
The half-guide solution is represented in Figure \ref{fig:half_guide_omega} for different values of $\omega$, when $\varphi = 1$.

\paragraph{Dependence with respect to the boundary data} The goal of this part is to see how the half-line and the half-guide solutions depend on the boundary data $\varphi$. To do so, we choose three different datas:
\begin{equation}
	\varphi_1(s) = 1, \quad \varphi_2(s) = \cos (2\pi s), \quad \textnormal{and} \quad %
	\varphi_3(s) = 1 - \mathbbm{1}_{[1/3, 2/3]}(s). %
\end{equation}
 We set $\omega = 8 + 0.25\,\icplx$, and we display results obtained with the quasi-1D method, knowing that the 2D method yields the same conclusions. The computations are carried out using $\mathbb{P}^1$ Lagrange finite elements, with mesh steps $h = h_\cut = 2 \times 10^{-3}$.

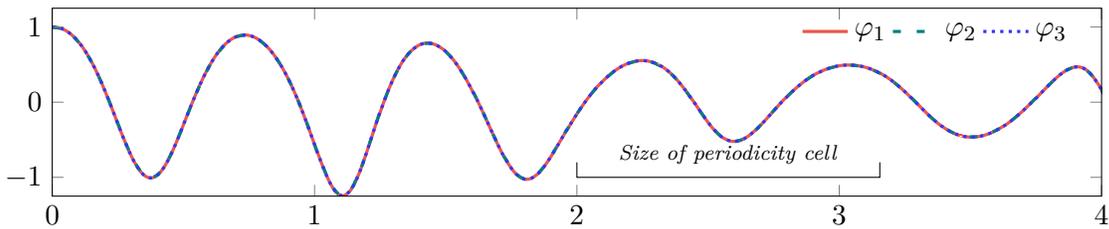
\begin{figure}[H]
	\centering
	\begin{tikzpicture}
		\centering
		\begin{axis}[
			width=0.9\textwidth, height=2.5cm,
			ymin=-1.25, ymax=1.25,
			xmin=0, xmax=4,
			xtick={0, 1, 2, 3, 4},
			scale only axis,
			ticklabel style = {scale=1},
			enlargelimits=false,
			axis on top,
			legend style = {draw=none, fill=none, legend columns = 3, nodes={scale=1, transform shape}},
			]
			\ifthenelse{\boolean{afficherGraphes}}{
			\addplot [color=rougeTerre, smooth, very thick] table [x={x}, y={u_real_phi1}] {\uPhiQuasiD};
			\addplot [color=teal, smooth, very thick, loosely dashed] table [x={x}, y={u_real_phi3}] {\uPhiQuasiD};
			\addplot [color=bleuDoux, smooth, very thick, dotted] table [x={x}, y={u_real_phi2}] {\uPhiQuasiD};
			\draw (axis cs:2, -0.8) -- (axis cs:2, -1) -- (axis cs:3.1547, -1) -- (axis cs:3.1547, -0.8);
			\draw (axis cs:2.5774, -0.95) node[above, align=center] {\textit{\scriptsize Size of periodicity cell}};
			\legend{$\varphi_1$, $\varphi_2$, $\varphi_3$}
			}{}
		\end{axis}
	\end{tikzpicture}
	\caption{Real part of the half-line solution computed using the quasi-1D approach, with $\mathbb{P}^1$ Lagrange finite elements and $h = 2 \times 10^{-3}$, and for different values of $\varphi$.\label{fig:half_line_solution_phi}}
\end{figure}

\vspace{-0.5\baselineskip}
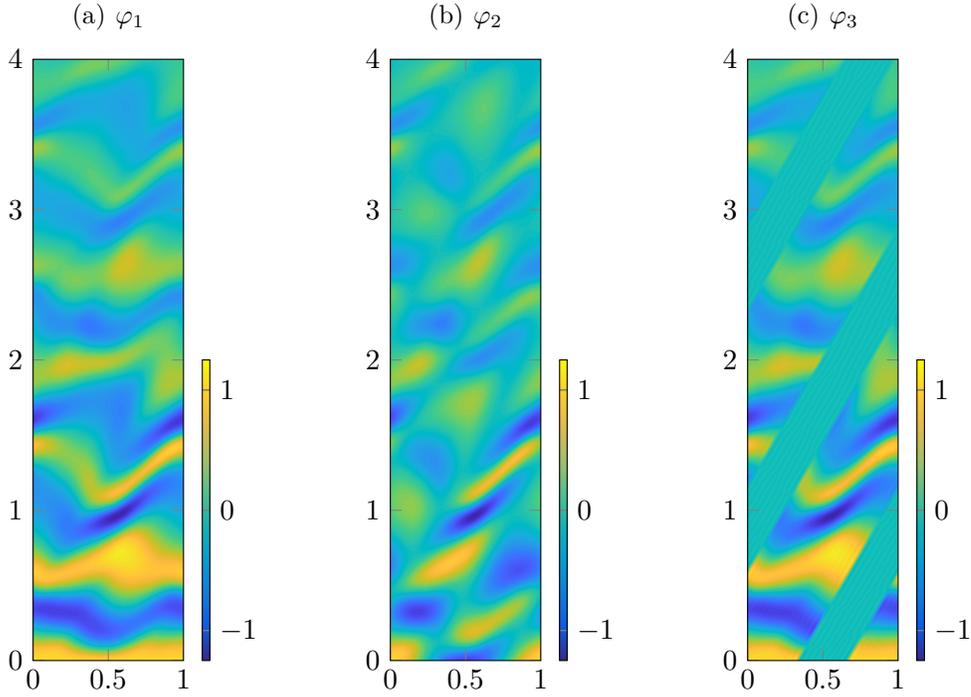
\begin{figure}[H]
	\centering
	\begin{tikzpicture}
		\begin{groupplot}[
			group style={
			group name=half_guide_plot,
			group size=3 by 1,
			horizontal sep=2.75cm,
			},
			enlargelimits=false,
			axis on top,
			width=2cm, height=8cm,
			scale only axis,
			xtick = {0, 0.5, 1}, ytick = {0, 1, 2, 3, 4},
			ticklabel style = {scale=1},
			disabledatascaling,
			colorbar,
			point meta min=-1.25, point meta max=1.25,
			colorbar style={
			yticklabel style={
			/pgf/number format/precision=2,
			/pgf/number format/fixed,
			},
			ytick={-1, 0, 1},
			ticklabel style = {scale=1},
			major tick length=0.05*\pgfkeysvalueof{/pgfplots/parent axis width},
			scaled y ticks=false,
			at = {(2.25cm, 0)},
			anchor=south west,
			width  = 0.05*\pgfkeysvalueof{/pgfplots/parent axis  width},
			height = 0.50*\pgfkeysvalueof{/pgfplots/parent axis height},
			},
			]
			\nextgroupplot[]
			\addplot graphics [xmin=0, xmax=1, ymin=0, ymax=4]{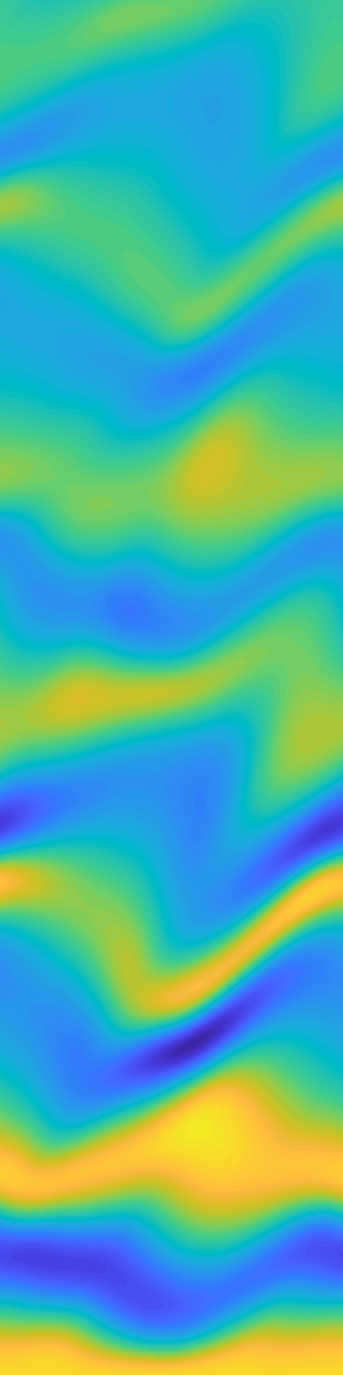};
			\nextgroupplot[]
			\addplot graphics [xmin=0, xmax=1, ymin=0, ymax=4]{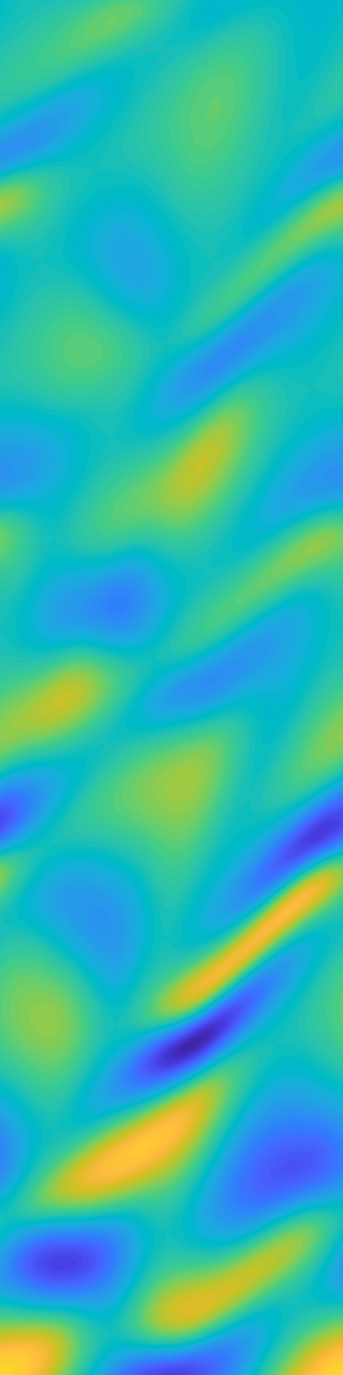};
			\nextgroupplot[]
			\addplot graphics [xmin=0, xmax=1, ymin=0, ymax=4]{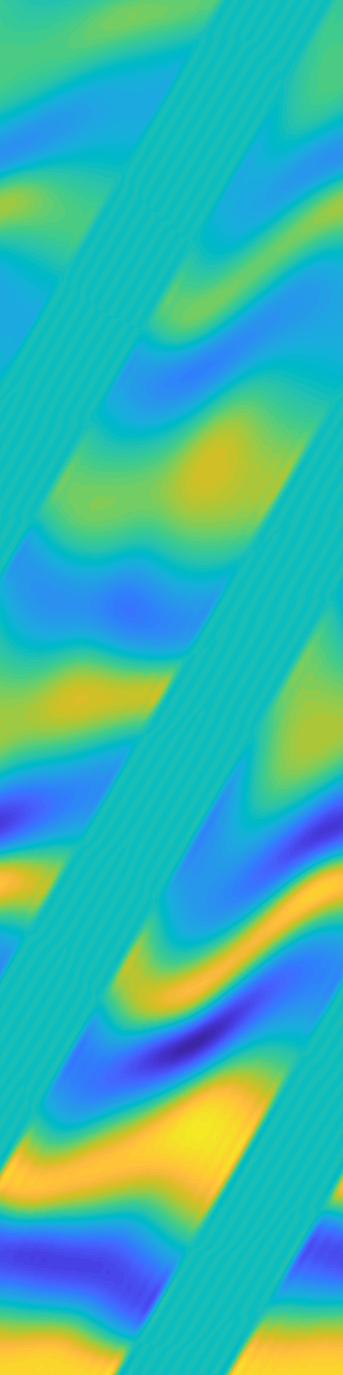};
		\end{groupplot}
		\node [text width=3cm, above=6pt] at (half_guide_plot c1r1.north) {\subcaption{$\varphi_1$\label{fig:half_guide_solution:a}}};
		\node [text width=3cm, above=6pt] at (half_guide_plot c2r1.north) {\subcaption{$\varphi_2$\label{fig:half_guide_solution:b}}};
		\node [text width=3cm, above=6pt] at (half_guide_plot c3r1.north) {\subcaption{$\varphi_3$\label{fig:half_guide_solution:c}}};
	\end{tikzpicture}
	\caption{Real part of the half-guide solution computed using the quasi-1D approach, with $\mathbb{P}^1$ Lagrange finite elements and $h = 2 \times 10^{-3}$, and for different values of $\varphi$.\label{fig:half_guide_solution}}
\end{figure}

 \noindent
 As expected, and as Figures \ref{fig:half_line_solution_phi} and \ref{fig:half_guide_solution:a}--\ref{fig:half_guide_solution:c} show, the aspect of half-guide solution changes extensively with respect to the boundary data, whereas the half-line solution looks invariant.

\subsubsection{The whole line problem}\label{sec:results:whole_line_problem}
The solutions $u^\pm_\cut$ of the half-line problems \eqref{eq:half_line_problems_0} allow one to compute the DtN coefficients $\lambda^\pm$, to solve \eqref{eq:interior_problem}, and then to compute the solution $u$ of Problem  \eqref{eq:whole_line_problem} using \eqref{eq:solution_of_whole_line_problem}. Recall that the coefficients $\mu$, $\rho$, and the source term $f$ are shown in Figure \ref{fig:coefficients_mu_rho}. The solution of \eqref{eq:whole_line_problem} is represented in Figure \ref{fig:whole_line_solution} for different values of $\omega$.
\begin{figure}[H]
	\begin{subfigure}[b]{0.99\textwidth}
		\centering
		\caption{$\omega = 8 + 0.25\,\icplx$}
		\begin{tikzpicture}
			\centering
			\begin{axis}[
				width=0.91\textwidth, height=2.5cm,
				ymin=-1.25, ymax=1.25,
				xmin=-6, xmax=6,
				xtick={-6, -4, -2, 0, 2, 4, 6},
				scale only axis,
				ticklabel style = {scale=1},
				enlargelimits=false,
				axis on top,
				legend style = {draw=none, fill=none, legend columns = 2, nodes={scale=1, transform shape}},
				]
				\ifthenelse{\boolean{afficherGraphes}}{
					\addplot [color=teal, smooth, thick] table [x={x}, y={u_real_omega1}] {\uLocPertQuasiD};
				}{}
			\end{axis}
		\end{tikzpicture}
	\end{subfigure}
\end{figure}
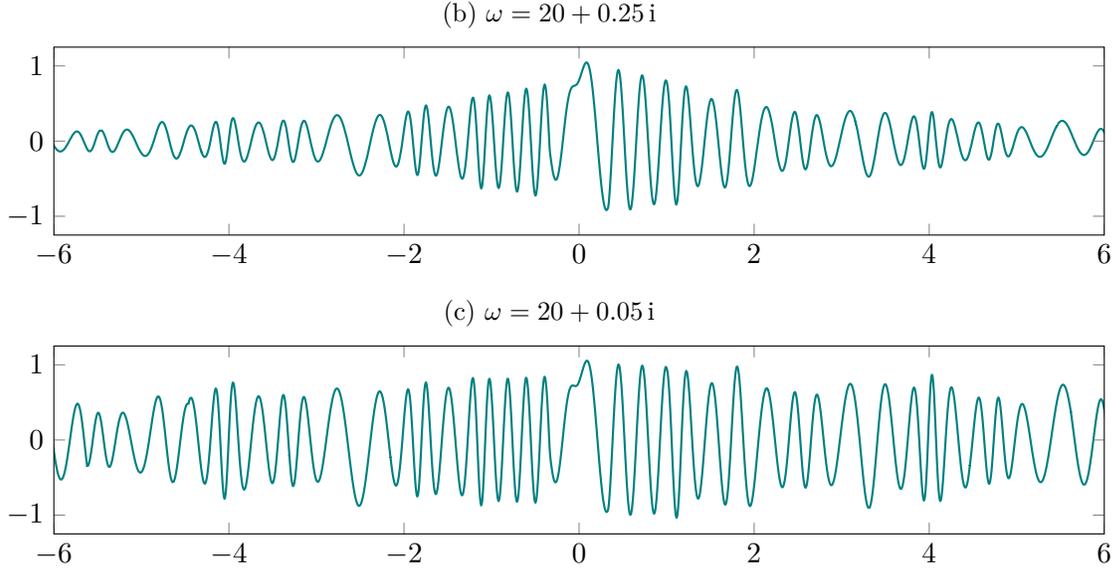
\begin{figure}[H]\ContinuedFloat
	%
	%
	\begin{subfigure}[b]{0.99\textwidth}
		\centering
		\caption{$\omega = 20 + 0.25\,\icplx$}
		\begin{tikzpicture}
			\centering
			\begin{axis}[
				width=0.91\textwidth, height=2.5cm,
				ymin=-1.25, ymax=1.25,
				xmin=-6, xmax=6,
				xtick={-6, -4, -2, 0, 2, 4, 6},
				scale only axis,
				ticklabel style = {scale=1},
				enlargelimits=false,
				axis on top,
				legend style = {draw=none, fill=none, legend columns = 2, nodes={scale=1, transform shape}},
				]
				\ifthenelse{\boolean{afficherGraphes}}{
					\addplot [color=teal, smooth, thick] table [x={x}, y={u_real_omega2}] {\uLocPertQuasiD};
				}{}
			\end{axis}
		\end{tikzpicture}
	\end{subfigure}
	\vfill
	\begin{subfigure}[b]{0.99\textwidth}
		\centering
		\caption{$\omega = 20 + 0.05\,\icplx$}
		\begin{tikzpicture}
			\centering
			\begin{axis}[
				width=0.91\textwidth, height=2.5cm,
				ymin=-1.25, ymax=1.25,
				xmin=-6, xmax=6,
				xtick={-6, -4, -2, 0, 2, 4, 6},
				scale only axis,
				ticklabel style = {scale=1},
				enlargelimits=false,
				axis on top,
				legend style = {draw=none, fill=none, legend columns = 2, nodes={scale=1, transform shape}},
				]
				\ifthenelse{\boolean{afficherGraphes}}{
					\addplot [color=teal, smooth, thick] table [x={x}, y={ureal}] {\uLocPertQuasiDsuite};
				}{}
			\end{axis}
		\end{tikzpicture}
	\end{subfigure}
	\caption{Real part of the solution of \eqref{eq:whole_line_problem} computed using the quasi-1D approach, with $\mathbb{P}^1$ Lagrange finite elements and $h = 2 \times 10^{-3}$, and for different values of $\omega$.\label{fig:whole_line_solution}}
\end{figure}

\subsubsection{About the dependence with respect to the absorption}\label{sec:results:absorption}
\begin{figure}[ht!]
	\centering
	\begin{tikzpicture}
		\begin{groupplot}%
			[
				group style={
					group name=myplot,
			    group size=3 by 1,
			    horizontal sep=0.85cm,
					ylabels at=edge left,
			    yticklabels at=edge left,
				},
				width=0.25\textwidth, height=4.5cm,
				scale only axis,
				ymin = 2.46594e-04, ymax=1.07008e+00,
				ylabel={Relative errors},
				xmode=log, ymode=log,
				xtick={32, 64, 128, 256, 512},
				xticklabels = {32, 64, 128, 256, 512},
				enlargelimits=false,
				axis on top,
				legend style = {legend pos = south west, fill=none, legend columns = 1},
			]
			\nextgroupplot[legend to name={CommonLegendImagOmega}, legend style={legend columns=2}]
			%
			\addlegendimage{mark=*, bleuDoux}
			\addlegendimage{mark=square*, rougeTerre}
			\addlegendentry{2D method}
			\addlegendentry{Quasi-1D method}
			\plotcurveandregression{\erreurs}{nDOFs}{quasiD_omega1}{mark=square*, rougeTerre}{dashed, rougeTerre}{0.5}{|-}{pos=0.25, left}{\slopeQuasiDImaga}
			\plotcurveandregression{\erreurs}{nDOFs}{multiD_omega1}{mark=*, bleuDoux}{dashed, bleuDoux}{0.25}{-|}{pos=0.75, right}{\slopeMultiDa}
			\nextgroupplot[]
			\plotcurveandregression{\erreurs}{nDOFs}{quasiD_omega5}{mark=square*, rougeTerre}{dashed, rougeTerre}{0.5}{|-}{pos=0.25, left}{\slopeQuasiDImagb}
			\plotcurveandregression{\erreurs}{nDOFs}{multiD_omega5}{mark=*, bleuDoux}{dashed, bleuDoux}{0.25}{-|}{pos=0.75, right}{\slopeMultiDb}
			\nextgroupplot[]
			\addplot [mark=*, bleuDoux] table [x={nDOFs}, y={multiD_omega8}] {\erreurs};
			\plotcurveandregression{\erreurs}{nDOFs}{quasiD_omega8}{mark=square*, rougeTerre}{dashed, rougeTerre}{0.75}{|-}{pos=0.25, left}{\slopeQuasiDImagc}
			\plotcurveandregression{\erreurs}{nDOFs}{multiD_omega8}{mark=*, bleuDoux}{dashed, bleuDoux}{0.5}{-|}{pos=0.75, right}{\slopeMultiDc}
		\end{groupplot}

		\draw ($0.5*(myplot c1r1.south)+0.5*(myplot c3r1.south) - (0, 0.5)$) node[below] {Discretization parameter $1/h$};
		\draw ($0.5*(myplot c1r1.south)+0.5*(myplot c3r1.south) - (0, 1.1)$) node[below] {\ref{CommonLegendImagOmega}};
		\node [text width=3cm, above=6pt] at (myplot c1r1.north) {\subcaption{$\omega = 8 + 0.25\,\icplx$\label{fig:erreurs_imagomega:a}}};
		\node [text width=3cm, above=6pt] at (myplot c2r1.north) {\subcaption{$\omega = 8 + 0.01\,\icplx$\label{fig:erreurs_imagomega:b}}};
		\node [text width=3cm, above=6pt] at (myplot c3r1.north) {\subcaption{$\omega = 8 + 0.001\,\icplx$\label{fig:erreurs_imagomega:c}}};
	\end{tikzpicture}
	\caption{Relative error in $H^1$ norm of the half-line solution for different values of $\omega$.\label{fig:erreurs_imagomega}}
\end{figure}
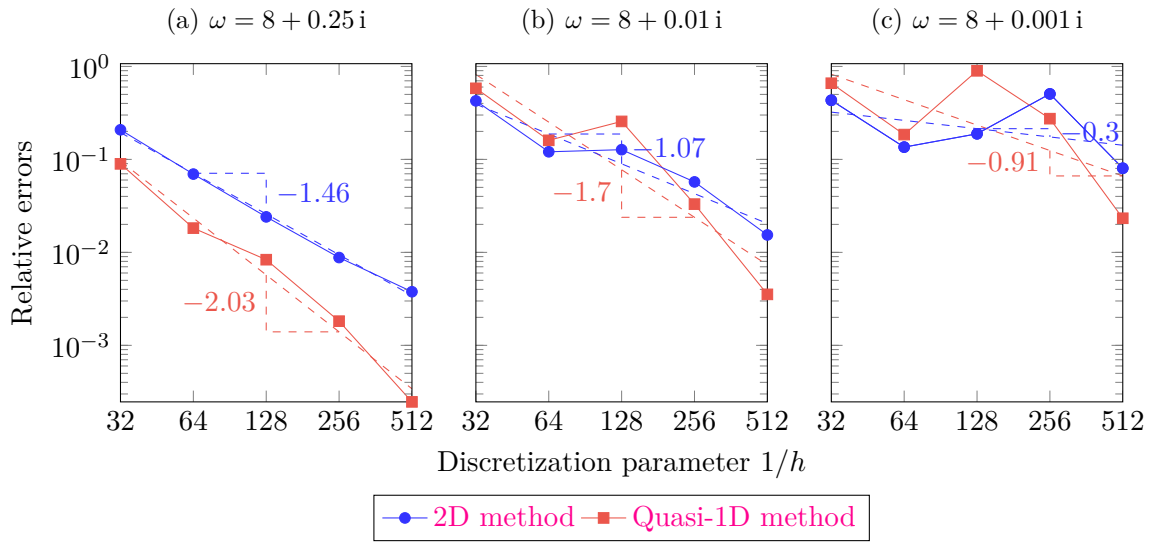

We come back to the numerical resolution of Problem \eqref{eq:half_line_problem}, and we study the convergence of the $2$D and quasi-1D methods depending on the absorption, especially when it tends to $0$. As in Section \ref{sec:results:half_line_guide}, the solutions are computed with Lagrange finite elements of degree $1$. The relative error $\varepsilon(u^+_\cut)$ defined \eqref{eq:relative_error} is represented in Figure \ref{fig:erreurs_imagomega} for both the $2$D and the quasi-1D method, and for different values of $\Imag \omega$.

\vspace{1\baselineskip}\noindent
As Figure \ref{fig:erreurs_imagomega} shows, the error deteriorates with $\Imag \omega$. It would mean that the numerical method becomes less efficient as the absorption decreases. This issue is closely related to the well-posedness of the local cell problems with Dirichlet boundary conditions when $\Imag \omega = 0$. In fact, for the elliptic Helmholtz equation, it is known (see \cite[Section 3.2.1.1]{flissthese} for instance) that the local cell problems are well-posed except for a \emph{countable} set of frequencies which correspond to the eigenvalues of the associated differential operator. In our case however, as the differential operator has a non-elliptic principal part, it also has a continuous spectrum, and one can show that when $\mu_p$ and $\rho_p$ are non-constant, the local cell problems are well-posed \emph{only for frequencies in a bounded set} (that can even be empty). An alternative to avoid this problem is to use a Robin-to-Robin operator instead of the DtN operator, which would involve solving well-posed local cell problems with Robin boundary conditions, as it is done in \cite{fliss2010exact} for periodic media. This will be done in a forthcoming paper for quasiperiodic media.

\subsubsection{About the spectral approximation of the propagation operator}\label{sec:results:spectral_approximation_P}
As explained in Subsection \ref{sec:discrete_Riccati_equation}, the discrete propagation operator $\mathcal{P}_h$ is computed by means of its eigenpairs. %
{
In this section, the eigenvalues of $\mathcal{P}_h$ are compared with the spectrum of the exact propagation operator which, according to Proposition \ref{prop:properties_P}, is a circle of radius
\[\displaystyle
M_{\log}(p_\cut) = \exp \Big(\int_0^1 \log |p_\cut(s)|\; d s \Big), \quad \textnormal{with} \quad p_\cut(s) = u^+_{s-\cuti_1/\cuti_2, \cut}(1/\sin \cuti_2).
\]
To compute this radius, $u^+_{s, \cut}$ is approximated by the unique function $u^+_{s, \cut, L}$ that satisfies \eqref{eq:half_line_problems} on a truncated domain $(0, L)$, with $\smash{u^+_{s, \cut, L}(L) = 0}$. One can show similar estimates to \eqref{eq:estimate_reference_solution}, and if $L$ is chosen large enough (for instance, if $L$ satisfies \eqref{eq:criterion_length_truncated_domain}), then $\smash{u^+_{s, \cut, L}}$ can be used as a reference solution. In practice, $\smash{u^+_{s, \cut, L}}$ is computed for several $s$, and finally the integral that defines $M_{\log}(p_\cut)$ is evaluated using a rectangular quadrature rule.
}

\vspace{1\baselineskip} \noindent
The spectra of $\mathcal{P}_h$ and $\mathcal{P}$ are shown in Figure \ref{fig:spectrum_P} for $\omega = 8 + 0.25\, \icplx$, and for different values of the discretization parameter $h$ (with $h_\cut = h$ for the quasi-1D method). Figure \ref{fig:eigenvalues_in_band} represents the number $N_h$ of eigenvalues of $\mathcal{P}_h$ that are close by $5\%$ to $\sigma(\mathcal{P})$, namely
\begin{equation}
	N_h = \# \bigg\{\lambda_h \in \sigma(\mathcal{P}_h)\ \ \Big/\ \ \bigg| \frac{|\lambda_h| - M_{\log}(p_\cut)}{M_{\log}(p_\cut)} \bigg| \leq 5\%\bigg\}.
\end{equation}
In Figure \ref{fig:eigenvalues_in_band}, one sees that $N_h$ increases with $1/h$, which means that more and more eigenvalues of $\mathcal{P}_h$ are close to $\sigma(\mathcal{P})$ when $h$ decreases. In other words, a finer discretization leads to a better approximation of the spectrum. The number $N_h$ of such eigenvalues also seems to increase linearly with $1/h$ (up to a subsequence for the quasi-1D method). Finally, note that $N_h$ is higher with the quasi-1D method than with the 2D method.

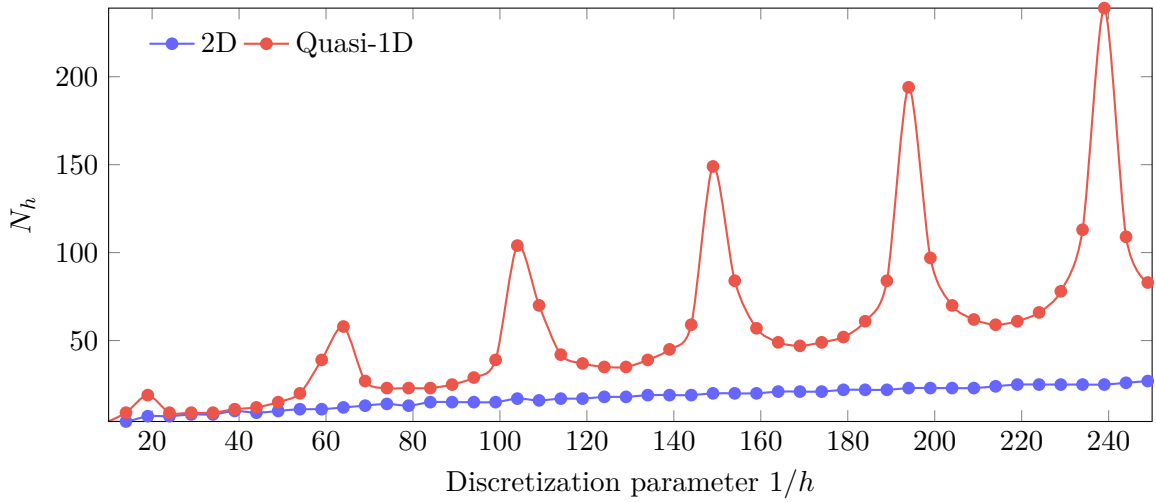
\begin{figure}
	\begin{tikzpicture}
		\begin{axis}[
			width=0.895\textwidth,
			height=5.5cm,
			xmin=10, xmax=250,
			xlabel={Discretization parameter $1/h$},
			ylabel={$N_h$},
			scale only axis,
			ticklabel style = {scale=1},
			enlargelimits=false,
			axis on top,
			legend style = {legend pos = north west, draw=none, fill=none, legend columns = 2, nodes={scale=1, transform shape}},
			]
			\addplot [color=bleuDoux!75, smooth, thick, mark=*] table [x={nDOFs}, y={numeigs_multiD}] {\numeigsP};
			\addplot [color=rougeTerre, smooth, thick, mark=*] table [x={nDOFs}, y={numeigs_quasiD}] {\numeigsP};
			\legend{$2$D, Quasi-$1$D}
		\end{axis}
	\end{tikzpicture}
	\caption{Number of eigenvalues of $\mathcal{P}_h$ that are close by $5\%$ to $\sigma(\mathcal{P})$ with respect to $h$.\label{fig:eigenvalues_in_band}}
\end{figure}
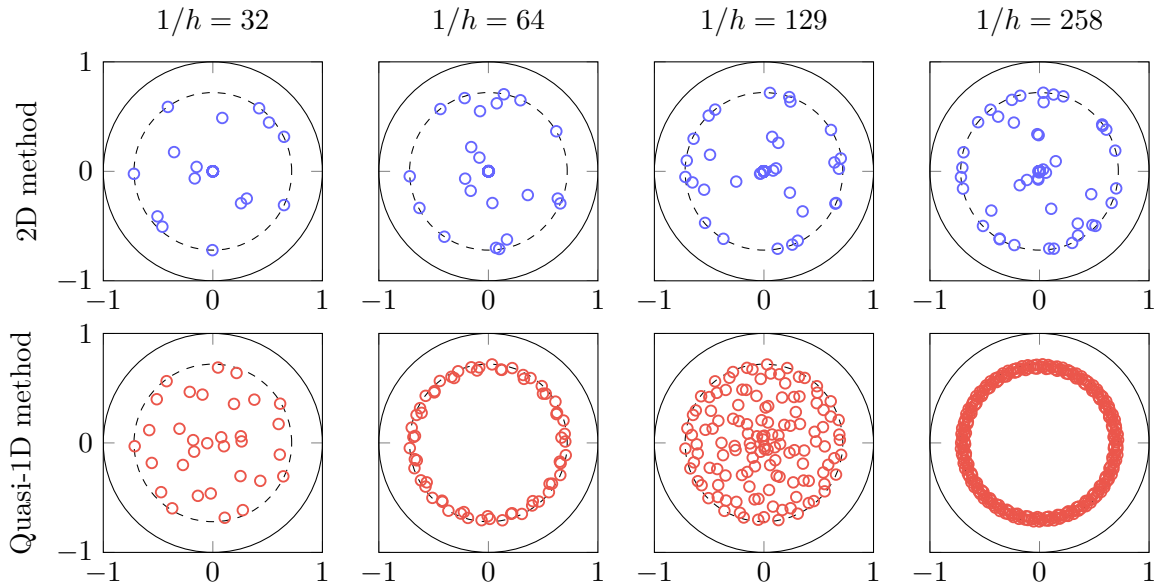
\begin{figure}
	\begin{tikzpicture}
		\begin{groupplot}[
			group style={
		    group size=4 by 2,
		    horizontal sep=0.75cm,
				vertical sep=0.7cm,
				ylabels at=edge left,
		    yticklabels at=edge left,
			},
			width=0.29\textwidth, height=0.29\textwidth,
			xmin=-1, xmax=1, ymin=-1, ymax=1,
			xtick={-1, 0, 1}, ytick={-1, 0, 1},
			enlargelimits=false,
			axis equal,
			]
			\nextgroupplot[ylabel={$2$D method}, title={$1/h = 32$}]
			\draw (axis cs:0, 0) circle [radius=1];
			\ifthenelse{\boolean{afficherGraphes}}{%
			\addplot[only marks, bleuDoux!75, line width=0.25mm, mark=o, unbounded coords=jump] table[x={real_eig_multi32}, y={imag_eig_multi32}] {\eigsP};
			}{}
			\draw [dashed] (axis cs:0, 0) circle [radius=0.719461];
			\nextgroupplot[title={$1/h = 64$}]
			\draw (axis cs:0, 0) circle [radius=1];
			\ifthenelse{\boolean{afficherGraphes}}{%
			\addplot[only marks, bleuDoux!75, line width=0.25mm, mark=o, unbounded coords=jump] table[x={real_eig_multi64}, y={imag_eig_multi64}] {\eigsP};
			}{}
			\draw [dashed] (axis cs:0, 0) circle [radius=0.719461];
			\nextgroupplot[title={$1/h = 129$}]
			\draw (axis cs:0, 0) circle [radius=1];
			\ifthenelse{\boolean{afficherGraphes}}{%
			\addplot[only marks, bleuDoux!75, line width=0.25mm, mark=o, unbounded coords=jump] table[x={real_eig_multi129}, y={imag_eig_multi129}] {\eigsP};
			}{}
			\draw [dashed] (axis cs:0, 0) circle [radius=0.719461];
			\nextgroupplot[title={$1/h = 258$}]
			\draw (axis cs:0, 0) circle [radius=1];
			\ifthenelse{\boolean{afficherGraphes}}{%
			\addplot[only marks, bleuDoux!75, line width=0.25mm, mark=o, unbounded coords=jump] table[x={real_eig_multi258}, y={imag_eig_multi258}] {\eigsP};
			}{}
			\draw [dashed] (axis cs:0, 0) circle [radius=0.719461];
			\nextgroupplot[ylabel={Quasi-$1$D method}]
			\draw (axis cs:0, 0) circle [radius=1];
			\ifthenelse{\boolean{afficherGraphes}}{%
			\addplot[only marks, rougeTerre, line width=0.25mm, mark=o, unbounded coords=jump] table[x={real_eig_quasi32}, y={imag_eig_quasi32}] {\eigsP};
			}{}
			\draw [dashed] (axis cs:0, 0) circle [radius=0.719461];
			\nextgroupplot[]
			\draw (axis cs:0, 0) circle [radius=1];
			\ifthenelse{\boolean{afficherGraphes}}{%
			\addplot[only marks, rougeTerre, line width=0.25mm, mark=o, unbounded coords=jump] table[x={real_eig_quasi64}, y={imag_eig_quasi64}] {\eigsP};
			}{}
			\draw [dashed] (axis cs:0, 0) circle [radius=0.719461];
			\nextgroupplot[]
			\draw (axis cs:0, 0) circle [radius=1];
			\ifthenelse{\boolean{afficherGraphes}}{%
			\addplot[only marks, rougeTerre, line width=0.25mm, mark=o, unbounded coords=jump] table[x={real_eig_quasi129}, y={imag_eig_quasi129}] {\eigsP};
			}{}
			\draw [dashed] (axis cs:0, 0) circle [radius=0.719461];
			\nextgroupplot[]
			\draw (axis cs:0, 0) circle [radius=1];
			\ifthenelse{\boolean{afficherGraphes}}{%
			\addplot[only marks, rougeTerre, line width=0.25mm, mark=o, unbounded coords=jump] table[x={real_eig_quasi258}, y={imag_eig_quasi258}] {\eigsP};
			}{}
			\draw [dashed] (axis cs:0, 0) circle [radius=0.719461];
		\end{groupplot}
	\end{tikzpicture}
	\caption{Eigenvalues of the discrete propagation operator (circle-shaped markers) compared to the spectrum of the exact propagation operator (circle in dashed line) for $\omega = 8 + 0.25\, \icplx$, and for different values of the discretization parameter.\label{fig:spectrum_P}}
\end{figure}
%
%

\vspace{1\baselineskip} \noindent
As Figure \ref{fig:spectrum_P} shows, the eigenvalues of $\mathcal{P}_h$ are all included in the disk of radius $\rho(\mathcal{P})$, but one observes some spectral pollution. This is a classical phenomenon when one approximates the spectrum of an operator which is neither compact nor self-adjoint. What is striking however, is that the pollution behaviours are very different depending on the method used.

On one hand, the eigenvalues obtained with the 2D approach tend to accumulate to $0$. A likely explanation for this phenomenon is that solving the local cell problems on 2D meshes does not take their directional structure into account. Since the location of the eigenvalues of $\mathcal{P}_h$ is similar to the one obtained in the elliptic case, for which $\mathcal{P}$ is compact (see \cite[Theorem 3.1]{jolyLiFliss}), we believe the 2D method somehow regularizes the half-guide problem \eqref{eq:half_guide_problem} by introducing an elliptic (discrete) approximation of the corresponding differential operator.

On the other hand, with the quasi-1D approach, the spectrum of $\mathcal{P}_h$ “oscillates” as the discretization parameter $h$ tends to $0$. This phenomenon has to do with the particular nature of $\mathcal{P}$ which is a weighted translation operator. We strongly suspect that one can extract a subsequence $(\mathcal{P}_{h'})$ whose spectrum converges towards $\sigma(\mathcal{P})$ in a sense to be defined precisely, as it is suggested by the peaks in Figure \ref{fig:eigenvalues_in_band}. The investigation of this assumption as well as the construction of such a subsequence are subject to ongoing works.

\vspace{1\baselineskip} \noindent
\textcolor{surligneur}{With both approaches, it has been observed numerically that the eigenfunctions associated to the spurious eigenvalues were highly oscillating functions that were badly approximated by the discretization, whereas the components of the half-guide solution with respect to these eigenfunctions are very small. This might explain why the spectral pollution does not have a visible influence on the approximation of the half-guide and the half-line solutions, as the errors in Figure \ref{fig:erreurs} seem to suggest.}

\section{Perspectives and ongoing works}
A numerical method has been proposed to solve Helmholtz equation in $1$D unbounded quasiperiodic media. Using the presence of absorption, we justified that this equation could be lifted onto a higher-dimensional problem which, in turn, can be solved using a Dirichlet-to-Neumann approach. For the discretization, we presented a multi-dimensional method, as well as a so-called quasi one-dimensional method. As shown by numerical simulations, both methods provide a suitable approximation of the solution as long as there is absorption. However, the quasi-1D method proved to be more efficient than the 2D method, as it takes the anisotropy of the problems involved into account.

\vspace{1\baselineskip}\noindent
The method presented opens up numerous perspectives, and raises multiple questions that are subject to ongoing works. For instance, it would be interesting to approximate efficiently the spectrum of the propagation operator, even though the spectral pollution seems to have no major impact on the efficiency of the overall method. Another key extension concerns the case where the absorption tends to $0$. This extension, which will be presented in a subsequent paper, involves replacing the DtN method by a Robin-to-Robin method as explained in Section \ref{sec:results:half_line_guide}, and finding a way to characterize the propagation operator which is no longer uniquely defined.

\vspace{1\baselineskip}\noindent
Finally, an approach which is similar to the one presented in this paper can be used to study the propagation of waves in presence of a $2$D periodic half-space when the interface does not lie in any direction of periodicity, or in presence of two $2$D periodic half-spaces with non-commensurable periods.


\printbibliography

\end{document}